\documentclass{siamart250211}
\usepackage{graphicx} 

\newcommand{\eps}{{\displaystyle \varepsilon}}
\newcommand{\bx}{\textbf{x}}
\newcommand{\A}{\mathcal{A}}
\newcommand{\bc}{\textbf{c}}
\newcommand{\by}{\textbf{y}}
\newcommand{\bz}{\textbf{z}}
\newcommand{\br}{\textbf{r}}
\newcommand{\bn}{\textbf{n}}
\newcommand{\bnu}{\bs{\nu}}
\newcommand{\bsub}{\begin{subequations}}
\newcommand{\esub}{\end{subequations}$\!$}
\newcommand{\ds}[0]{\displaystyle}
\newcommand{\bs}[0]{\boldsymbol}

\newcommand{\bigoh}{\mathcal{O}}
\newcommand{\G}{\mathcal{G}}

\newcommand{\Gint}{G^{\textrm{int}}}
\newcommand{\Gext}{G^{\textrm{ext}}}
\newcommand{\Rint}{R^{\textrm{int}}}
\newcommand{\Rext}{R^{\textrm{ext}}}

\newcommand{\al}[1]{\textcolor{black}{#1}}
\newcommand{\edit}[1]{\textcolor{black}{#1}}

\usepackage[sort&compress,numbers]{natbib}
\usepackage{graphicx,booktabs}
\usepackage{amsmath}
\usepackage{amssymb,subfigure}
\usepackage{comment}
\usepackage{amsfonts}
\usepackage{pgfplots}
\usepackage{color}
\usepackage[normalem]{ulem}
\definecolor{sectioncolor}{rgb}
{0.9, 0.3, 0.8}
\definecolor{hillencolor}{rgb}{0.50, 0.03, 0.46}
\definecolor{paintercolor}{rgb}{0.1, 0.7, 0.2}
\definecolor{mypink}{rgb}{0.9, 0.3, 0.6}
\definecolor{lightgray}{rgb}{0.1, 0.1, 0.1}

\usepackage{titlesec}
\pgfplotsset{compat=1.18}

\begin{document}

\title{High-order integral methods for the Neumann Green's function: applications to capture and signaling problems in two dimensions.}
\date{\today}
\author{Sanchita Chakraborty\thanks{Department of Applied and Computational Mathematics and Statistics, University of Notre Dame, Notre Dame, IN, 46656, USA. {\tt schakra2@nd.edu}}\and
Jeremy Hoskins\thanks{Department of Statistics, University of Chicago, USA and NSF-Simons National Institute for Theory and Mathematics in Biology, Chicago, IL {\tt jeremyhoskins@uchicago.edu}}
\and Alan E. Lindsay\thanks{Department of Applied and Computational Mathematics and Statistics, University of Notre Dame, Notre Dame, IN, 46656, USA. {\tt a.lindsay@nd.edu}}
}

\maketitle

\begin{abstract}
  We present a high order numerical method for the solution of the Neumann Green's function in two dimensions. For a general closed planar curve, our computational method resolves both the interior and exterior Green's functions with the source placed either in the bulk or on the surface -- yielding four distinct functions. Our method exactly represents the singular nature of the Green's function by decomposing the singular and regular components. In the case of the interior function, we exactly prescribe an integral constraint which is necessary to obtain a unique solution given the arbitrary constant solution associated with Neumann boundary conditions. Our implementation is based on a fast integral method for the regular part of the Green's function which allows for a rapid and high order discretization for general domains. We demonstrate the accuracy of our method for simple geometries such as disks and ellipses where closed form solutions are available. To exhibit the usefulness of these new routines, we demonstrate several applications to open problems in the capture of Brownian particles, specifically, how the small traps or boundary windows should be configured to maximize the capture rate of Brownian particles. 
  \end{abstract}

\label{firstpage}



\pagestyle{myheadings}
\markboth{S.~Chakraborty, J.~Hoskins, A.E.~Lindsay}{ Integral methods for planar Green's functions}

\section{Introduction}
\al{This paper concerns the high-order numerical solution of four distinct Neumann Green's functions in general two dimensional geometries. The Neumann Green's function (referred to as the \emph{Neumann function} in classical potential theory) is a central quantity in the study of numerous phenomena, including fluids, electrostatics, optics, heat conduction and Brownian motion \cite{CourantHilbert1953,Garabedian1964,BergmanSchiffer1953}.}

\begin{figure}
    \centering
    \includegraphics[width=0.95\textwidth]{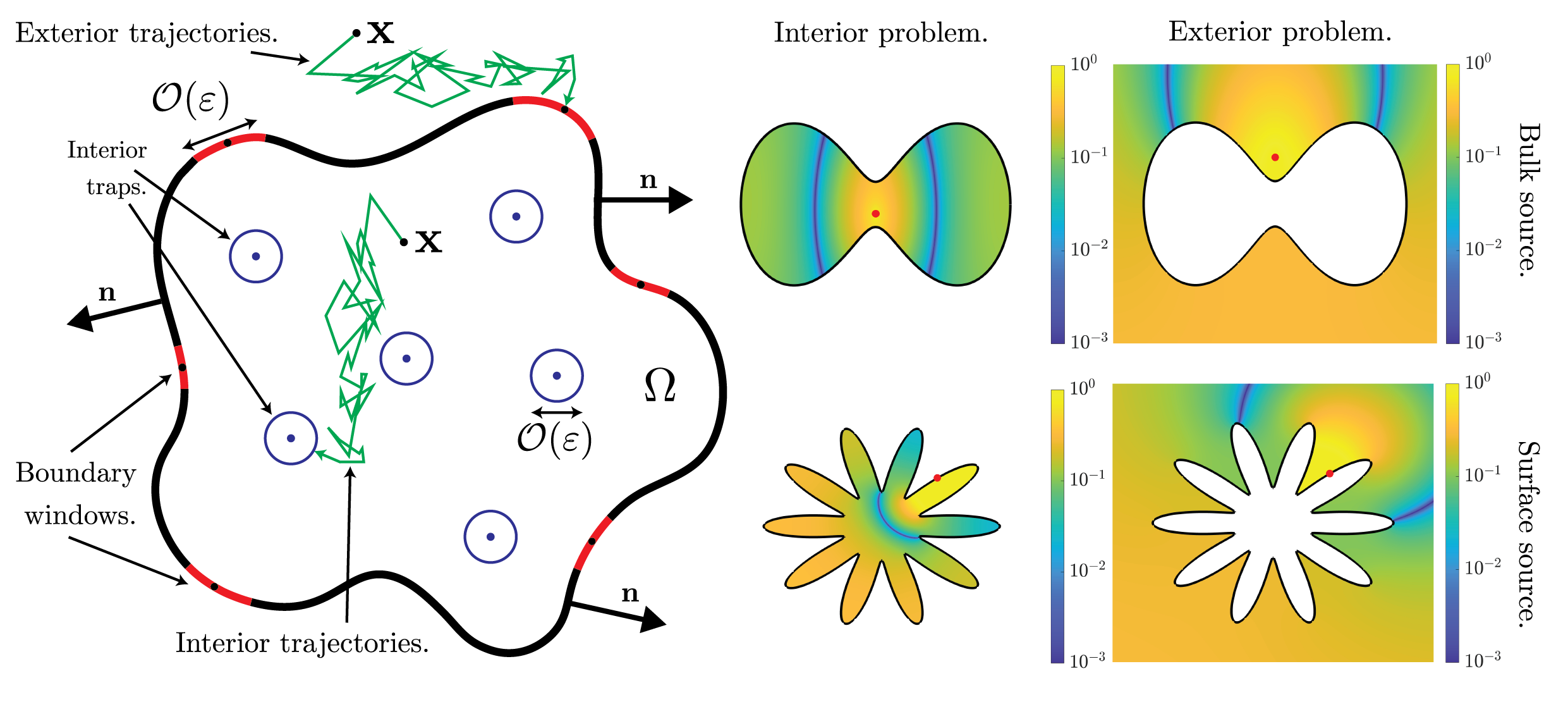}
    \caption{We consider Brownian motion both interior and exterior to a planar geometry $\Omega$. For the interior diffusion, a long standing problem is to find configurations of interior traps or locations of boundary windows that minimize the capture times. For the exterior problem, we seek the capture rate at boundary windows on $\partial\Omega$.  To address these questions, we develop numerical solutions of the interior $\Gint$ and exterior $\Gext$ Neumann Green's function for both surface sources (on $\partial\Omega$) and bulk (off $\partial\Omega$) sources (red dots). In total, we provide high accuracy methods for evaluating four Green's functions as seen in the above plots of $|\Gint_{s,b}|$ and $|\Gext_{s,b}|$. \label{fig:intro}}
\end{figure}

The setting of this study is a closed region $\Omega\in\mathbb{R}^2$ with an oriented boundary $\partial\Omega$ and outward facing normal vector $\textbf{n}$ (see Fig.~\ref{fig:intro}). The interior Green's function, $\Gint(\bx;\by)$ for $\bx\in\Omega$ and $\by\in\Omega$, solves the problem
\bsub\label{eqn_neumG}
\begin{gather}
 \label{eqn_neumG_a} \Delta \Gint = \frac{1}{|\Omega|}- \delta(\bx-\by), \quad \bx\in\Omega;\\[5pt]
\label{eqn_neumG_b} \qquad  \partial_{\bn} \Gint \equiv  \nabla \Gint\cdot \textbf{n} = 0,\quad \bx\in\partial\Omega;\\[5pt]
 \label{eqn_neumG_c} \int_{\Omega} \Gint(\bx;\by)d\bx = 0, \qquad \by\in\Omega.
\end{gather}
\esub
In equation \eqref{eqn_neumG_a}, the right hand side is modified with the term $1/|\Omega|$ therefore satisfying the solvability condition. The condition $\int_{\Omega} \Gint(\bx;\by)d\bx = 0$ in \eqref{eqn_neumG_c} removes the non-uniqueness associated with translation by a constant function and specifies that the solution has zero average. The solutions of \eqref{eqn_neumG} are singular as $\bx\to\by$ and related to the free-space Green's function
\begin{equation}\label{eqn:FS_Greens}
    \al{G_0(\bx;\by) =  -\ds\frac{1}{2\pi} \log|\bx-\by|.}
\end{equation}
When posed in $\bx\in\mathbb{R}^2\setminus\Omega$ with $\by\in\mathbb{R}^2\setminus\Omega$, the exterior Green's function $\Gext$ solves
\bsub\label{eqn_neumGExt}
\begin{gather}
 \label{eqn_neumGExt_a} \Delta \Gext = - \delta(\bx-\by), \quad \bx\in \mathbb{R}^2\setminus\Omega;\\[5pt]
\label{eqn_neumGExt_b} \qquad  \partial_{\bn} \Gext \equiv  \nabla \Gext\cdot \bn = 0,\quad \bx\in\partial\Omega.
\end{gather}
\esub
\al{The interior and exterior functions are further delineated by the position of the source. The Dirichlet trace of the Neumann function is called a \emph{surface} Green's function. When $\by\in\partial\Omega$, the surface Green's functions is henceforth labeled $G^{w}_s$ for $w\in\{\text{int},\text{ext}\}$. In the case where $\by\notin\partial\Omega$, the solution is referred to as the \emph{bulk} Green's function and is henceforth labeled $G^{w}_b$ for $w\in\{\text{int},\text{ext}\}$.
For $w\in\{\text{int},\text{ext}\}$, these four Green's functions have the singular-regular decompositions
\bsub\label{eqn:SingReg}
\begin{align}
    G^{w}_s(\bx;\by) &= 2 G_0(\bx;\by) + R^{w}_s(\bx;\by), \qquad \by\in\partial\Omega;\\[5pt]
    G^{w}_b(\bx;\by) &=  \phantom{2}G_0(\bx;\by) + R^{w}_b(\bx;\by),\qquad \by\notin\partial\Omega.
\end{align}
\esub
}
\al{The factor of two in the surface case reflects the fact that a source at $\by \in \partial\Omega$ subtends a half-plane angle of $\pi$ rather than the full $2\pi$ angle subtended by an interior source, an effect familiar from the jump relations of layer potentials \cite{kress}.} The functions $R^w_{s,b}(\bx;\by)$ for $w\in\{\text{int},\text{ext}\}$ are the \emph{regular parts} of the Green's function and bounded as $\bx\to\by$.  

The aim of this paper is to present a high-order boundary integral solution of these four functions. Our method exactly reproduces the singular-regular decompositions \eqref{eqn:SingReg} and, in the case of the interior Green's function, exactly satisfies the zero mean condition \eqref{eqn_neumG_c}. We will demonstrate that our method approximates $R^{w}_{s,b}(\bx;\by)$ (and hence $G^{w}_{s,b}(\bx;\by)$) for $w\in\{\text{int},\text{ext}\}$ to relative errors of around $\bigoh(10^{-12})$\footnote{Code available at \url{https://github.com/alanlindsay/FastGreensFunctions} }.

The primary motivation for this work is several problems in the theory of planar Brownian motion. For the interior problem, a long standing optimization problem asks: how should a collection of small internal traps or small boundary windows be arranged in order to minimize the expected capture time of a diffusing particle? The mean first passage time (MFPT) of a diffusing particle to a small target is a fundamental quantity in stochastic processes, with wide-ranging applications in cellular biology, physics, and material science \cite{BressloffSchumm2022,Miles2020,Newby2013,Camley2025,Simpson_2021}. Problems of this type — often referred to as narrow capture (NCP) or narrow escape (NEP) problems \cite{Schuss2014} — arise naturally in models of intracellular signaling \cite{Bressloff2024,BERNOFF2023}, receptor–ligand interactions \cite{CHOU2014,GIBBONS2010}, and diffusion-limited chemical reactions \cite{Grebenkov2018Nature}. 

The MFPT $T(\bx)$ gives the average time taken for a particle of diffusivity $D$, initialized at $\bx\in\Omega$, to reach an absorbing interior trap. The MFPT $T(\bx)$ satisfies \cite{REDNER2001,CHOU2014} the mixed boundary value problem
\bsub\label{eq:MFPT_main}
\begin{gather}
\label{eq:MFPT_main_a}    D\Delta T +1 = 0, \qquad \bx\in\Omega;\qquad \partial_{\bn} T = 0, \qquad \bx\in\partial\Omega;\\[5pt]
\label{eq:MFPT_main_b}    \eps\partial_{\bn}T + \kappa_jT = 0, \qquad \bx\in\partial\Omega_{\eps_j}\qquad j = 1,\ldots,N.
\end{gather}
\esub
In the context of the NCP, the set $\partial\Omega_{\eps} = \cup_{j=1}^N \Omega_{\eps_j}$ is a collection of $N$ small interior traps of extent $\mathcal{O}(\eps)$ centered at $\bx_j$, such that $\Omega_{\eps_j}\to \bx_j$ as $\eps\to0$ for $k=1,\ldots,N$. The boundary condition \eqref{eq:MFPT_main_b} encodes a range of reactivities $0<\kappa_j<\infty$ for the absorber. When $\kappa_j\to\infty$, the trap becomes absorbing so that any particle reaching $\partial\Omega_{\eps_j}$ reacts rapidly. Conversely for $\kappa_j\ll1$, the trap is weakly reactive so that particles reaching $\partial\Omega_{\eps_j}$ are mostly reflected \cite{BressloffSchumm2022,PCB2025,CengizLawley2024,GrebenkovWard2025a}. In the case of the NEP, the capture set is a collection of absorbing or partially absorbing boundary windows (see Fig.~\ref{fig:intro}).

A specific quantity that determines the overall capture rate for a domain $\Omega$ with traps centered at $\{\bx_1,\ldots,\bx_N\}$ is the global mean first passage time (GMFPT). This is a scalar quantity determined by
\begin{equation}\label{eqn:GMFPT}
    \tau = \tau(\bx_1,\ldots,\bx_N) = \frac{1}{|\Omega|}\int_{\Omega}T(\bx)d\bx,
\end{equation}
corresponding to an MFPT averaged over a uniform distribution of starting locations. 

\al{Direct numerical solution of \eqref{eq:MFPT_main} to obtain the GMFPT $\tau$ would require, for each candidate trap configuration $\{\bx_1,\ldots,\bx_N\}$, a fresh boundary integral equation on $\partial\Omega \cup \partial\Omega_{\varepsilon_1} \cup \cdots \cup \partial\Omega_{\varepsilon_N}$ in which the small trap boundaries must be geometrically resolved. As the trap radius $\varepsilon \to 0$, the panel counts needed to resolve features of size $\mathcal{O}(\varepsilon)$ grow correspondingly, and the system must be reassembled and resolved for every new configuration. Moreover, optimization requires multiple evaluations of $\tau$ and its derivatives which can be very expensive to obtain for high accuracy. Our method is able to approximate the Green's function and its derivatives to similar accuracy, thus avoiding the need for differencing.}

Matched asymptotic analysis performed in the limit as $\eps\to0$ is a powerful methodology that allows for the systematic reduction of the continuous problem \eqref{eq:MFPT_main} to a discrete one and a decoupling of the trap locations with the global geometry. The method constructs solutions of \eqref{eq:MFPT_main} by considering \lq\lq local\rq\rq\ problems in the vicinity of each trap where $|\bx-\bx_j| = \mathcal{O}(\eps)$ and a global or \lq\lq outer\rq\rq\ problem valid away from the traps where all the local behaviors are integrated. The details of this process vary in different dimensions and depending on whether the absorbers are centered interior to the domain ($\bx_j\in\Omega$) or on its boundary ($\bx_j\in\partial\Omega$). Across these geometric settings, a commonality is that the influence between traps is encoded in a discrete energy of form
\begin{equation}\label{eq:GreensMat}
    q(\bx_1,\ldots,\bx_N) = \bc^T \mathcal{G}\bc, \qquad \mathcal{G}_{i,j} = \left\{ \begin{array}{rcl}
    \Rint(\bx_i;\bx_i), & i = j;\\[4pt]
    \Gint(\bx_i;\bx_j), & i\neq j,
    \end{array} \right.
\end{equation}
where $\bc\in\mathbb{R}^N$ is a fixed vector that encodes variations in trap shapes and size. \al{The matrix $\mathcal{G}$ defined by \eqref{eq:GreensMat} is symmetric due to the reciprocity property of the Green's function $\Gint(\bx_i;\bx_j) = \Gint(\bx_j;\bx_i)$ for $i\neq j$. This  property follows from the self-adjointness of $-\Delta$ on $\Omega$ with Neumann boundary conditions and ensures that the quadratic form \eqref{eq:GreensMat} is well-defined.} Hence in the limit of small traps, the crucial information regarding the optimal configuration of traps is determined by the Green's function $\Gint(\bx;\by)$ and its regular part $\Rint(\bx;\by)$ which satisfy equation \eqref{eqn_neumG}. For identical traps, the optimization problem reduces to minimizing the functional
\begin{equation}\label{eq_Discrete}
    p(\bx_1,\ldots,\bx_N) = \sum_{i=1}^N \Big[ \Rint(\bx_i;\bx_i) + \sum_{ \substack{j=1 \\ j\neq i}}^N \Gint(\bx_j; \bx_i)\Big].
\end{equation}
In more general applications to singular solutions of nonlinear PDEs, it is long known that the critical points of \eqref{eq_Discrete} establish the limiting singular behavior for a wide class of nonlinear PDEs, see for example \cite{ChenWard2011,Davila2012,delPino2005,WWK09,TXKTW2017}.

The Green's function \eqref{eqn_neumG} encodes global information about the shape of the domain and hence closed form solutions for the Green's function are very valuable for determining optimizing configurations. Unfortunately, they are relatively scarce. In the two dimensional unit disk setting, the bulk regular part $\Rint_b(\bx;\by)$ was derived (cf.~\cite{KTW_2005}) to be
\begin{equation}\label{eqn:GreenRegDisk}
    \Rint_b(\bx;\by) = -\frac{1}{2\pi} \left[ \frac12 \log(1+ |\bx|^2 |\by|^2 - 2 \bx\cdot \by) - \frac{1}{2}(|\bx|^2 + |\by|^2) + \frac{3}{4} \right].
\end{equation}
This expression was utilized together with \eqref{eq_Discrete} to enable the identification of optimal annular ring configurations of traps with the possibility of a center trap for numbers $N = 1,\ldots, 25$ \cite{KTW_2005}.

In the case where diffusion is occurring exterior to a two dimensional geometry, the MFPT is not finite. In this setting, a more appropriate static quantity is the so-called \lq\lq splitting probabilities\rq\rq\ or steady-state fluxes. For a collection of boundary windows $\partial\Omega_{\eps_j}$, the splitting probabilities $\{\phi_k(\bx)\}_{k=1}^N$ give the probability that a particle initially at $\bx\in\mathbb{R}^2\setminus\Omega$ will be absorbed at the $k^{th}$ window. These quantities satisfy \cite{Venu2015,Lindsay2023,Miles2020}
\bsub\label{eq:intro_splitting}
\begin{gather}
\label{eq:intro_splitting_a} \Delta \phi_k = 0, \quad \bx\in\mathbb{R}^2\setminus\Omega; \qquad \phi_k \quad \text{bounded as } |\bx|\to\infty;\\[5pt]
\label{eq:intro_splitting_b} \phi_k = \delta_{jk}, \quad \bx\in\partial\Omega_{\eps_j}, \quad j = 1,\ldots,N; \qquad \partial_{\bn} \phi_k = 0, \quad \bx\in\partial{\Omega}\setminus\partial\Omega_{\eps}.
\end{gather}
\esub
The quantities $\{\phi_k(\bx)\}_{k=1}^N$ have been utilized in several recent studies \cite{Holcman2021,DOBRAMYSL2018,Miles2020,BERNOFF2023,BL2025,Lindsay2023} of directional sensing, the cellular process whereby cells detect the directionality of chemical cues from the noisy arrival of diffusing signaling molecules at membrane receptors. The problem \eqref{eq:intro_splitting} is posed exterior to the domain $\Omega$ with absorbing windows placed on $\partial\Omega$. As we discuss in Sec.~\ref{sec:splitting}, the asymptotic solution of this problem as $\eps\to0$ is derived in terms of the exterior surface Green's function $\Gext_{s}$ solving \eqref{eqn_neumGExt}.

For elliptical and rectangular geometries, rapidly convergent series based on re-summation of the Fourier series solution \cite{Sarafa2021,Gilbert_Cheviakov_2023} have recently been developed which provide the necessary accuracy to employ finite differencing for approximation of the gradient and Hessian of the objective function - an essential ingredient for reliable optimization. This new analysis identified GMFPT minimizing configurations and bifurcations of co-linear configurations \cite{Cheviakov2024}. In contrast to the exact solution  \eqref{eqn:GreenRegDisk} for the disk, we cannot easily use these rapidly convergent series to find exact bifurcation values or closed forms for optimizing arrangements. Hence our high accuracy numerical algorithm to approximate solutions to \eqref{eqn_neumG} and \eqref{eqn_neumGExt} has utility comparable to series solutions for the investigation of optimal trap arrangements, with the advantage of being applicable to general geometries.

In a recent survey article \cite{Cheviakov2024}, the authors mention the paucity of high fidelity solutions to \eqref{eqn_neumG} as a major impediment to further study of this central problem in stochastic processes, a bottleneck that has hindered the broader application of asymptotic–numerical methods for trap optimization in realistic biological and physical settings.

In the present work, we fill this gap by developing a boundary integral formulation that is able to provide accurate solutions to \eqref{eqn_neumG} and \eqref{eqn_neumGExt} for general two dimensional geometries $\Omega$ with both bulk and surface sources (cf.~Fig.~\ref{fig:intro}). Our method is based on two steps; 1) a separation of the solution into a prescribed singular component and a regular component that our method resolves for, and 2) an exact implementation of the average constraint $\int_{\Omega}\Gint(\bx;\by)d\bx = 0$ for the interior problem. Our approach combines a boundary integral representation for the regular part with fast multipole acceleration and high order quadrature to achieve both accuracy and computational efficiency. Moreover, the explicit nature of the boundary integral solution ansatz allows us to calculate the gradient and Hessian of the solution to similar precision, without any need for differencing.

We validate our approach against known analytical benchmarks (disk and ellipse geometries) and demonstrate its performance on deformed and non-convex domains. We demonstrate several examples of the optimization of trap configurations in general geometries. Beyond MFPT analysis, our framework provides a foundation for tackling related questions in eigenvalue optimization, homogenization theory, and biophysical modeling of diffusive search.

The organization of the paper is as follows. In Sec.~\ref{sec:MFPT} we briefly review results from existing MFPT theory for small traps in two dimensions and the discrete optimization problems which arise from them. In Sec.~\ref{sec:NumericalMethods} we describe our boundary integral approach to the solution of \eqref{eqn_neumG} and our approach to integrating with existing  routines to find minima of \eqref{eq_Discrete}. The accuracy of our method is validated against existing results for the disk and ellipse geometry and known optimal configurations in these settings. In Sec.~\ref{sec:Results} we present results on the optimal locations of interior traps and boundary windows for a variety of general geometries, including barbell, dumbbell, ellipses and a random domain. In Sec.~\ref{sec:Discussion} we conclude with some potential areas for future development and applications of these methods.

\section{MFPT theory for small absorbers in two spatial dimensions.}\label{sec:MFPT}

In this section, we recap upon results in the asymptotic analysis of \eqref{eq:MFPT_main} in the limit as $\eps\to0$ in two dimensions $d=2$. We consider the scenario where the traps are interior to the domain - the so called narrow capture problem.

\subsection{Two Dimensional analysis}

 For the case $d=2$, it has been determined \cite{Sarafa2021,Gilbert_Cheviakov_2023,Venu2015,GrebenkovWard2025a} that as $\eps\to0$ that the solution of \eqref{eq:MFPT_main} admits the expansion $T(\bx) = T_0(\bx;\bnu) + \mathcal{O}(\eps)$ where the problem $T_0$ solves
\bsub\label{eqn:T0_intro}
 \begin{gather}
\label{eqn:T0_intro_a}     D\Delta T_0 = -1, \quad \bx\in\Omega; \qquad \partial_{\bn} T_0 = 0, \quad \bx\in\partial \Omega;\\[5pt]
\label{eqn:T0_intro_b}     T_0 \sim S_j\nu_j \log |\bx-\bx_j| + S_j, \qquad \bx\to\bx_j, \qquad j = 1,\ldots,N.
 \end{gather}
\esub
Here $S_j$ are constants to be found by matching and  $\nu_j=-1/\log \eps d_j$. The constants $d_j$ for $j=1,\ldots,N$ are the logarithmic capacitance determined by the shape $\mathcal{A}_c$ of each trap by the problem
\bsub
\begin{gather}
    \Delta v_c = 0; \quad \bx\in\mathbb{R}^2\setminus \A_c; \qquad \partial_{\bn} v_c + \kappa_c v_c = 0, \quad \bx \in \partial\A_c\\[4pt]
    v_c = \log|\bx| - \log d_c, \qquad |\bx| \to \infty.
\end{gather}
\esub
A process for determining the functional relationship between the logarithmic capacitance and the reactivity parameter, $d_c = d_c(\kappa_c)$, was recently deployed in \cite{GrebenkovWard2025a}. In terms of the Green's function \eqref{eqn_neumG}, the solution of \eqref{eqn:T0_intro} is given by
\begin{equation}\label{eqn:T0_sol}
    T_0 = \al{-2\pi }\sum_{j=1}^NS_j \nu_j\Gint(\bx;\bx_j) + \tau_0.
\end{equation}
The value of $\tau_0$, known as the global MFPT, is determined as
\begin{equation}
    \tau_0 = \frac{1}{|\Omega|}\int_{\Omega}T_0(\bx)d\bx \al{\, ,}
\end{equation}
and gives a broad measure of capture statistics averaged uniformly over all starting locations. From matching the solution \eqref{eqn:T0_sol} to the local behavior \eqref{eqn:T0_intro_b}, a linear system for the unknowns is
\begin{equation}\label{eq:LAS}
    \big[ I + 2\pi \mathcal{G} \mathcal{V} \big]\textbf{S} = \tau_0 \textbf{e}, \qquad \sum_{j=1}^N\nu_j S_j = \frac{|\Omega|}{2\pi D}.
\end{equation}
\al{Here} $\mathcal{V}$ is a diagonal matrix with entries $\bnu=(\nu_1,\ldots,\nu_N)$ along the main diagonal and $I$ is the $N\times N$ identity matrix. In terms of the Green's matrix $\G$ defined in \eqref{eq:GreensMat}, the linear system \eqref{eq:LAS} can be solved for the strengths $\textbf{S} = (S_1,\ldots,S_N)$ and the GMFPT $\tau_0$. However, we can take the approximation one step further by applying the asymptotic inverse $[I+ 2\pi\mathcal{G} \mathcal{V}]^{-1}\sim [I- 2\pi\mathcal{G} \mathcal{V}]$ which \al{is valid for $\|\mathcal{G} \mathcal{V}\|\ll1$}. This yields the following solution
\begin{equation}\label{eqn:exp_approx}
    \tau_0 = \tilde\tau_0 + \mathcal{O}(\nu_M^2), \qquad \tilde\tau_0 = \frac{|\Omega|}{2\pi D\bar{\bnu} N}\left[ 1+ \frac{2\pi}{N\bar{\bnu}}\bnu^{T}\mathcal{G}\bnu  \right],
\end{equation}
where $\bar{\bnu} = N^{-1}\sum_{j=1}^N\nu_j$ and $\nu_M = \max_{1\leq j\leq N} \nu_j$. If we specialize to the case of identical traps, such that $\nu_j = \nu$ for $j = 1,\ldots,N$, then we have that
\begin{equation}\label{eqn:tau_0}
    \tilde\tau_0  = \frac{|\Omega|}{2\pi D\nu N}\left[ 1+ \frac{2\pi\nu}{N}p(\bx_1,\ldots,\bx_N)  \right] \qquad p(\bx_1,\ldots,\bx_N) = \sum_{i=1}^N\sum_{j=1}^N \mathcal{G}_{i,j}.
\end{equation}
\al{We note that the expansion \eqref{eqn:exp_approx} requires $2\pi \nu_M \| \mathcal{G} \| \ll 1$, while $\| \mathcal{G} \|$ grows roughly as $\mathcal{O}(N)$ for well-separated traps). In principle, this approximation can break down in the regime $N \gg \nu_M^{-1}$ even when the exact system \eqref{eq:LAS} remains well-conditioned. Empirically we observe negligible discrepancy up to $N = 25$. This is in agreement with previous optimization studies \cite{Cheviakov2024,Sarafa2021} that considered minimizers of both $\tau_0$ derived from the linear system \eqref{eq:LAS} and the explicit formulation $\tilde{\tau}_0$ derived in \eqref{eqn:tau_0} with remarkably good agreement observed between the two values.  In our results section, we focus on minimizers of $p(\bx_1,\ldots, \bx_N)$ which leads to a simplified  $\eps$ independent presentation. }

\subsection{The Neumann eigenvalue problem for the Laplacian}

A problem closely related to the narrow capture problem is the identification of trap configurations that maximize \cite{Bramburger2025,KTW_2005} the principal eigenvalue $\lambda_0$ of the problem
\bsub
\begin{gather}
    \Delta u + \lambda u = 0, \qquad \bx\in\Omega\setminus\Omega_{\eps};\qquad \partial_{\bn} u = 0, \quad \bx\in\partial\Omega;\\[5pt] \qquad  
    \eps\partial_{\bn} u + \kappa_j u = 0, \qquad \bx\in\partial\Omega_{\eps j}, \qquad j = 1,\ldots,N.
\end{gather}
\esub
In the limit as $\eps\to0$ for $N$ identical traps, the limiting behavior (cf.~\cite{Sarafa2021,KTW_2005}) of $\lambda_0$ is
\begin{equation}\label{eqn:}
    \lambda_0 \sim \frac{2 \pi N\nu}{|\Omega|} - \frac{4\pi^2 \nu^2}{|\Omega|} p(\bx_1,\ldots,\bx_N) + \bigoh(\nu^3),
\end{equation}
where the function $p(\bx_1,\ldots,\bx_N)$ is again given in \eqref{eq_Discrete}. 

\subsection{Splitting probabilities and directional sensing.}\label{sec:splitting}

In recent applications of asymptotic theory to chemical signaling, several authors have investigated how a cell can infer the source of a diffusive signal from localized receptor activity \cite{Holcman2021,DOBRAMYSL2018,Miles2020,Lindsay2023}. Specifically, let us define $\phi_k(\bx)$ to be the probability that a signaling molecule originally at $\bx\in\mathbb{R}^2\setminus\Omega$ first arrives at the $k^{th}$ surface receptor. We can consider the steady state flux into each receptor, also known as the \emph{splitting probabilities}. The splitting probabilities satisfy the exterior Laplace mixed value problem
\bsub\label{eq:SplittingProbs}
\begin{gather}
\label{eq:SplittingProbs_a} \Delta \phi_k = 0, \quad \bx\in\mathbb{R}^2\setminus\Omega; \qquad \phi_k \quad \text{bounded as } |\bx|\to\infty;\\[5pt]
\label{eq:SplittingProbs_b} \phi_k = \delta_{jk}, \quad \partial\Omega_{\eps_j}, \quad j = 1,\ldots,N; \qquad \partial_{\bn} \phi_k = 0, \quad \bx\in\partial{\Omega}\setminus\cup_{j=1}^N\partial\Omega_{\eps_j}\al{\,.}
\end{gather}
\esub
The process of developing an asymptotic solution to \eqref{eq:SplittingProbs} as $\eps\to0$ is similar and described in details in many works. As is commonly performed in other analysis, the key step is to replace each receptor by a local singularity condition. Specifically, we expand the solution of \eqref{eq:SplittingProbs} in a \lq\lq sum-of-all logs\rq\rq\ expansion $\phi_k(\bx) = \phi_k^{\ast}(\bx;{\bs \nu}) + o(1)$ where 
\bsub\label{eq:SplittingProbs_leading}
\begin{gather}
\label{eq:SplittingProbs_leading_a} \Delta \phi_k^{\ast} = 0, \quad \bx\in\mathbb{R}^2\setminus\Omega; \qquad \phi_k^{\ast} \quad \text{bounded as } |\bx|\to\infty;\\[5pt]
\label{eq:SplittingProbs_leading_b}
\partial_{\bn} \phi_k^{\ast} = 0, \quad \bx\in\partial{\Omega}\setminus\cup_{j=1}^N\{\bx_j\}\\[5pt]
\label{eq:SplittingProbs_leading_c} \phi_k^{\ast} \sim \delta_{jk} + S_{jk}\nu_j\log|\bx-\bx_j| + S_{jk}, \quad \mbox{as} \quad \bx\to\bx_{j}\, \quad j = 1,\ldots,N; \qquad 
\end{gather}
\esub
The solution of this problem is expressed as
\begin{equation}\label{eq:splitting_asy}
    \phi_k^{\ast}(\bx) = -\pi\sum_{j=1}^N S_{jk}\nu_j \Gext_s(\bx;\bx_j) + \bar{\phi}_k,
\end{equation}
where $\bar{\phi}_k = \lim_{|\bx|\to\infty}\phi_k^{\ast}$ is a constant arising form the homogeneous solution which describes the splitting probabilities averaged over all initial locations. The matching procedure as $\bx\to\bx_j$, together with the solvability condition $\sum_{j=1}^N \nu_jS_{jk} = 0 $ yields a system of equations for the unknowns $(S_{1k},\ldots,S_{Nk},\bar{\phi}_k)$. The system can be compactly represented as
\bsub\label{eq:splitting_linear}
\begin{equation}
\begin{bmatrix}
\mathcal{I} + \pi \, \mathcal{G}_s \, \mathcal{V}  & -\textbf{e}^T \\[5pt]
\bs{\nu}^T & 0
\end{bmatrix}
\begin{bmatrix}
\textbf{S}_k\\[5pt]
\bar{\phi}_k
\end{bmatrix}
= -
\begin{bmatrix}
\textbf{e}_k\\[5pt]
0
\end{bmatrix},
\qquad
\begin{array}{c}
\textbf{S}_k = (S_{1k},\ldots, S_{Nk}), \\[5pt] 
\bs{\nu} = (\nu_1,\ldots,\nu_N),
\end{array}
\end{equation}
\begin{gather}
\begin{array}{c}
\textbf{e} = (1,1,\ldots, 1), \\[5pt] \textbf{e}_k = (0,\ldots \underbrace{ 1 }_{k^{\textrm{th}}}  \ldots,  0 ),
\end{array}
\quad
[\mathcal{G}_s]_{i,j} = 
\left\{
\begin{array}{cl}
\Rext_s (\bx_i) & i = j\\[5pt]
\Gext_s (\bx_i;\bx_j) & i\neq j
\end{array}\right.
\quad
[\mathcal{V}]_{i,j} = 
\left\{
\begin{array}{cl}
\nu_i & i = j\\[5pt]
0 & i\neq j
\end{array}\right.
\label{eq:EqnEGV}
\end{gather}
\esub
and $\mathcal{I}$ is the $N\times N$ identity matrix. The key step is to now determine the Green's matrix $\mathcal{G}_s$ \al{which will be performed with our numerical method.}

\section{Numerical methods.}\label{sec:NumericalMethods}

In this section we describe and validate our numerical methods, including our boundary integral formulation for the solutions of \eqref{eqn_neumG} and \eqref{eqn_neumGExt} together with the optimization routines.

\subsection{Boundary integral methods.}\label{sec:BIEM}
 We first reduce the problem to a boundary integral equation, an approach which is commonplace in potential theory (cf.~\cite{kress,rachhgreengard, AMMARI201266,Fryklund2023}). Our presentation focuses on the interior Green's function \eqref{eqn_neumG} which, compared to the counterpart exterior function $\Gext$, has two additional features to be accounted for. First the source term $1/|\Omega|$ in \eqref{eqn_neumG_a} and second, the integral constraint \eqref{eqn_neumG_c}. To begin our process of reduction to an integral equation, we introduce for $\bx = (x_1,x_2)$ the particular solution
\begin{equation}\label{eqn:v}
    v(\bx) = \frac{x_1^2+x_2^2}{4}, \qquad \Delta v  = 1.
\end{equation}
Following this, we separate the solution of \eqref{eqn_neumG} into singular, particular and homogeneous components. \al{For a bulk source $\by\notin\partial\Omega$,} the homogeneous solution is defined as
\begin{equation}\label{eqn:wp}
    w(\bx;\by) = \Gint(\bx;\by) - \frac{1}{|\Omega|}v(\bx)-G_0(\bx;\by) .
\end{equation}
\al{In the case of a surface source $\by\in\partial\Omega$, the free space contribution in \eqref{eqn:wp} is $-2G_0(\bx;\by)$. In the following discussion, we specialize to the bulk solution.} The equation for the homogeneous component $w(\bx;\by)$ satisfies
\begin{align}\label{eqn:BEMHomogeneous}
    \begin{cases}
        \Delta w(\bx;\by) = 0, \quad &\bx \in \Omega,\\[5pt]
        \partial_{\bn} w(\bx;\by) = -\ds\frac{1}{|\Omega|}\partial_{\bn} v(\bx)-\partial_{\bn} G_0(\bx;\by), \quad & \bx \in \partial \Omega,\\[5pt]
        \ds\int_\Omega \Big[ w(\bx;\by)+\frac{1}{|\Omega|}v(\bx)+G_0(\bx;\by)\,{\rm }\Big]d\bx = 0.
    \end{cases}
\end{align}

The solution ansatz for the problem \eqref{eqn:BEMHomogeneous} is given by
\begin{equation}\label{eqn:Integral_pot}
    w(\bx;\by) =-\frac{1}{2\pi} \int_{\partial \Omega} \log|\bx-\bz|\,\sigma(\bz;\by)\,{\rm d}S(\bz)+\alpha,
\end{equation} 
for an unknown density $\sigma$ and constant $\alpha(\by),$ depending on $\by$ but independent of $\bx$. Taking the limit as $\bx\to\partial\Omega$, we obtain the following integral equation for the unknown density $\sigma$ 
\begin{align}\label{eq:int_eq}
    \frac{1}{2}\sigma(\bx;\by) -\frac{1}{2\pi}\int_{\partial \Omega} \frac{(\bx-\bz)\cdot \bn}{|\bx-\bz|^2}\sigma(\bz;\by)\,{\rm d}S(\bz) =  -\frac{1}{|\Omega|}\partial_{\bn} v(\bx)-\partial_{\bn} G_0(\bx;\by). 
\end{align}
Here $\bn$ denotes the \edit{outward pointing} normal derivative of $\partial\Omega$ at $\bx\in\partial\Omega$. We note that the left-hand side has a one-dimensional nullspace \cite{kress}. By the divergence theorem, the right-hand side integrates to zero which guarantees solvability. We are then free to impose the additional constraint $\int_{\partial \Omega} \sigma(\bz;\by)\,{\rm d}S(\bz) = \alpha(\by).$
We thus obtain
\begin{align}\label{eq:int_eq2a}
    \frac{1}{2}\sigma(\bx;\by) &-\frac{1}{2\pi}\int_{\partial \Omega} \frac{(\bx-\bz)\cdot \bn}{|\bx-\bz|^2}\sigma(\bz;\by)\,{\rm d}S(\bz) + \int_{\partial \Omega}\sigma(\bz;\by)\,{\rm d}S(\bz) -\alpha(\by) \nonumber\\
    &=  -\frac{1}{|\Omega|}\partial_{\bn} v(\bx)-\partial_{\bn} G_0(\bx,\by).
\end{align}
We refer the reader to \cite{sifuentes,kress} and the references therein for related treatments of this nullspace.

\subsection*{Implementation of the integral constraint.} Turning to the integral constraint on $\Gint,$ defined in \eqref{eqn_neumG_c}. We recall from \eqref{eqn:v} that $\Delta v(\bx) \equiv 1,$ and apply Green's second identity to obtain that
\begin{align}
\nonumber    0 &= \int_\Omega \Gint(\bx;\by)\,{\rm d}\bx= \int_{\Omega} (\Delta v(\bx))\, \Gint(\bx;\by)\,{\rm d}\bx,\\
    &= \frac{1}{|\Omega|}\int_\Omega v(\bx)\,{\rm d}\bx-v(\by)+ \int_{\partial \Omega}  \partial_{\bn}v (\bx) \, \Gint(\bx;\by)\,{\rm d}S(\bx)
\end{align}
In the above expression, the term $v(\by)$ arises from integration against either a bulk or surface Green's function. In the end, we wish for all calculations to be performed on $\partial \Omega$, hence we note from the divergence theorem that
$$\int_{\Omega} v(\bx)\,{\rm d}\bx = \frac{1}{4}\int_{\Omega} (x_1^2+x_2^2) d\bx = \frac{1}{12}\int_{\partial \Omega} \bn \cdot (x_1^3,x_2^3)\,{\rm d}S(\bx).$$

Let $\mu(\by)$ be defined by
\begin{align}
    \mu(\by)=\frac{1}{|\Omega|}\int_\Omega v(\bx)\,{\rm d}\bx-v(\by)+ \int_{\partial \Omega} \partial_{\bn} v (\bx)\,\left(\frac{1}{|\Omega|}v(\bx) + G_0(\bx;\by)\right)\,{\rm d}S(\bx).
\end{align}
Then, the integral constraint becomes $\int_{\partial \Omega} \partial_{\bn} v(\bx)\, w(\bx;\by)\,{\rm d}S(\bx) = -\mu(\by).$
Expanding, gives
$$-\frac{1}{2\pi}\int_{\partial \Omega} \partial_{\bn} v(\bx) \int_{\partial \Omega} \log|\bx-\bz| \sigma(\bz;\by)\,{\rm d}S(\bz)\,{\rm d}S(\bx) +\alpha(\by) \int_{\partial \Omega} \partial_{\bn} v(\bx)\,{\rm d}S(\bx) = -\mu(\by).$$
We note that
$$\int_{\partial \Omega} \partial_{\bn} v(\bx)\,{\rm d}S(\bx) = \frac{1}{2}\int_{\partial\Omega} \bx\cdot \bn\,{\rm d}S(\bx) = |\Omega|,$$
and hence 
\begin{equation}\label{eqn:cond}
 \alpha(\by) = \frac{1}{2\pi|\Omega|}\int_{\partial \Omega} \partial_{\bn} v(\bx) \int_{\partial \Omega} \log|\bx-\bz| \sigma(\bz;\by)\,{\rm d}S(\bz)\,{\rm d}S(\bx) -\frac{1}{|\Omega|}\mu(\by). 
\end{equation}
Inserting this identity into (\ref{eq:int_eq2a}) gives
\begin{align}\label{eqn:inteq_final}
        \frac{1}{2}\sigma(\bx;\by) &-\frac{1}{2\pi}\int_{\partial \Omega} \frac{(\bx-\bz)\cdot \bn}{|\bx-\bz|^2}\sigma(\bz;\by)\,{\rm d}S(\bz) + \int_{\partial \Omega}\sigma(\bz;\by)\,{\rm d}S(\bz) \\
    &  -\frac{1}{2\pi|\Omega|}\int_{\partial \Omega}\left[\int_{\partial \Omega} \partial_{\bn} v(\br)  \log|\br-\bz|\,{\rm d}S(\br) \right] \sigma(\bz;\by)\,{\rm d}S(\bz)\\
    &=  -\frac{1}{|\Omega|}\partial_{\bn} v(\bx)-\partial_{\bn} G_0(\bx;\by)- \frac{1}{|\Omega|}\mu(\by).\nonumber 
\end{align}
Once the linear system \eqref{eqn:inteq_final} has been solved, the final solution for the regular part is then
\begin{equation}\label{eqn:Reg_Final}
    \Rint(\bx;\by) = \alpha(\by) + \frac{|\bx|^2}{4|\Omega|} -\frac{1}{2\pi} \int_{\partial\Omega} \log|\bx-\bz| \,\sigma(\bz;\by) {\rm d}S(\bz).
\end{equation}
For evaluation off surface $(\bx\notin\partial\Omega)$, explicit differentiation yields the gradient and Hessians as
\bsub\label{eq:Diffs}
\begin{align}
\label{eq:Diffs_a}\partial_{x_i} \Rint(\bx;\by) &= \frac{x_i}{2|\Omega|}-\frac{1}{2\pi} \int_{\partial\Omega} \frac{x_i-z_i}{|\bx-\bz|^2} \,\sigma(\bz;\by) {\rm d}S(\bz),\\[5pt]
\label{eq:Diffs_b}\partial^2_{x_ix_j} \Rint(\bx;\by) &=  \frac{\delta_{ij}}{2|\Omega|} -\frac{1}{2\pi}\int_{\partial\Omega} \left[ \frac{\delta_{ij}}{|\bx-\bz|^2} - \frac{2(x_i-z_i)(x_j-z_j)}{|\bx-\bz|^4} \right]\,\sigma(\bz;\by) {\rm d}S(\bz),
\end{align}
\esub
and the above integrals are smooth. For determining derivatives on surface $\bx\in\partial\Omega$, a more involved process of applying the normal derivative \eqref{eq:int_eq} together with spectral differentiation to compute tangential derivatives can be used to find these quantities.

\edit{We remark that when the function $\int_{\partial\Omega}\partial_{\bn} v(\mathbf{r})\log|\mathbf{r}-\mathbf{z}|\,dS(\mathbf{r})$ 
is not orthogonal to the (one-dimensional) kernel of the homogeneous adjoint operator 
in \eqref{eq:int_eq}, the reduced operator obtained by setting $\alpha\equiv 0$ is itself 
invertible, and the constant $\alpha$ may be omitted. We expect that this condition holds generically but 
not universally; for the unit disk, for instance, $\partial_{\bn} v \equiv \tfrac{1}{2}$ and 
$\int_{|\mathbf{r}|=1}\log|\mathbf{r}-\mathbf{z}|\,dS(\mathbf{r}) = 0$ for $\mathbf{z}\in\partial\Omega$, 
so the condition fails.}

\edit{\subsection*{Exterior problem.} The integral equation for the exterior problem is a known adaptation (see~\cite{kress} for example) of our approach which we include here for completeness. As above, we first decompose $\Gext$ into a sum of the free-space Green's function and a ``scattered field'' $w,$
\begin{align*}
    w(\bx;\by)= \Gext(\bx;\by) - G_0(\bx;\by). 
\end{align*}
We then represent $w$ via the ansatz
\begin{align}
    w(\bx;\by) = -\frac{1}{2\pi} \int_{\partial \Omega} \log|\bx-\bz| \sigma(\bz;\by)\,{\rm d}S(\bz),
\end{align}
for an unknown boundary density $\sigma.$ Inserting this expression into the BVP and taking appropriate limits to the boundary yields the integral equation
\begin{align}\label{eq:int_eq_ext}
    -\frac{1}{2}\sigma(\bx;\by) -\frac{1}{2\pi}\int_{\partial \Omega} \frac{(\bx-\bz)\cdot \bn}{|\bx-\bz|^2}\sigma(\bz;\by)\,{\rm d}S(\bz) = -\partial_{\bn} G_0(\bx;\by). 
\end{align}
We remark that the sign difference in the term $-\frac{1}{2} \sigma$ when compared to \eqref{eq:int_eq} is due to the convention that the normal $\textbf{n}$ points outwards from $\partial\Omega$ into $\mathbb{R}^2\setminus\Omega$ for both the interior and exterior problems. For simply connected domains, the invertibility of (\ref{eq:int_eq_ext}) is well-known~\cite{kress}. Additionally, applying the divergence theorem to the right-hand side of (\ref{eq:int_eq_ext}) gives $\int_{\partial \Omega} \partial_{\bn}\,G_0{\rm d}S = 0$ and hence $\int_{\partial \Omega} \sigma\,{\rm d}S = 0.$ In particular, $w = o(1)$ as $|\bx| \to \infty.$
}

\subsection*{Discretization of the integral equation}
In two dimensions we discretize the integral equation \eqref{eqn:inteq_final} using a standard patch-based collocation method, in which the unknowns of the system are values of the unknown density at a set of discretization nodes. For ease of exposition, we write the integral equation \eqref{eqn:inteq_final}  as $\sigma + \mathcal{K}[\sigma] = f,$ where $\mathcal{K}_{\by}$ is a compact weakly-singular integral operator. Let $K_{\by}$ denote its kernel. In two dimensions, the boundary $\partial \Omega$ is divided into panels, each of which is discretized using a $k^{th}$ order Gauss-Legendre rule. Let $\{\bx_j\}_{j=1}^N$ denote the collection of discretization nodes on $\partial \Omega.$ We then enforce the integral equation at the discretization nodes, yielding a set of $N$ equations
\begin{align}\label{eqn:partial_disc}\sigma(\bx_j;\by)+ \int_{\partial \Omega} K_{\by}(\bx_j; \bz) \sigma(\bz; \by)\,{\rm d}S(\bz) = f(\bx_j), \quad j=1,\ldots,N.
\end{align}
In order to fully discretize the system we require suitable quadrature nodes $w_{j,k}$ for computing the integrals on the left-hand side of \eqref{eqn:partial_disc}. For points which are far away, we can choose $w_{j,k}$ independent of $j$. Typically these quadrature weights $w_k$ are chosen so that the quadrature rule with nodes $\bx_j$ and weights $w_j$ accurately integrate sufficiently smooth functions on $\partial \Omega.$ For $\bx_j$ near $\bx_k,$ custom quadrature formulae based on a mixture of adaptive integration and \emph{generalized Gaussian quadratures} are used. Let $\mathcal{S}_j$ denote the set of $k$ such that $w_{j,k}$ requires such a correction. We note that if patches are approximately uniform in size, that $\mathcal{S}_j$ will typically be $\bigoh(1)$ for large $N$. Let $W_{j,k}$ denote the corresponding corrections. We can then replace \eqref{eqn:partial_disc} by the approximation
\begin{align}\label{eqn:fully_disc}
\sigma_j + \sum_{k \in \mathcal{S}_j} W_{j,k} \sigma_k + \sum_{k \notin \mathcal{S}_j} K_{\by}(\bx_j;\bx_k)\sigma_k w_k = f_j, \quad j=1,\ldots,N,
\end{align}
where $\sigma_j \approx \sigma(\bx_j;\by)$ and $f_j = f(\bx_j).$ Moreover, we note that the first term in \eqref{eqn:fully_disc} is sparse and can be applied in $\bigoh(N)$ operations. The second term can be computed using many standard fast algorithm techniques, including the fast multipole method~\cite{fmm2dlib} which reduces the cost of applying the discretized integral operators from $\bigoh(N^2)$ to $\bigoh(N)$, and fast direct methods \cite{Ho2020} which construct a highly-compressed factorization of the inverse of the system matrix. The implementation of these quadratures, together with geometric calculations, are performed using the integral equation package 
{\tt chunkIE}\footnote{available at \url{https://chunkie.readthedocs.io}} \cite{Askham_chunkIE_a_MATLAB_2024}. 

\subsection{Convergence studies}

In this section, we validate our method against existing known representations for the Green's functions in certain geometries, particularly disks and ellipses. \al{We validate the convergence of our method as quadrature order ($k$) and panel number increases. In the results shown in the below sections, we have monitored the condition number of the system matrix and verified a well-understood benefit \cite{kress} of the boundary integral equation approach - that the condition number approaches a finite value under refinement, rather than growing with the number of unknowns $N_{\text{pts}}$.}

\subsubsection{The disk domain.}

\begin{figure}[htbp]
    \centering
    \subfigure[Regular part $\Rint_b$.]{\includegraphics[width = 0.325\textwidth]{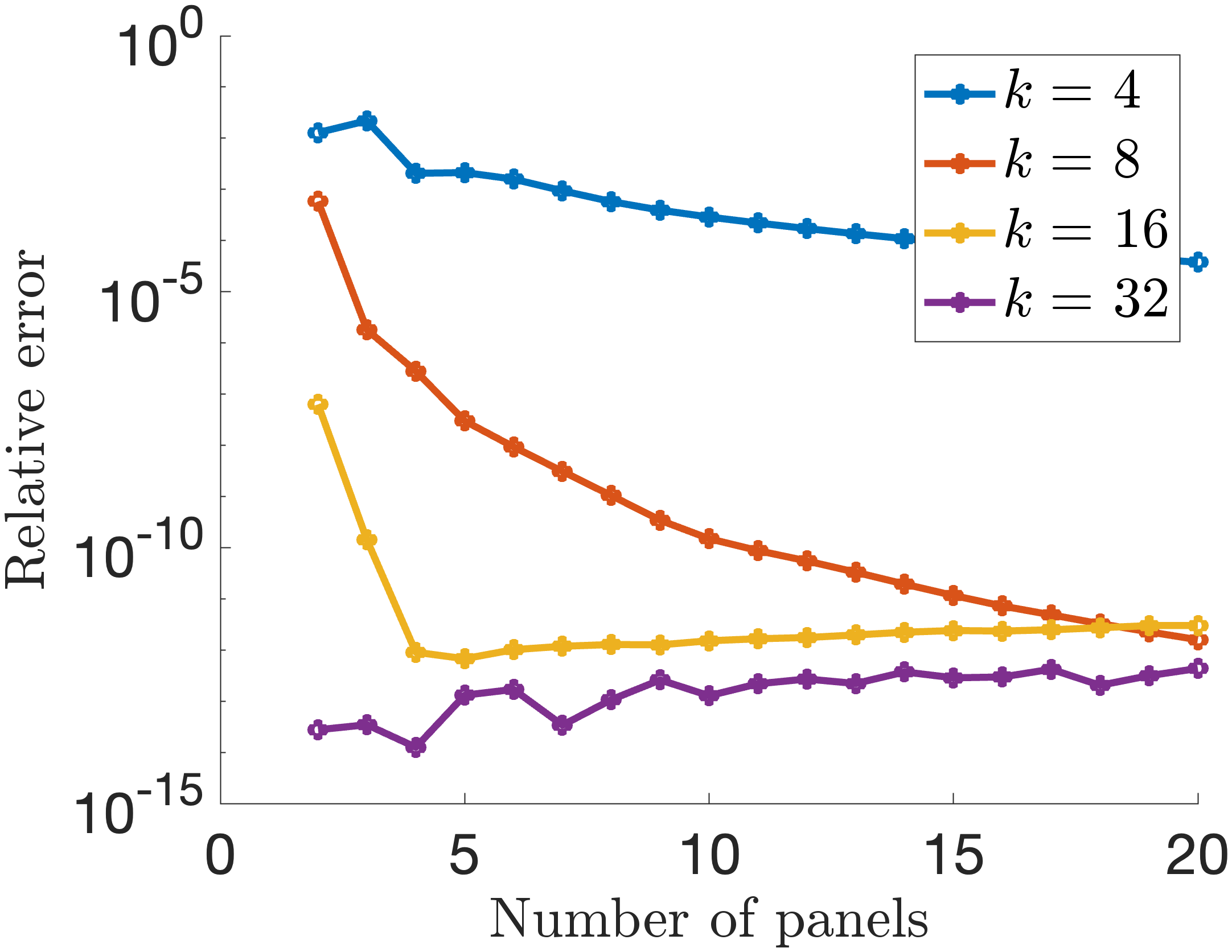}}
    \subfigure[Gradient $\|\nabla \Rint_b\|_2$.]{\includegraphics[width = 0.325\textwidth]{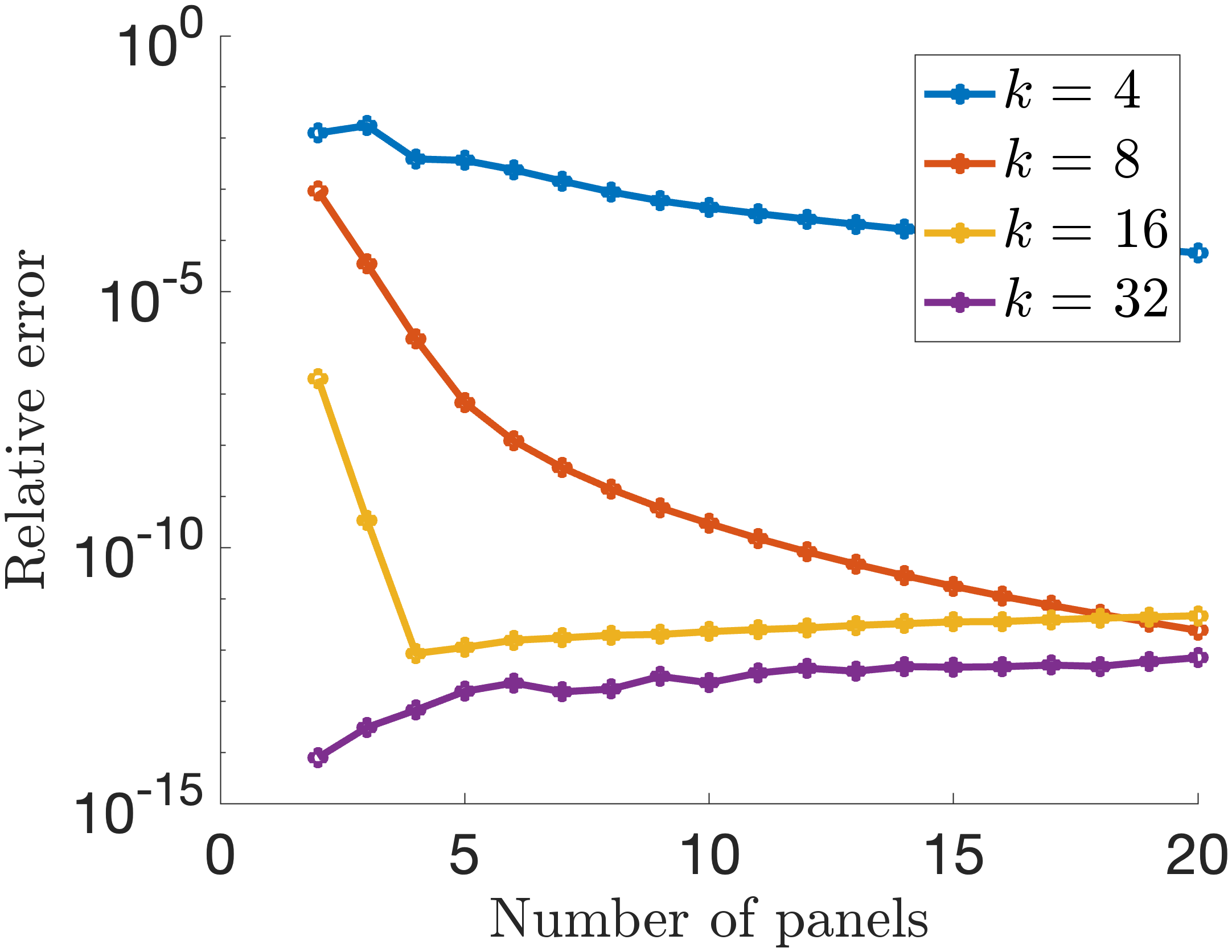}}
    \subfigure[$\|\nabla^2 \Rint_b\|_2$.]{\includegraphics[width = 0.325\textwidth]{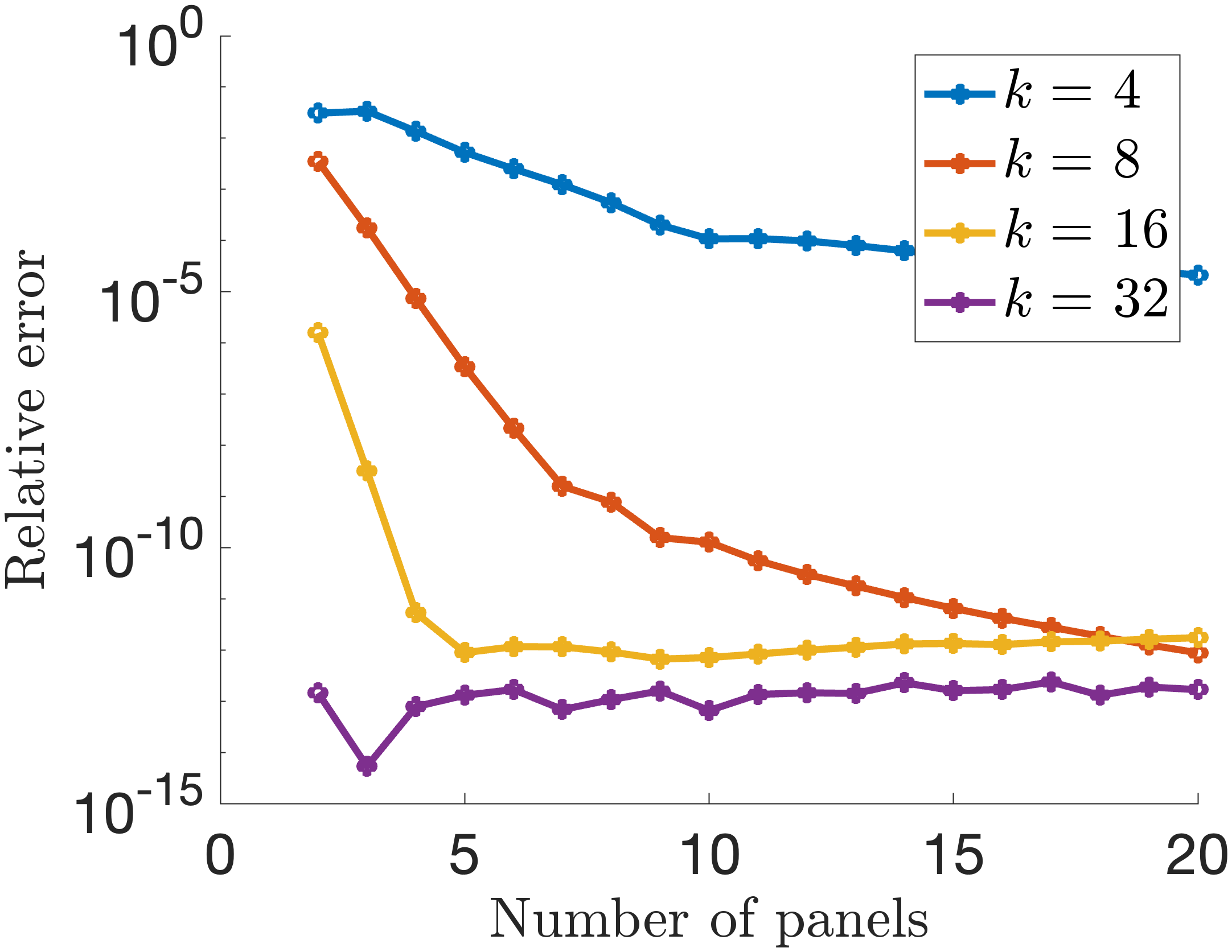}}
\caption{Convergence of the relative errors in the case of the unit disk against panel number and quadrature order $k$. Numerical approximation of the regular part $\Rint_b$ (Panel (a)), the gradient $\|\nabla \Rint_b\|_2$ (Panel (b)) and the Hessian $\|\nabla^2 \Rint_b \|_2$ (Panel (c)) for the point $\emph{\by} = (\frac14,\frac13)$.
\label{fig:ConvergenceDisk}}
\end{figure}

In the case of the unit disc, we have a closed form for the interior Green's function \cite{KTW_2005}. We can identify the regular part $\Rint(\bx;\by)$ together with its gradient and Hessian. 
In the case of the disk domain $\Omega = \{ \bx = (x_1,x_2) \ | \ x_1^2 + x_2^2 \leq1\}$ and source $\by = (y_1,y_2)$, we have \cite{KTW_2005} that
\bsub\label{eq:GreensDisk}
\begin{align}
\Gint_b(\bx;\by) &= \frac{-1}{2\pi} \log|\bx-\by| + \Rint_b(\bx;\by);\\[5pt]
\Rint_b(\bx;\by) &= -\frac{1}{2\pi} \left[ \frac12 \log(1+ |\bx|^2 |\by|^2 - 2 \bx\cdot \by) - \frac{1}{2}(|\by|^2 + |\bx|^2) + \frac{3}{4} \right].
\end{align}
\esub
From the expression \eqref{eq:GreensDisk}, we can easily calculate the gradient and Hessian. In Fig.~\ref{fig:ConvergenceDisk}, we plot convergence of the regular part $|\Rint_b|$, its gradient $\|\nabla \Rint_b\|_2$ and Hessian $\|\nabla^2 \Rint_b\|_2$ as a function of panel number and quadrature order, $k$ for $\bx=\by = [\frac14,\frac13]$. We observe excellent agreement between the numerical approximation and the exact solution \eqref{eq:GreensDisk}. We observe that for the highest order $(k=32)$ that the method is able to reproduce the solution and its derivatives up to a relative error of close to $\bigoh(10^{-15})$.

\begin{figure}[htbp]
\centering
\subfigure[Error in $\Rint_s(\bx;\by)$.]{\includegraphics[width = 0.325\textwidth]{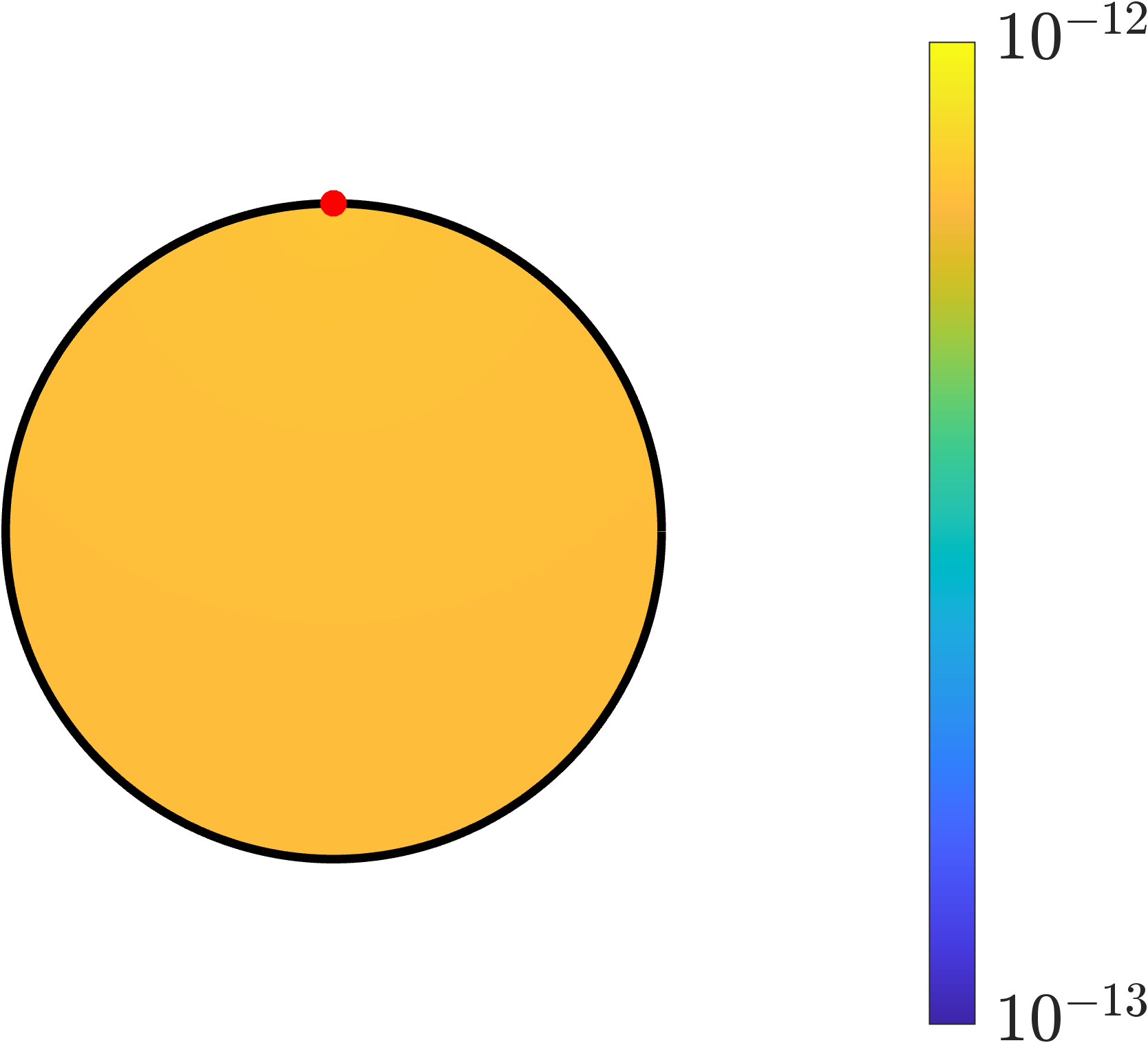}}
\hspace{1cm}
\subfigure[Error in \al{$\Rext_s(\bx;\by)$}.]{\includegraphics[width = 0.38\textwidth]{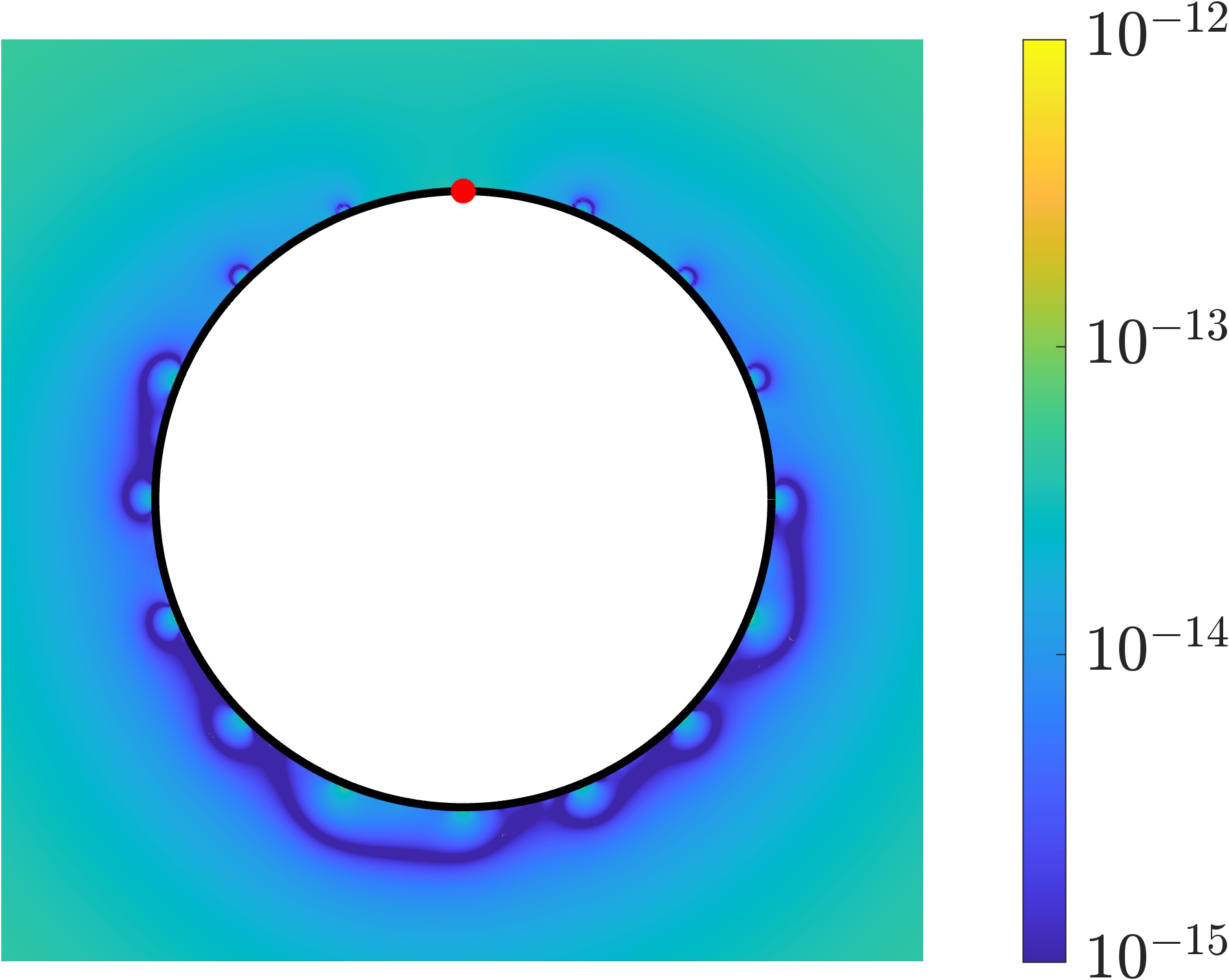}}
\caption{Errors for the surface Green's function for the disk. Panel (a) shows the absolute error in $\Rint_s(\emph{\bx;\by})$ with the exact solution given in \eqref{eq:GreensDiskSurf}. Panel (b) shows the absolute error in $\Rext_s(\emph{\bx;\by})$ with the exact solution given in \eqref{eq:GreensDiskSurfExt}. In both cases the source point is $\emph{\by} = [0,1]$ (solid red dot).
\label{fig:Disk_surf}}
\end{figure}

Additionally, the surface Green's function for the disk can be calculated from \eqref{eq:GreensDisk} by noticing that
\[1+|\bx|^2|\by|^2 - 2 \bx\cdot\by = (\bx|\by| - \by/|\by|)\cdot(\bx|\by| - \by/|\by|) = |\bx-\by|^2, \quad \mbox{for} \quad |\by|=1.
\]
Hence taking the direct limit as $\by\to\partial\Omega$ in \eqref{eq:GreensDisk}, we obtain
\begin{equation}\label{eq:GreensDiskSurf}
\Gint_s(\bx;\by) = \frac{-1}{\pi} \log|\bx-\by| + \Rint_s(\bx;\by);\qquad
\Rint_s(\bx;\by) = -\frac{1}{8\pi}  +  \frac{|\bx|^2}{4\pi}.
\end{equation}
Another case where explicit solutions are available is the external surface function $\Gext_s$ satisfying \eqref{eqn_neumGExt}, which was derived in \cite[Appendix B]{Lindsay2023} to be
\begin{equation}\label{eq:GreensDiskSurfExt}
\Gext_s(\bx;\by) = \frac{-1}{\pi} \log|\bx-\by| + \Rext_s(\bx;\by); \qquad
\Rext_s(\bx;\by) = \frac{1}{2\pi}\log|\bx|.
\end{equation}
In Fig.~\ref{fig:Disk_surf} we show a comparison between our numerical solution and the exact forms \eqref{eq:GreensDiskSurf} and \eqref{eq:GreensDiskSurfExt}. Similar to other examples, we observe that the boundary integral method resolves the solution to an accuracy around $\bigoh(10^{-12})$.

\subsubsection{The ellipse domain.} 
\begin{figure}[htbp]
    \centering
        \subfigure[$a=1.5$.]{\includegraphics[width = 0.325\textwidth]{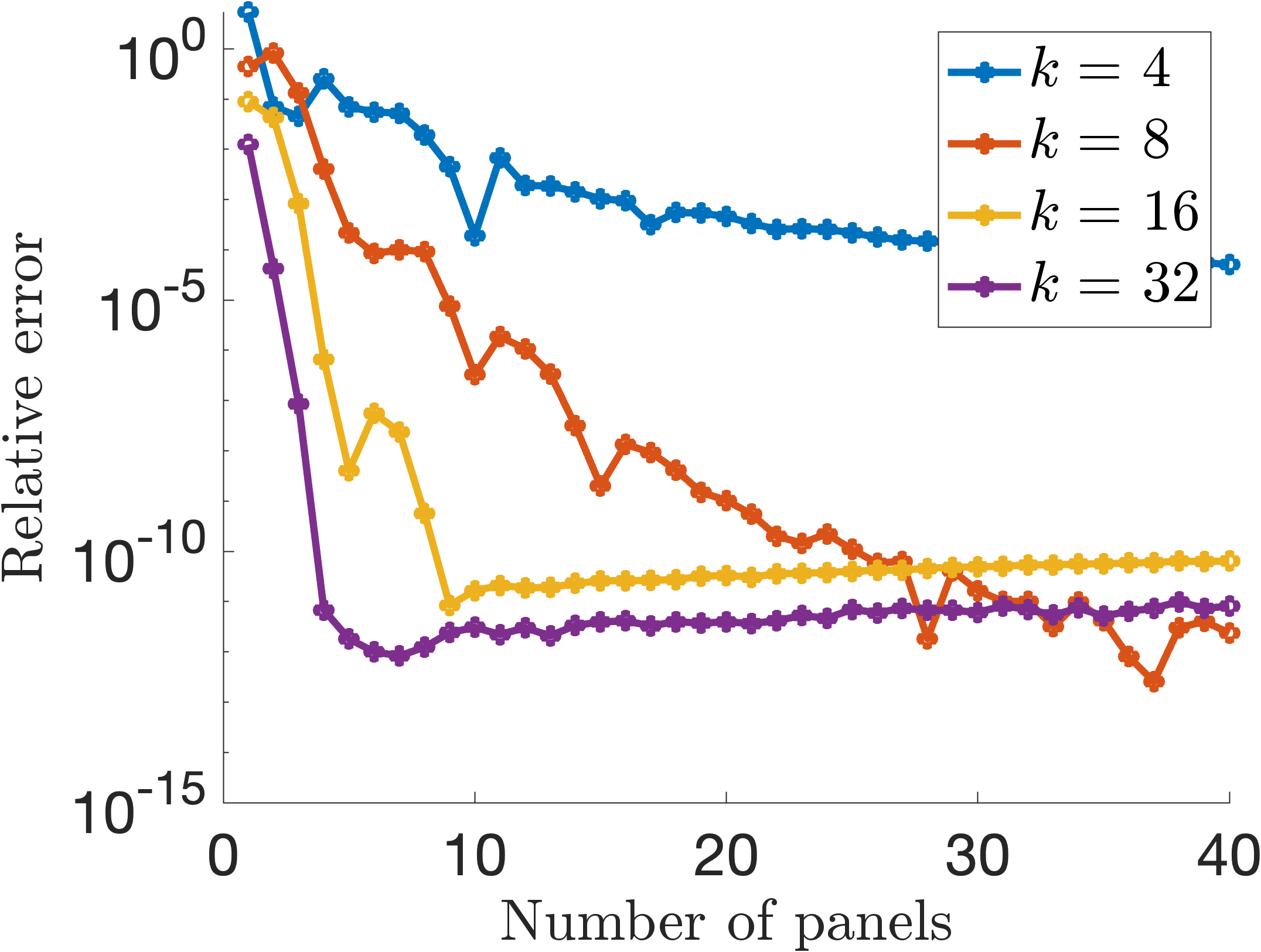}}
        \subfigure[$a=2.0$.]{\includegraphics[width = 0.325\textwidth]{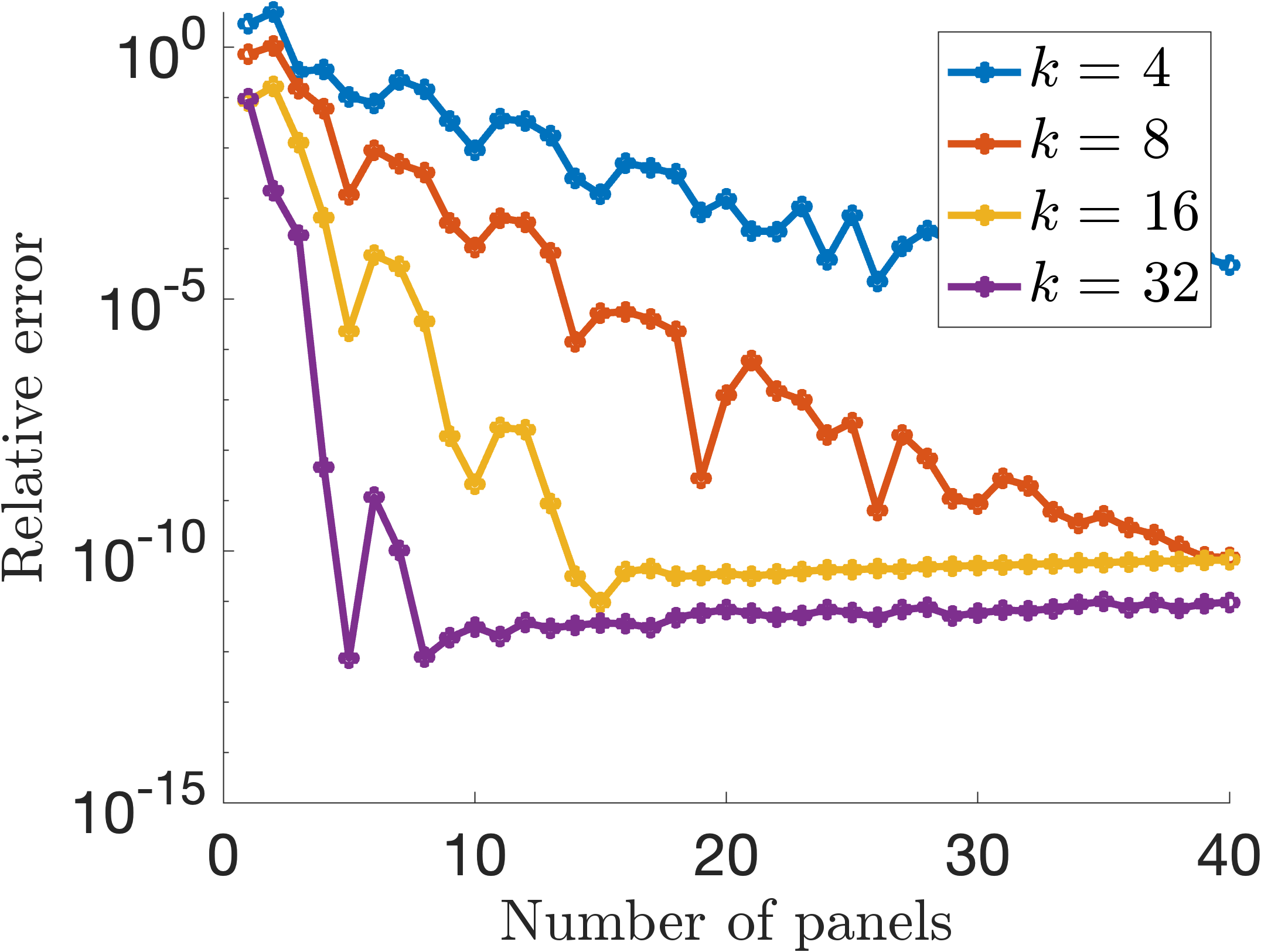}}
        \subfigure[$a=2.5$.]{\includegraphics[width = 0.325\textwidth]{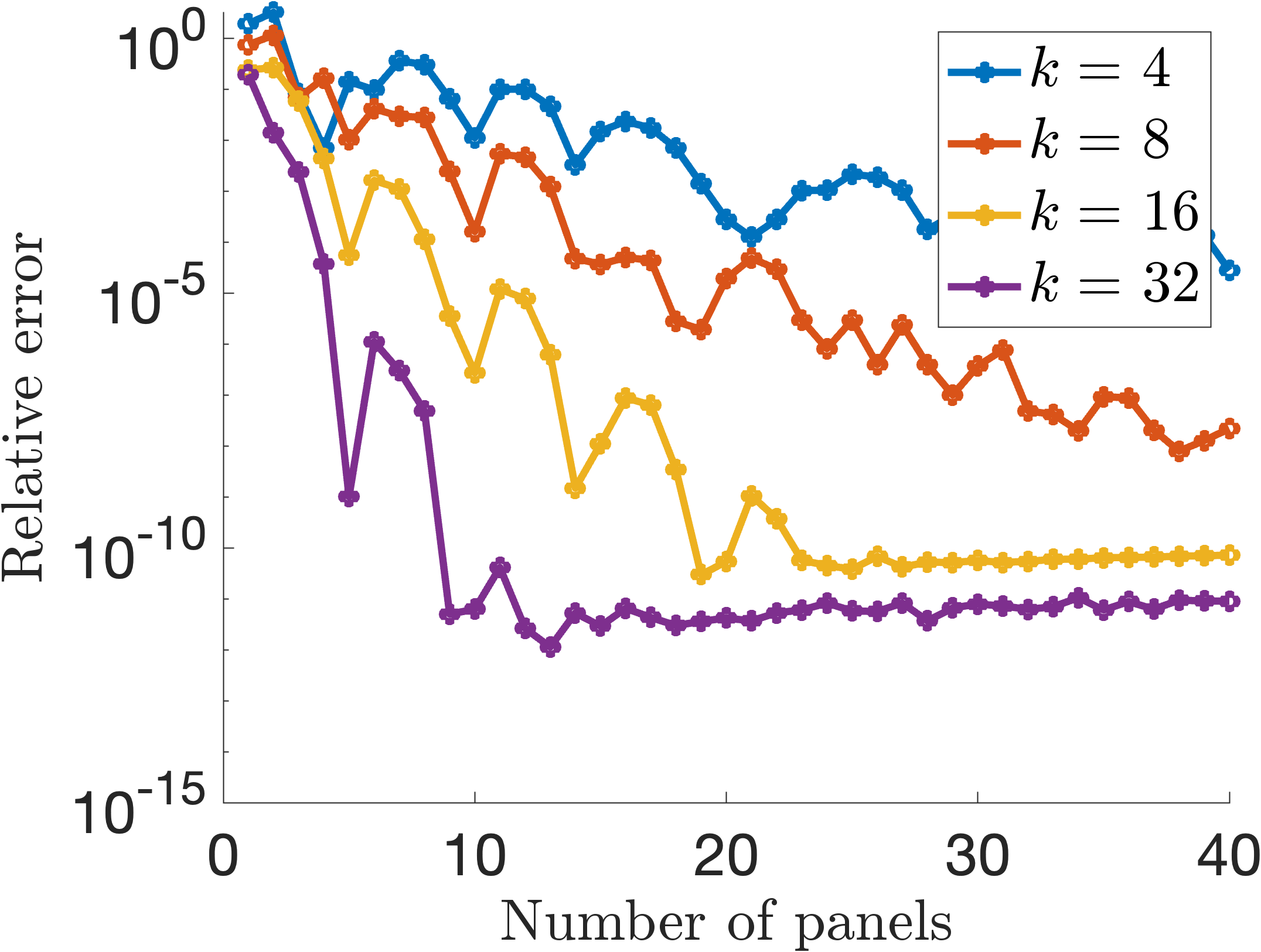}}
        \caption{Convergence of relative errors of $\Rint_b\emph{(\by;\by)}$ for elliptical domains $(x_1/a)^2 + (x_2/b)^2 \leq1$ and quadrature orders $k$. The point $\emph{\by} =[a/4,b/3]$ with $|\Omega| = \pi a b = \pi$ and various values of $a$ ($b=1/a$). \label{fig:ellipseConvergence}}
\end{figure}

In the case where the domain is an ellipse, a rapidly convergent series solution for interior bulk Green's function $\Gint_b$ was recently derived (cf.~\cite{Sarafa2021}). This approach was based on a mapping of the elliptical geometry to a rectangle followed by a series solution. This Fourier series was subsequently amenable to re-summation which yielded a rapid and accurate solution of \eqref{eqn_neumG}. In Fig.~\ref{fig:ellipseConvergence}, we plot agreement between this series solution and our numerical method for ellipses of fixed area $|\Omega| = \pi$ varying semi-major axis $a$. For each $a$, we observe an achieved accuracy of roughly $\bigoh(10^{-12})$.

\begin{figure}[htbp]
    \centering
        \subfigure[$\Rext_{s}(\by)$ and $a=1.5$.]{\includegraphics[width = 0.325\textwidth]{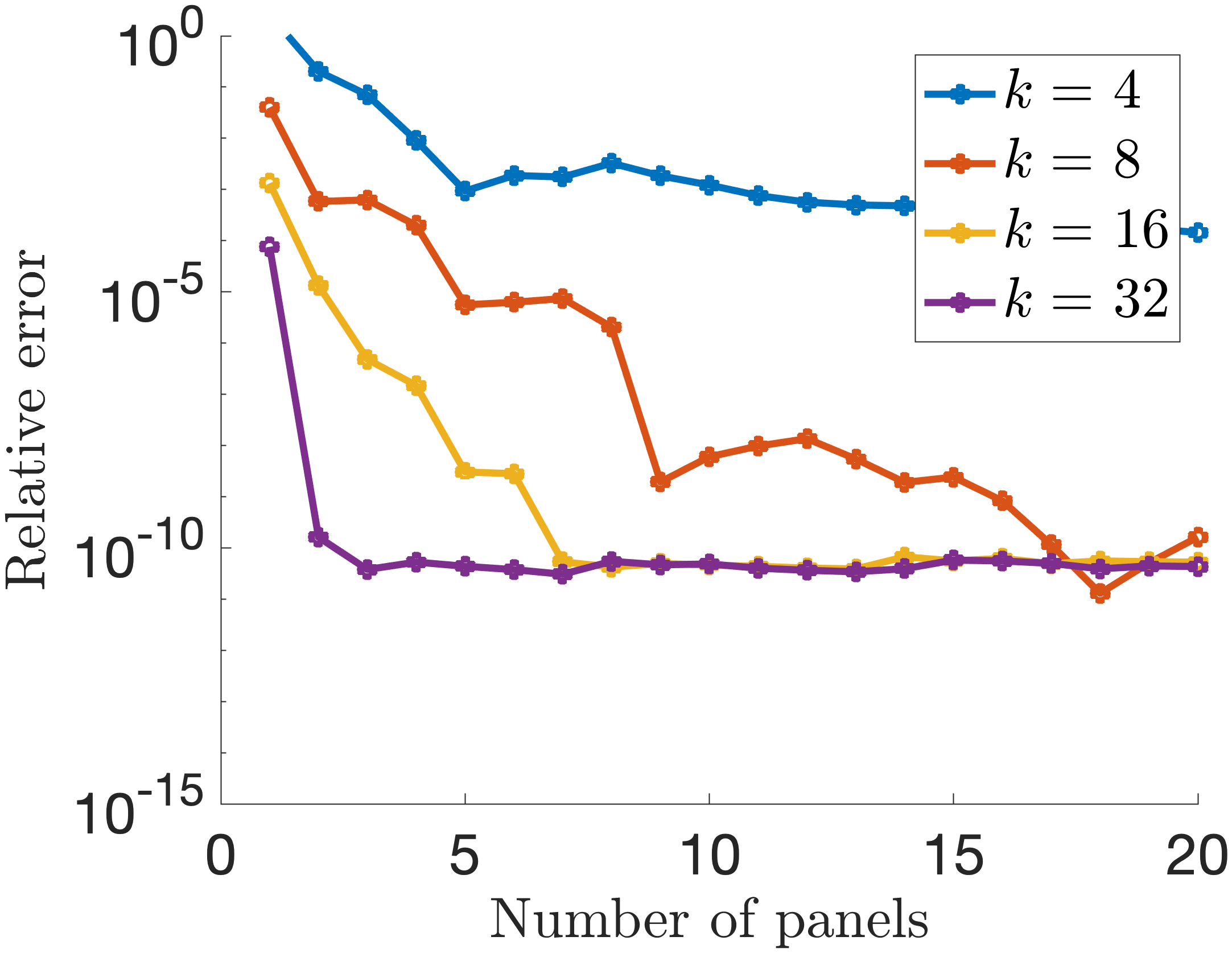}}
        \subfigure[$\Rext_{s}(\by)$ and $a=2.0$.]{\includegraphics[width = 0.325\textwidth]{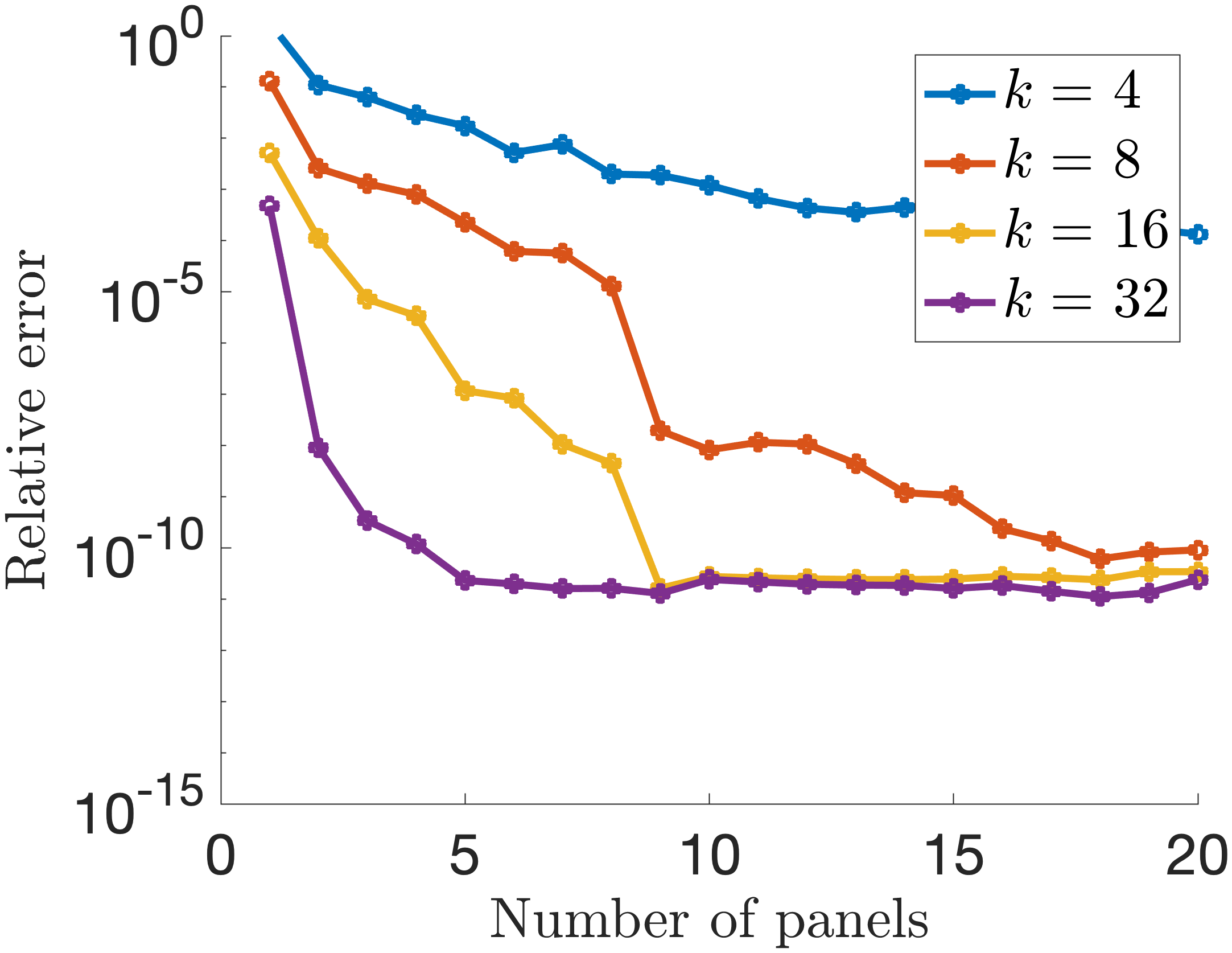}}
        \subfigure[$\Rext_{s}(\by)$ and $a=2.5$.]{\includegraphics[width = 0.325\textwidth]{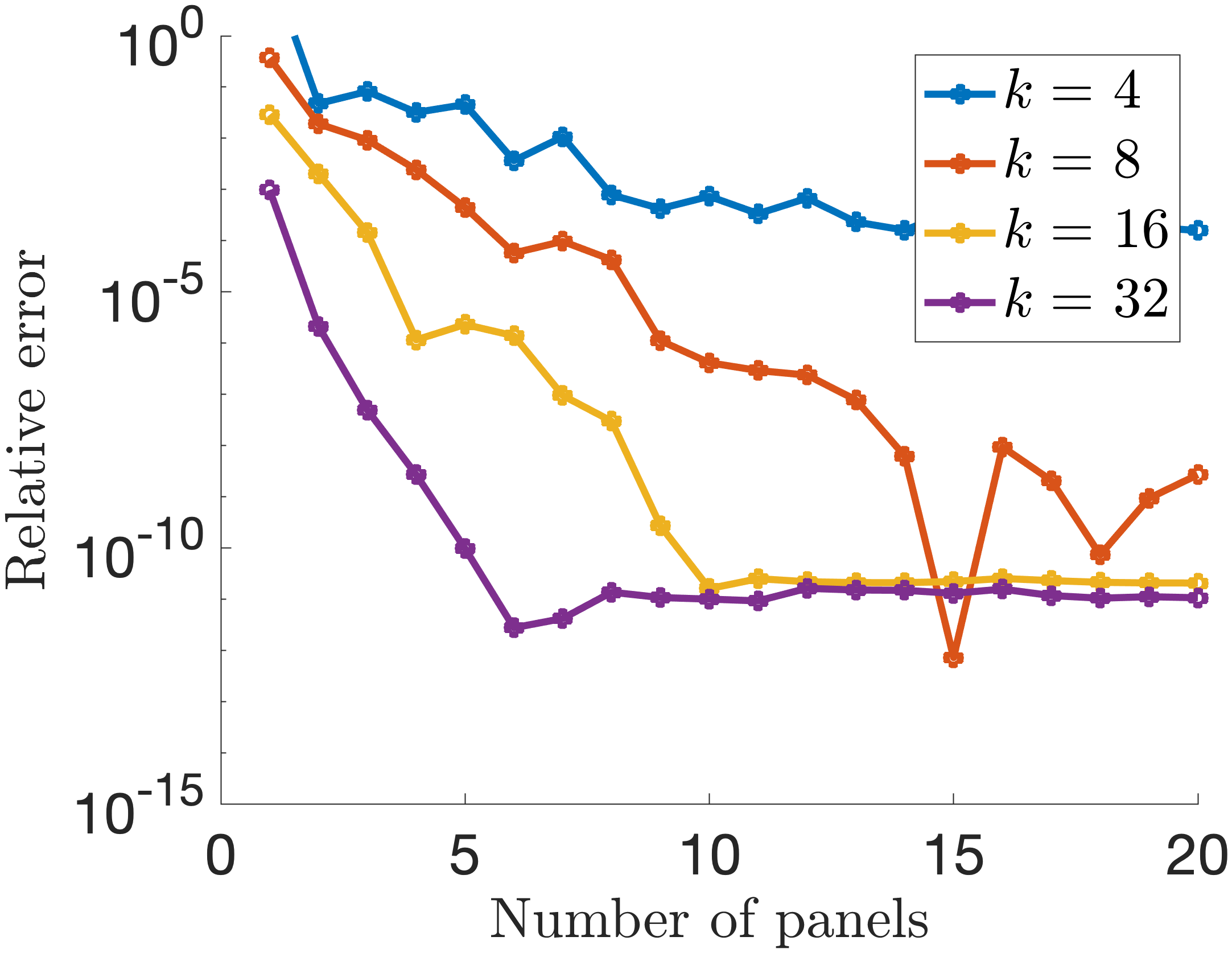}}\\
        \subfigure[$\Rint_{s}(\by)$ and $a=1.5$.]{\includegraphics[width = 0.325\textwidth]{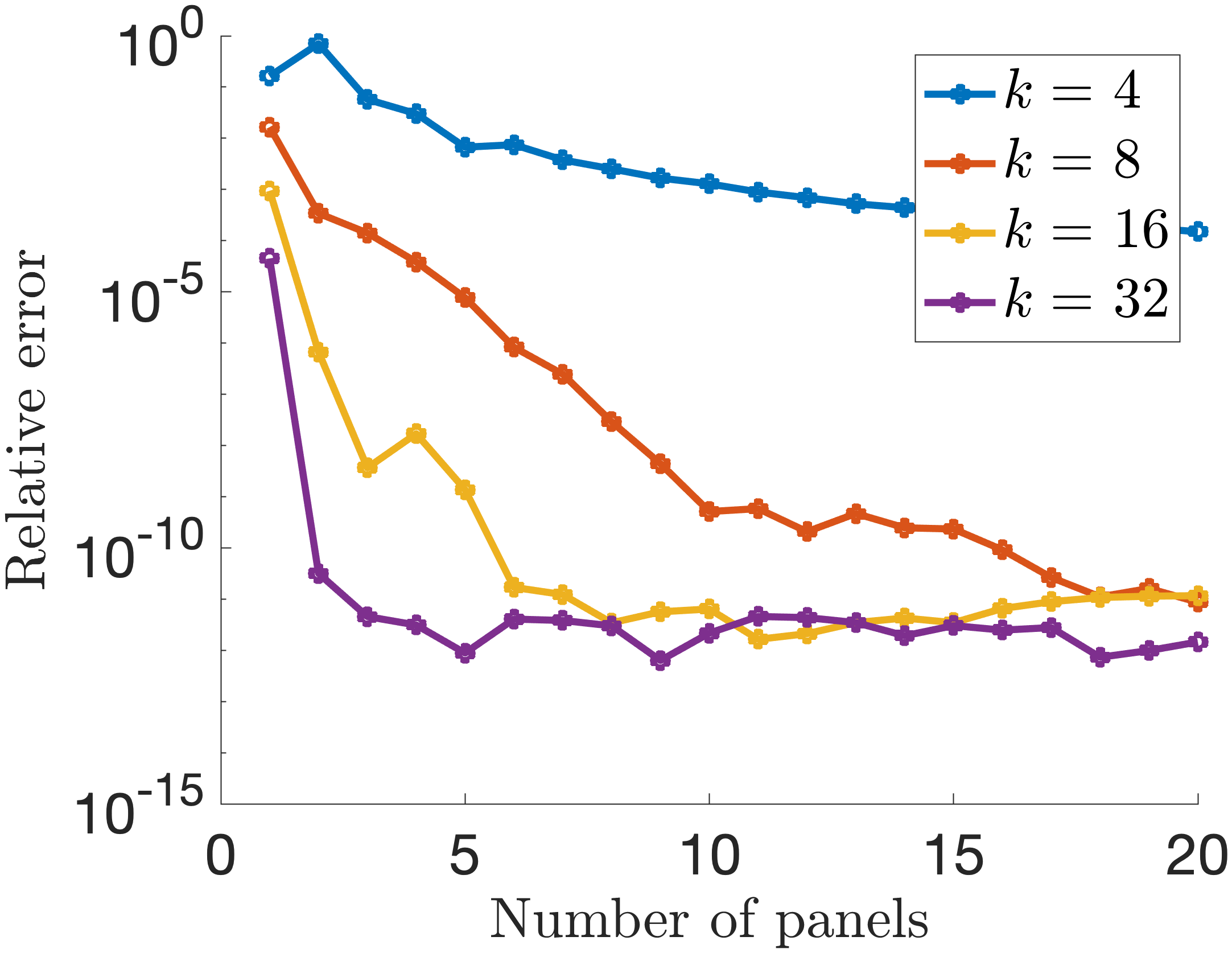}}
        \subfigure[$\Rint_{s}(\by)$ and $a=2.0$.]{\includegraphics[width = 0.325\textwidth]{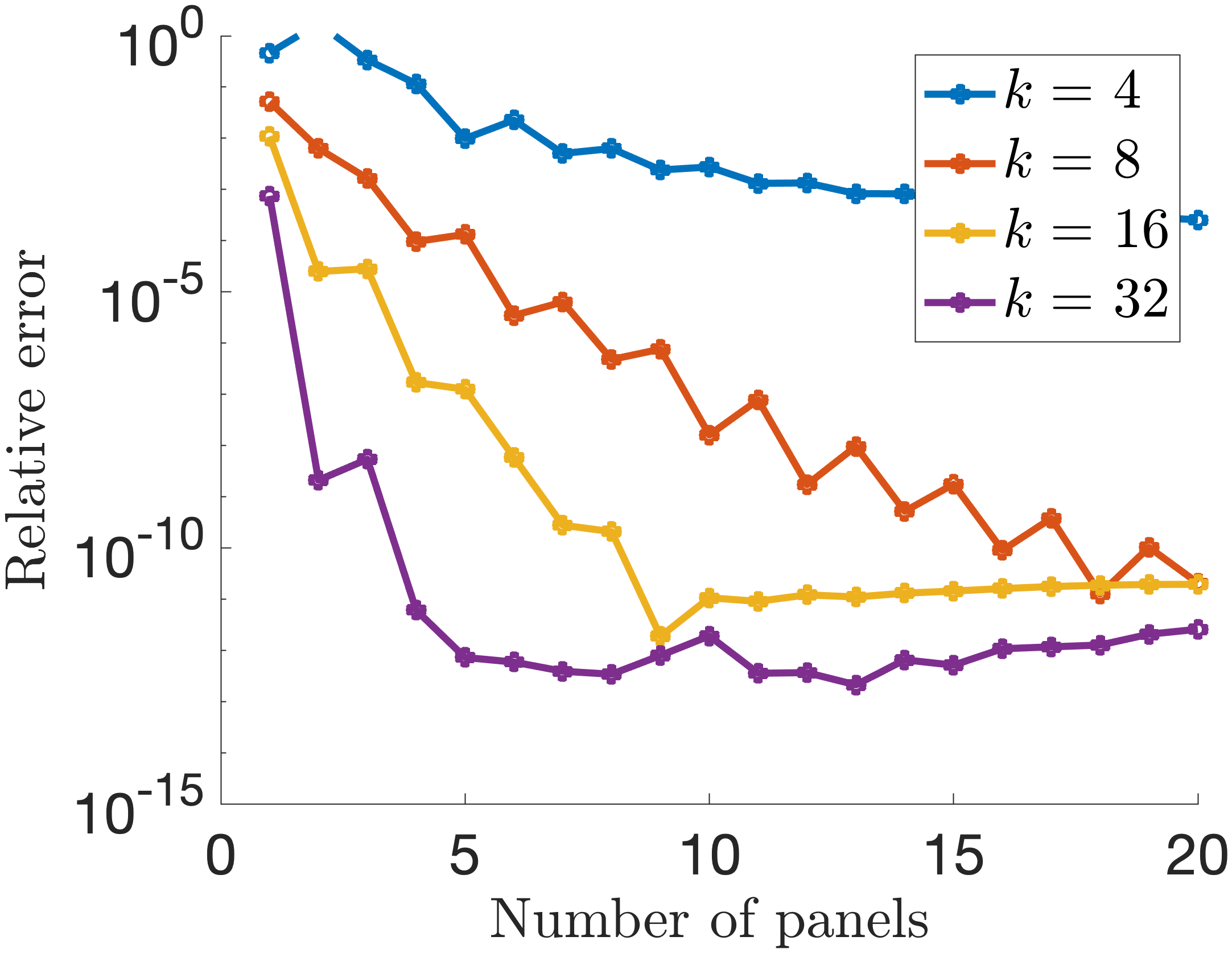}}
        \subfigure[$\Rint_{s}(\by)$ and $a=2.5$.]{\includegraphics[width = 0.325\textwidth]{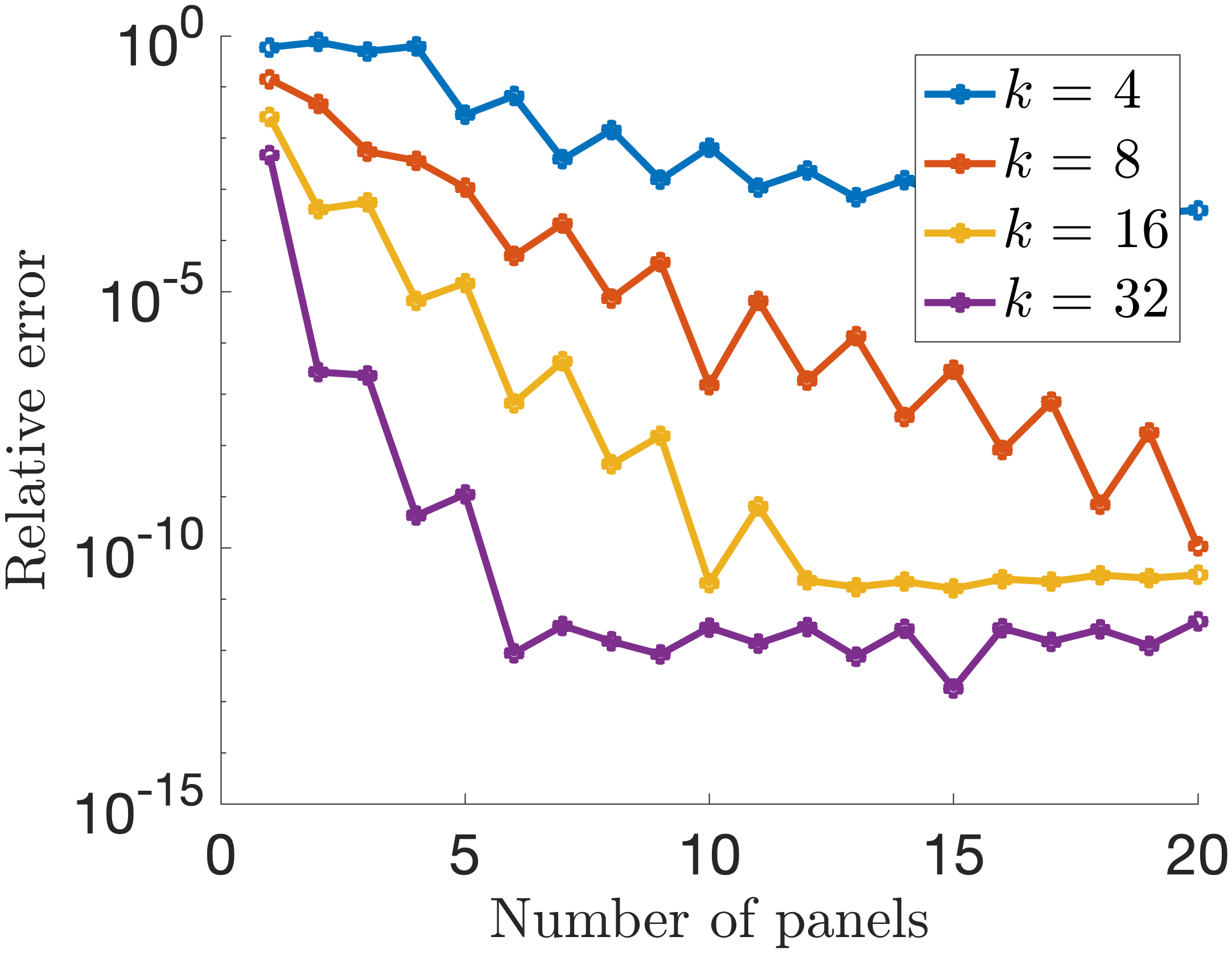}}
\caption{Convergence of the regular parts of the exterior $(\Rext_{s}$ - upper row) and interior $(\Rint_{s}$ - lower row) surface Green's functions for the elliptical domain $(x_1/a)^2 + (x_2/b)^2 = 1$. The area of the ellipses are fixed at $|\Omega| = \pi$ ($b = 1/a$) and from left to right results are presented for $a = 1.5$, $a= 2.0$ and $a= 2.5$. The method is able to reduce relative errors to roughly $\bigoh(10^{-12})$.  \label{fig:surfaceGreens}}
\end{figure}

In a more recent study \cite{GrebenkovWard2025a}, the interior surface ($\Rint_s$) and exterior surface ($\Rext_s$) Green's functions for an elliptical domain were derived. For a point $\by= (y_1,y_2)\in\partial\Omega$, the interior surface Green's function was derived (see \cite[Eqn.~(F17)]{GrebenkovWard2025a}) and its regular part $\Rint_{s}(\by) = \Rint_{s}(\by;\by)$ was found to be
\bsub\label{GreensEllipse_Surf}
\begin{align}\label{GreensEllipse_Surf_a}
   \nonumber \Rint_{s}(\by) &= \frac{|\by|^2}{2\pi a b} - \frac{3(a^2+b^2)}{16\pi a b} + \frac{1}{2\pi}\log \big( b^2 + (a^2-b^2) \sin^2\theta_0\big)\\
    & - \frac{2}{\pi} \sum_{n=1}^{\infty} \big( \log(1-\beta^{2n}) + \log | 1- \beta^{2n-1} e^{2 i \theta_0}|\big).
\end{align}
Here $\beta = (a-b)/(a+b)$ and $\tan\theta_0 = (ay_2)/(by_1)$ for $\by = (y_1,y_2)$. For the exterior counterpart, the regular part $\Rext_{s}(\by)$ of the surface Green's function was found (see \cite[Eqn.~(F30)]{GrebenkovWard2025a}) to be
\begin{equation}\label{GreensEllipse_Surf_b}
\Rext_{s}(\by) = \frac{1}{2\pi}\left[\log\Big(\frac{2}{a+b}\Big) + \log \big( b^2 + (a^2-b^2)\sin^2\theta_0\big)\right].
\end{equation}
\esub
In Fig.~\ref{fig:surfaceGreens} we demonstrate the convergence of the numerical method for both the exterior and interior regular parts \eqref{GreensEllipse_Surf}. We observe that the convergence rate is reduced somewhat for large $a$ as the ellipse becomes more elongated, however, with enough resolution, we are again able to reproduce the correct value to a relative error of around $10^{-12}$. \al{Overall, for quadrature order $k$, the error decays at a rate consistent with $\mathcal{O}(N_{\text{pts}}^{-k})$ before reaching around $10^{-12}$, as observed in Figs.~3.1, 3.3, 3.4.}

\subsection{Optimization routines.}\label{sec:Optim}

In this section, we describe our approach to the numerical minimization of the discrete energy $p(\bx_1,\ldots,\bx_N)$ described in \eqref{eq_Discrete}. The global minimization of this class of discrete energies is a long known computation challenge problem due to the large number of critical points, which grows in complexity as the dimensionality increases. In light of this, a number of optimization strategies combine local (e.g. Newton methods) and global search methods (e.g.~particle swarm optimization) \cite{Cheviakov2024,Gilbert_Cheviakov_2023,Sarafa2021}. In the present work, we are primarily interested in demonstrating the effectiveness of our numerical method when interfaced with standard optimization methods.

Our numerical method is implemented in \textsc{Matlab} and utilizes the boundary integral method package {\tt chunkIE}. This package handles geometry, implements high order quadrature routines and handles fast multipole method acceleration. Owing to this integration, the results we present here are obtained with the \textsc{Matlab} optimization routine {\tt fmincon}. We utilize the {\tt interior-point} algorithm and where possible specify the gradient of the objective function. In the case where we are optimizing over the centers of interior traps, the optimization routines have a tendency to take points outside the domain. To counteract this, we set a large value of the objective function when the points are outside the domain using the {\tt chunkerinterior} function which returns a boolean for whether $\bx\in\Omega$. In an effort to conduct a global search for minima, we initialize each routine with a random array of points in the domain and take the minima from many independent realizations of the optimization process. Other details can be viewed in our available code\footnote{}.

\section{Results}\label{sec:Results}

In this section we present the results of our integrated boundary integral and optimization method. We start by replicating some recent results in the literature \cite{Sarafa2021,Cheviakov2024,Gilbert_Cheviakov_2023}  followed by demonstrations of our method on some general geometries.

\begin{figure}[htbp]
    \centering
\includegraphics[width=0.19\textwidth]{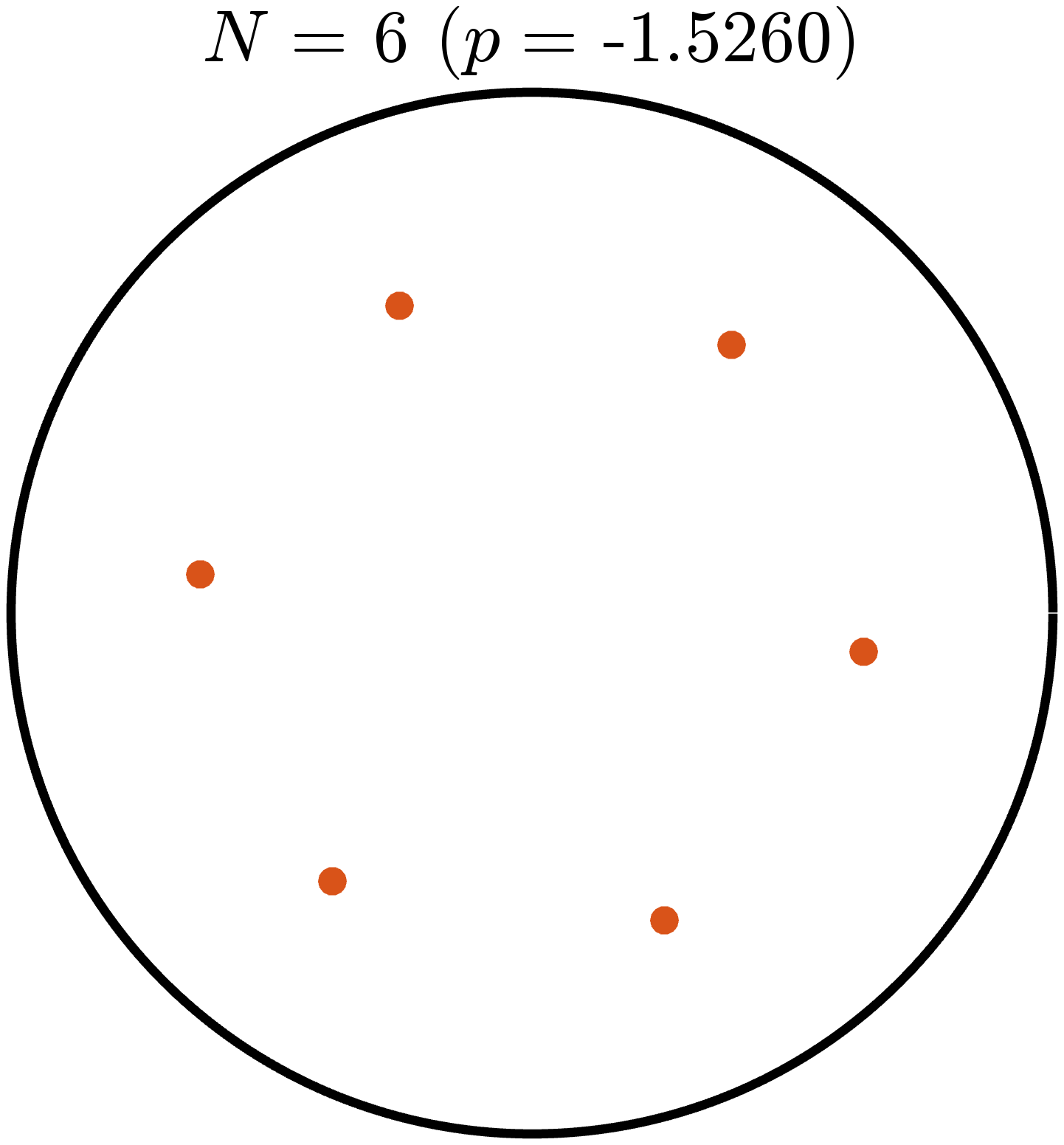}
\includegraphics[width=0.19\textwidth]{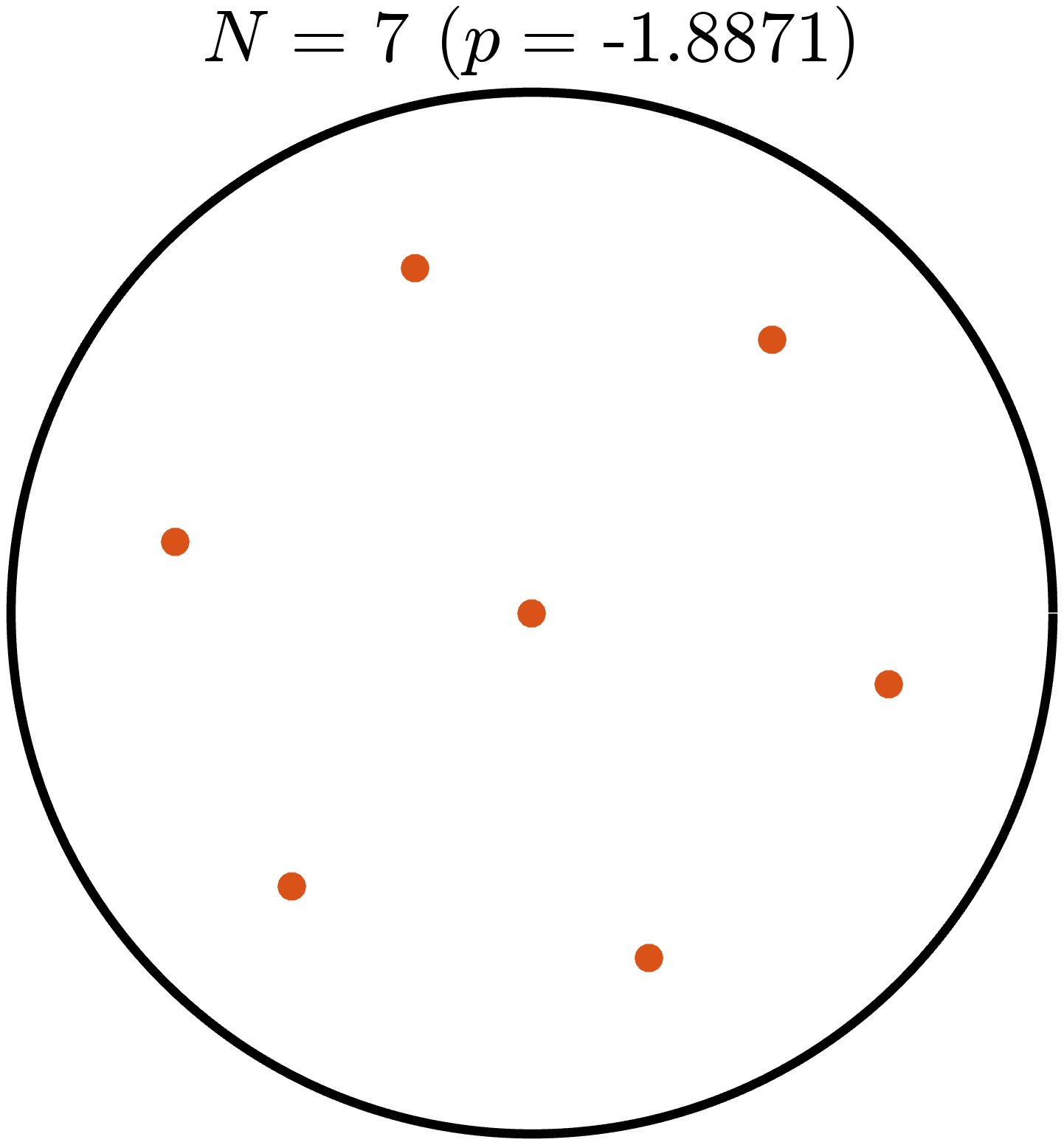}
\includegraphics[width=0.19\textwidth]{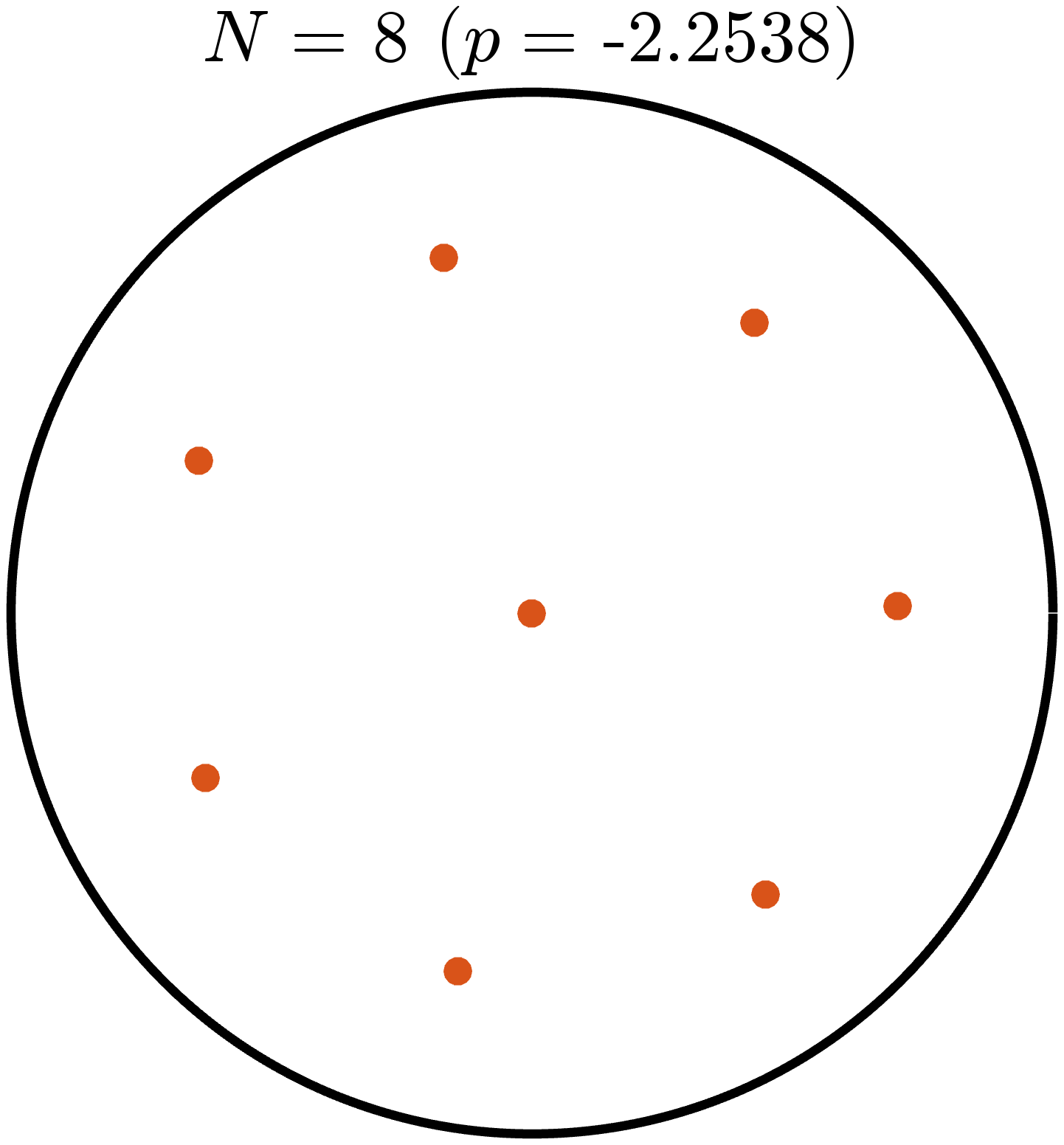}
\includegraphics[width=0.19\textwidth]{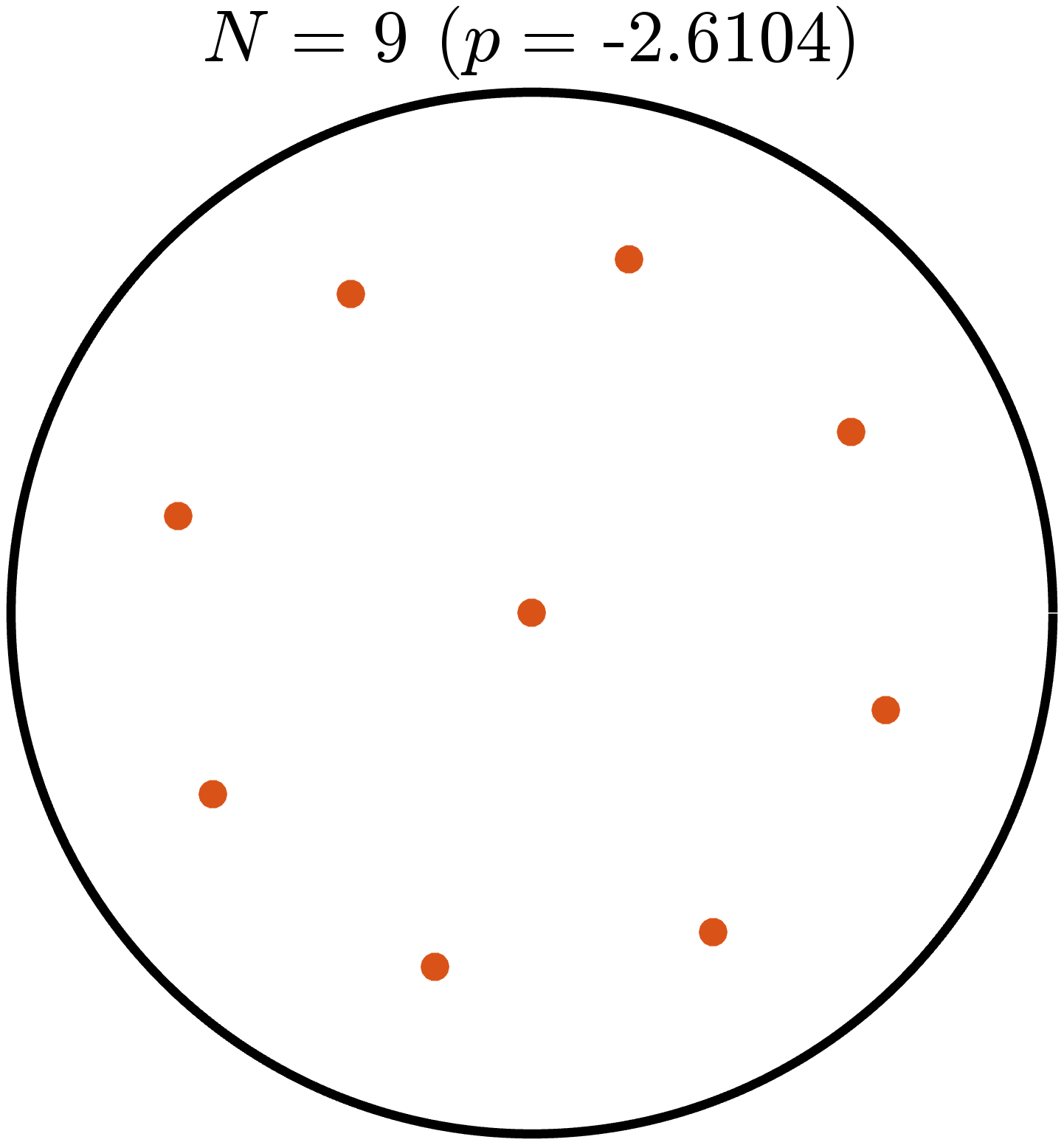}
\includegraphics[width=0.19\textwidth]{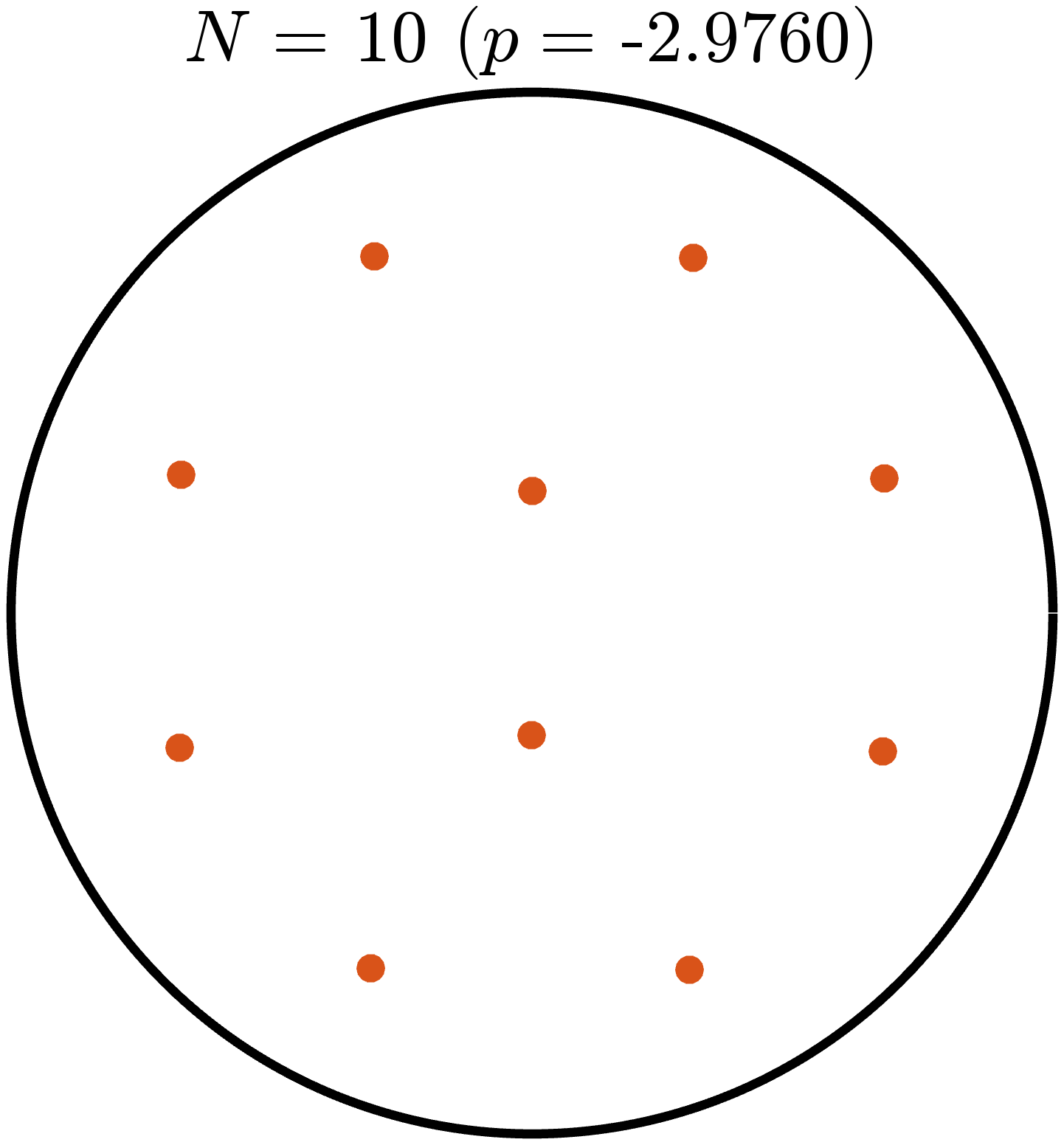}\\[5pt]
\includegraphics[width=0.19\textwidth]{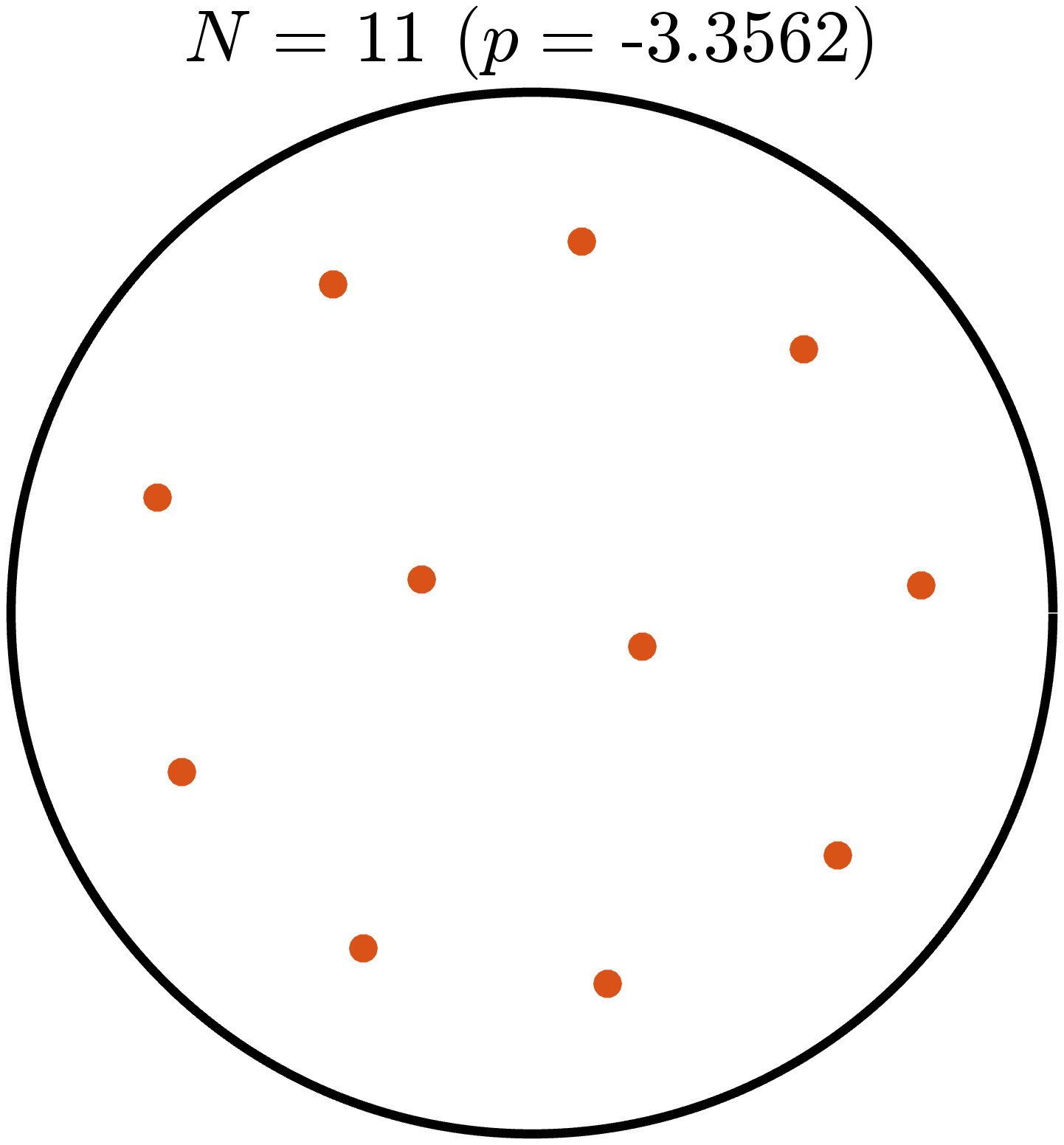}
\includegraphics[width=0.19\textwidth]{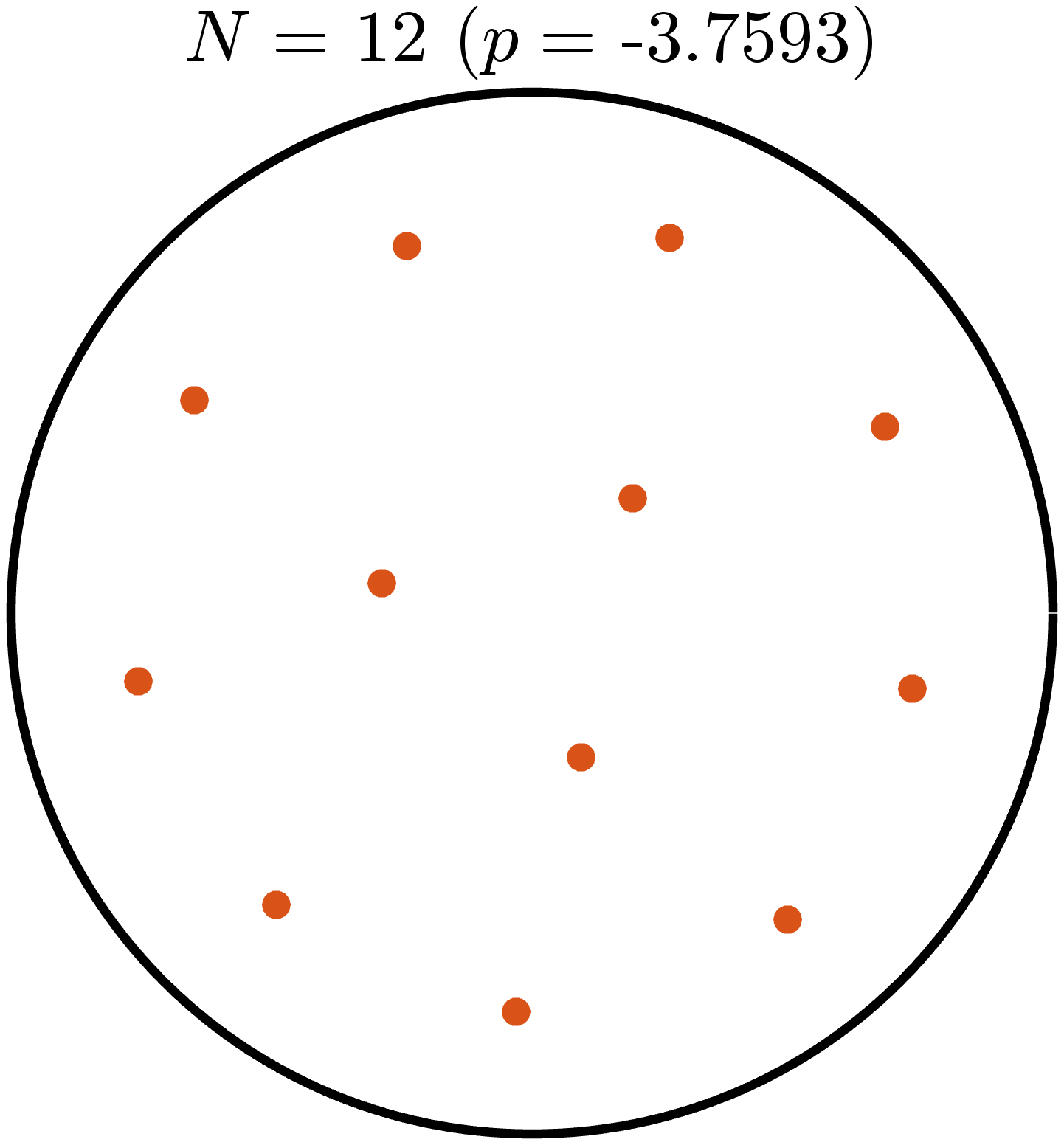}
\includegraphics[width=0.19\textwidth]{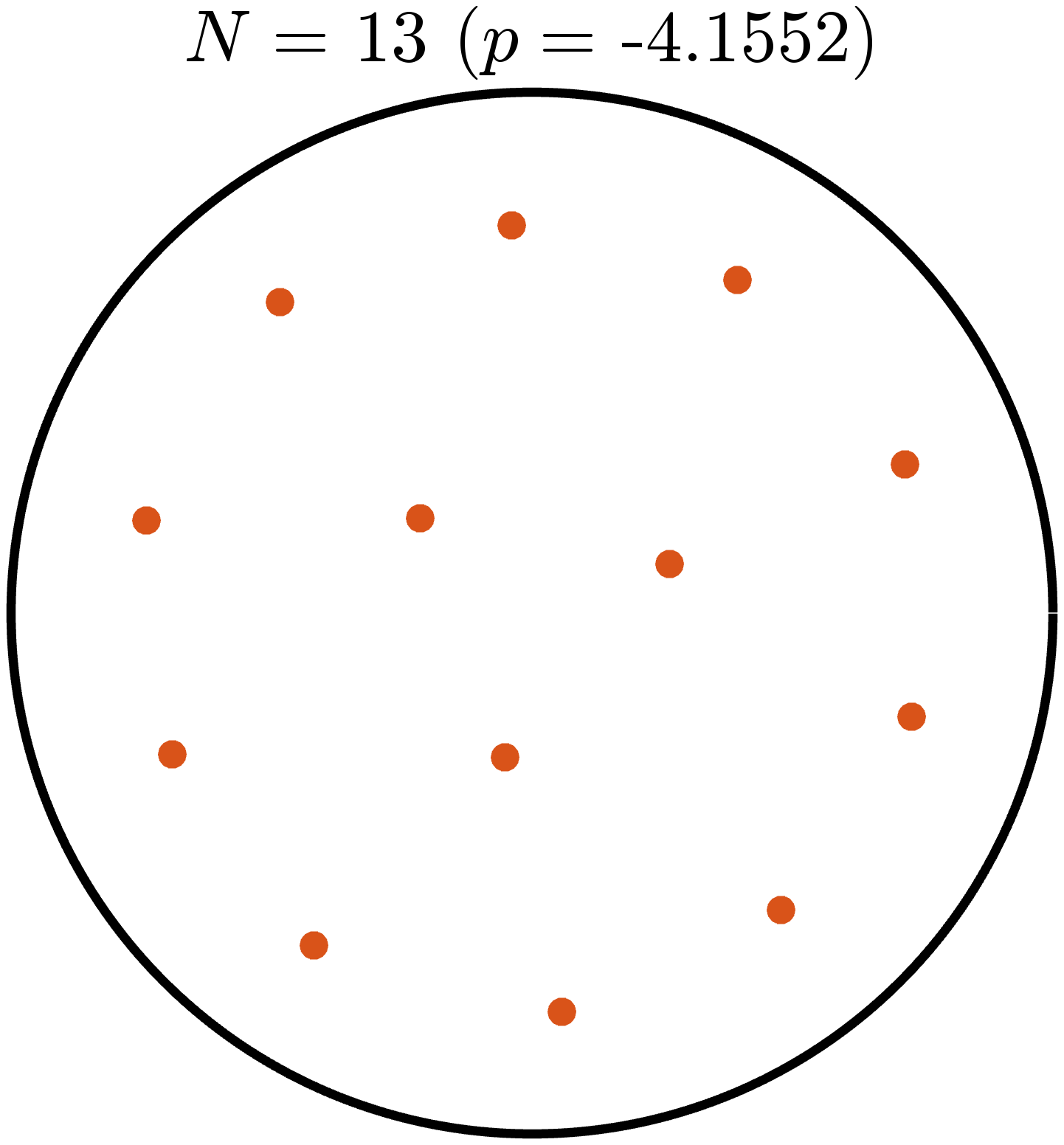}
\includegraphics[width=0.19\textwidth]{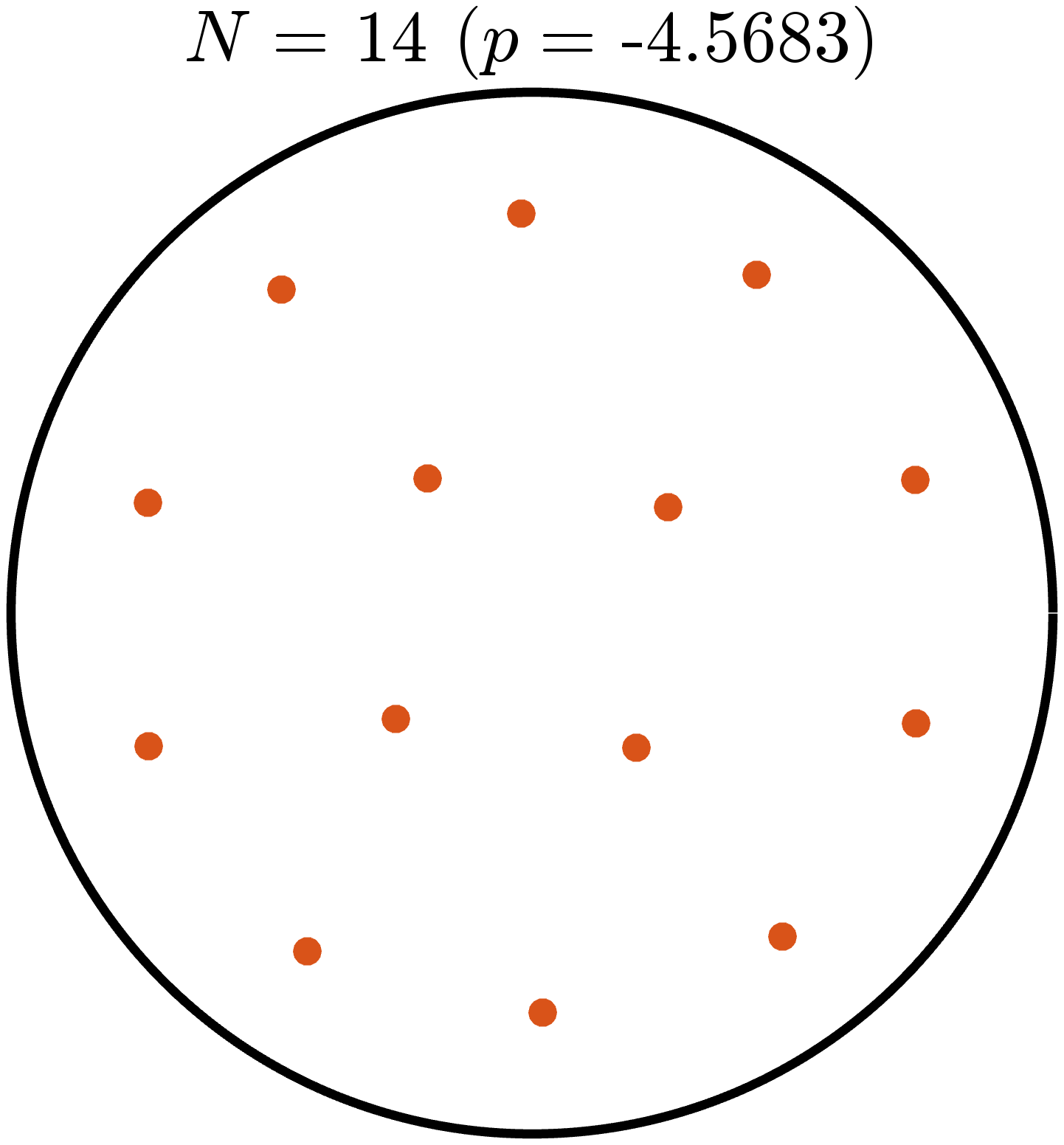}
\includegraphics[width=0.19\textwidth]{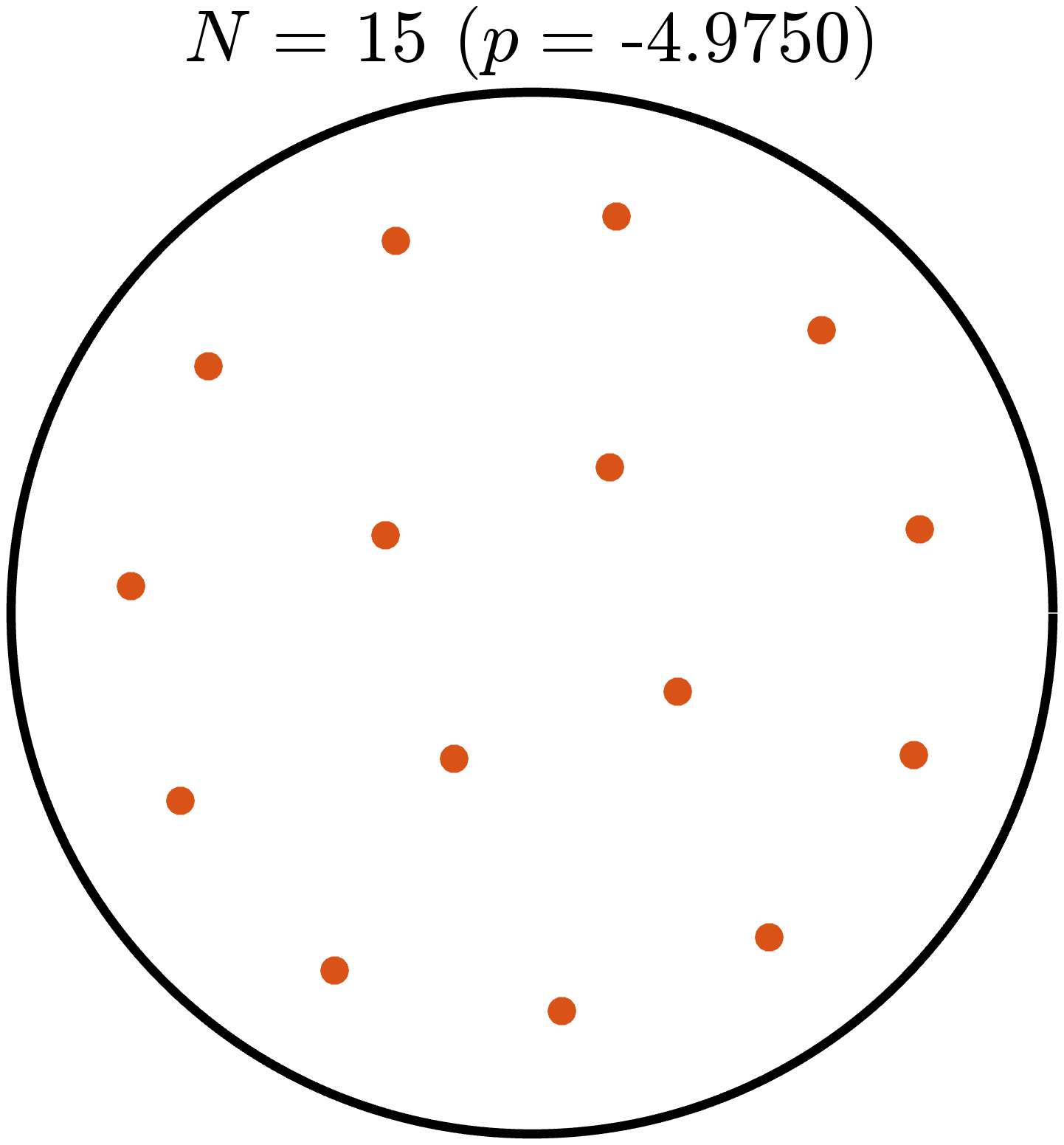}\\[5pt]
\includegraphics[width=0.19\textwidth]{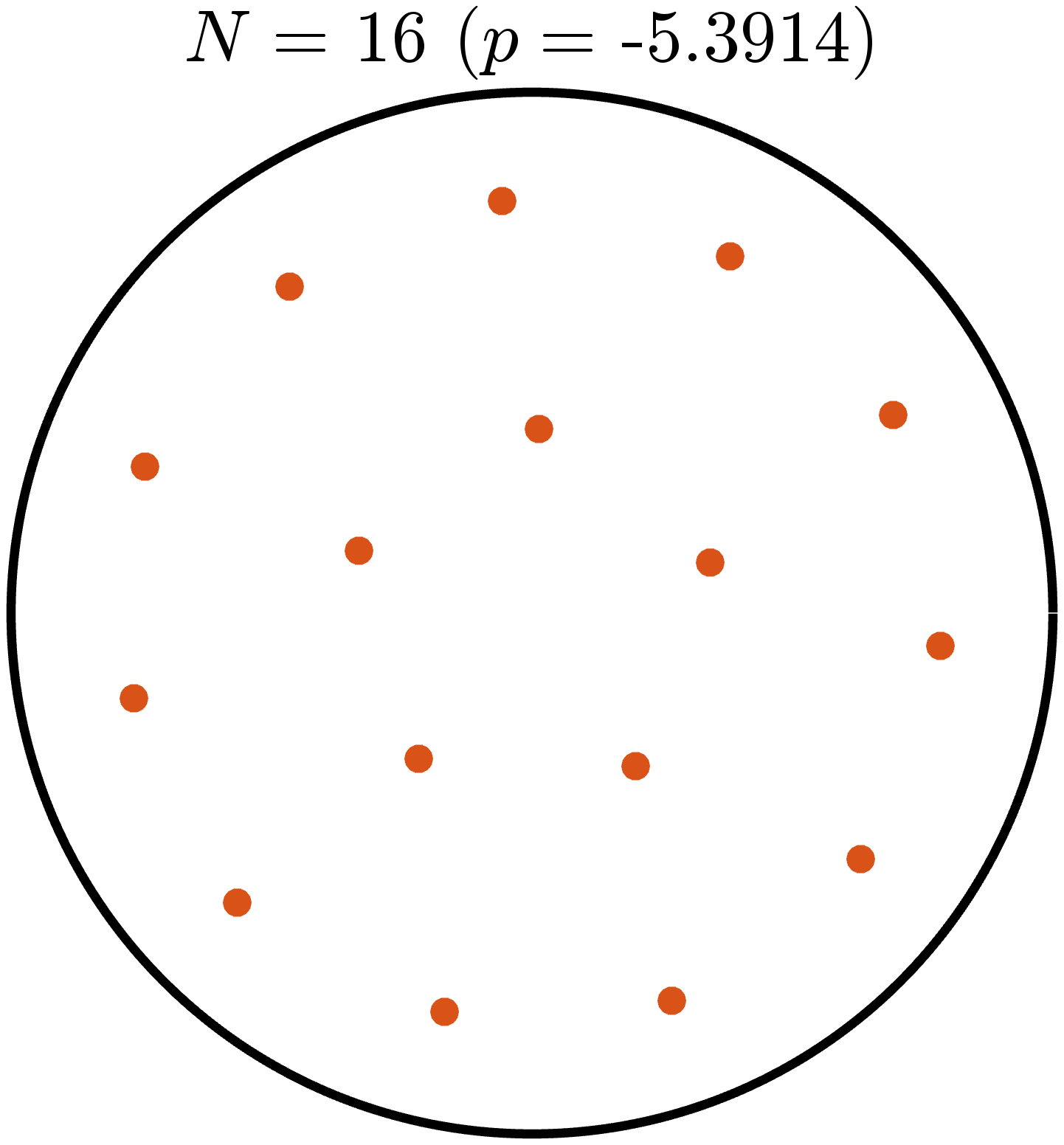}
\includegraphics[width=0.19\textwidth]{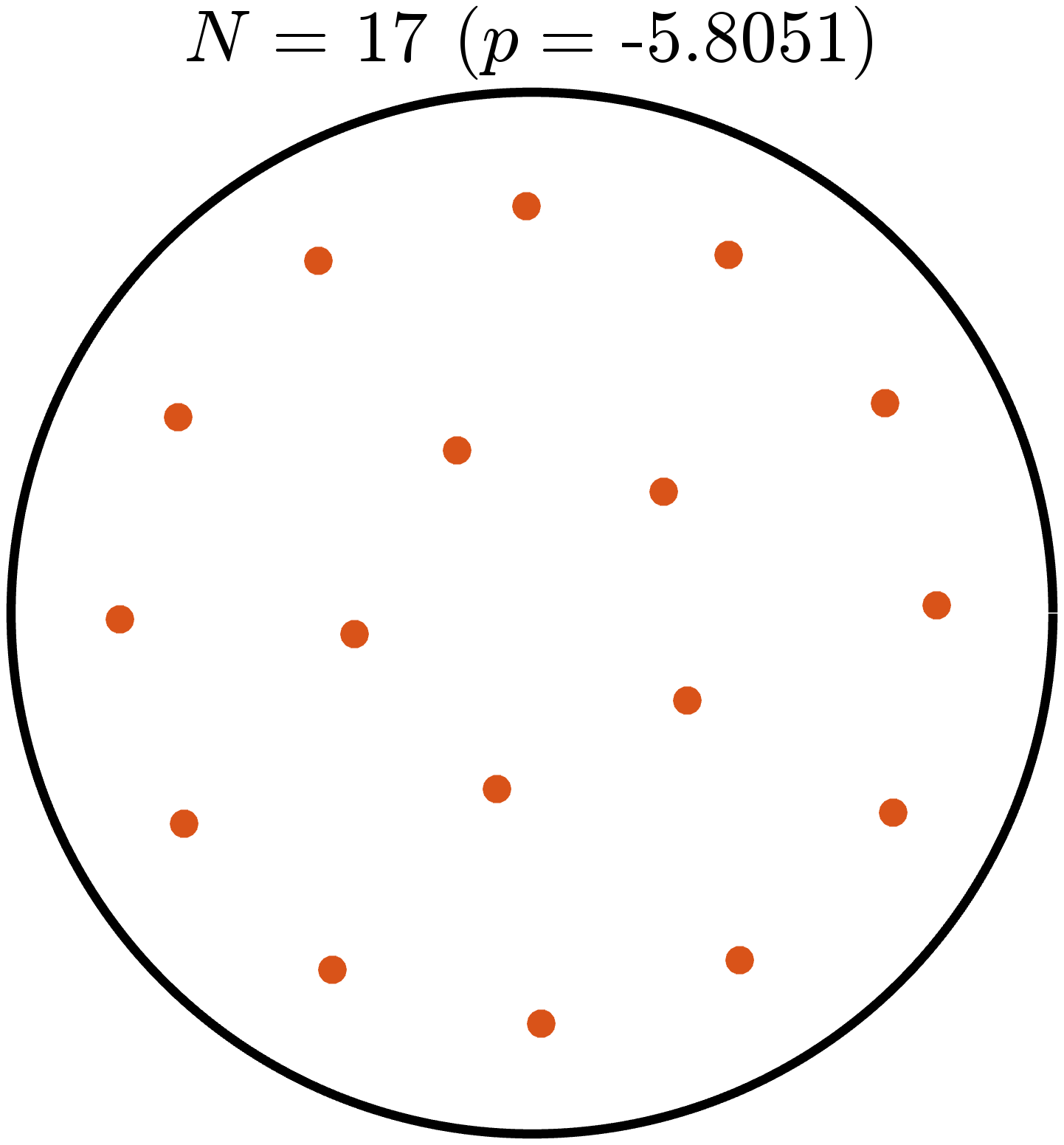}
\includegraphics[width=0.19\textwidth]{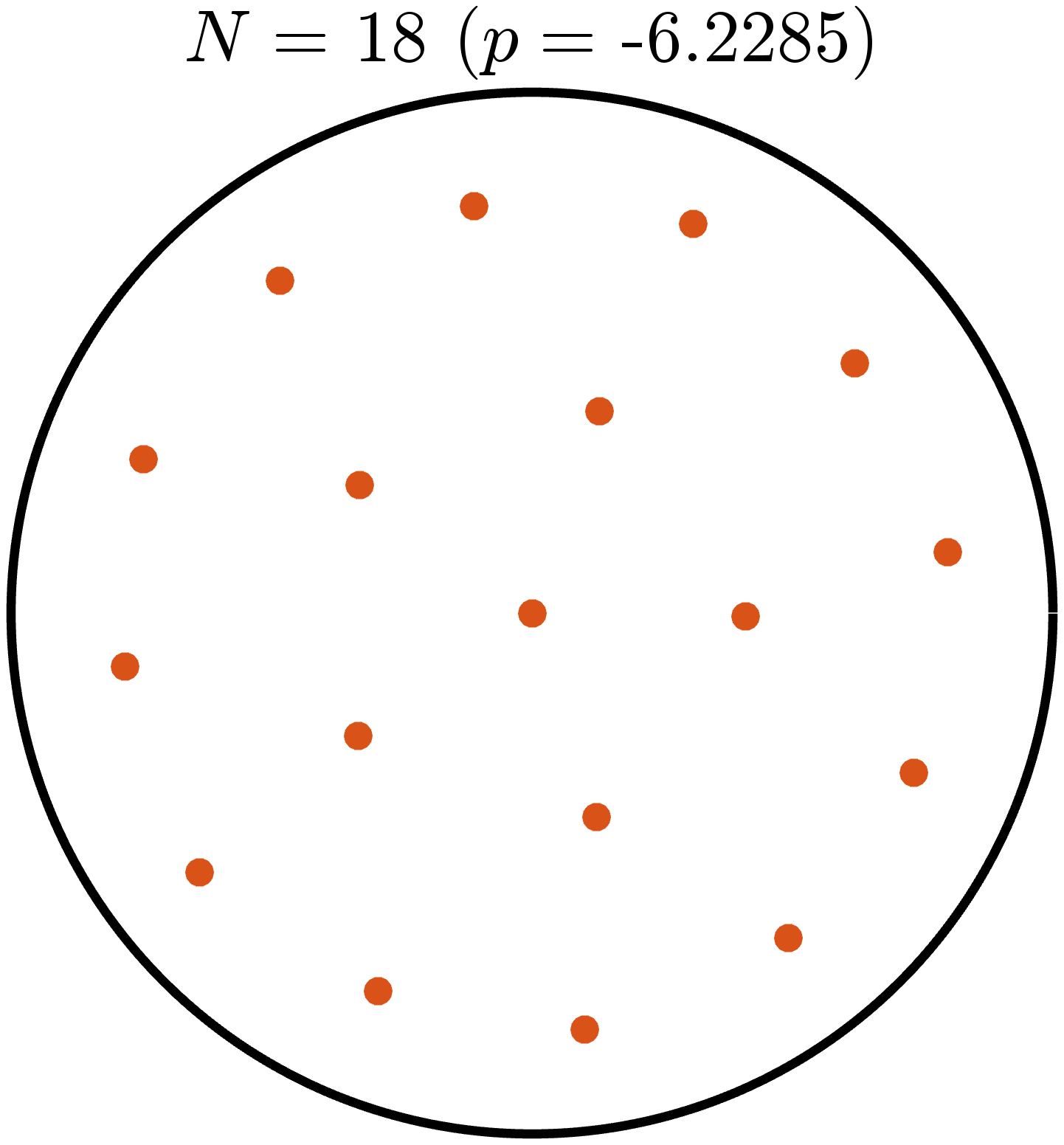}
\includegraphics[width=0.19\textwidth]{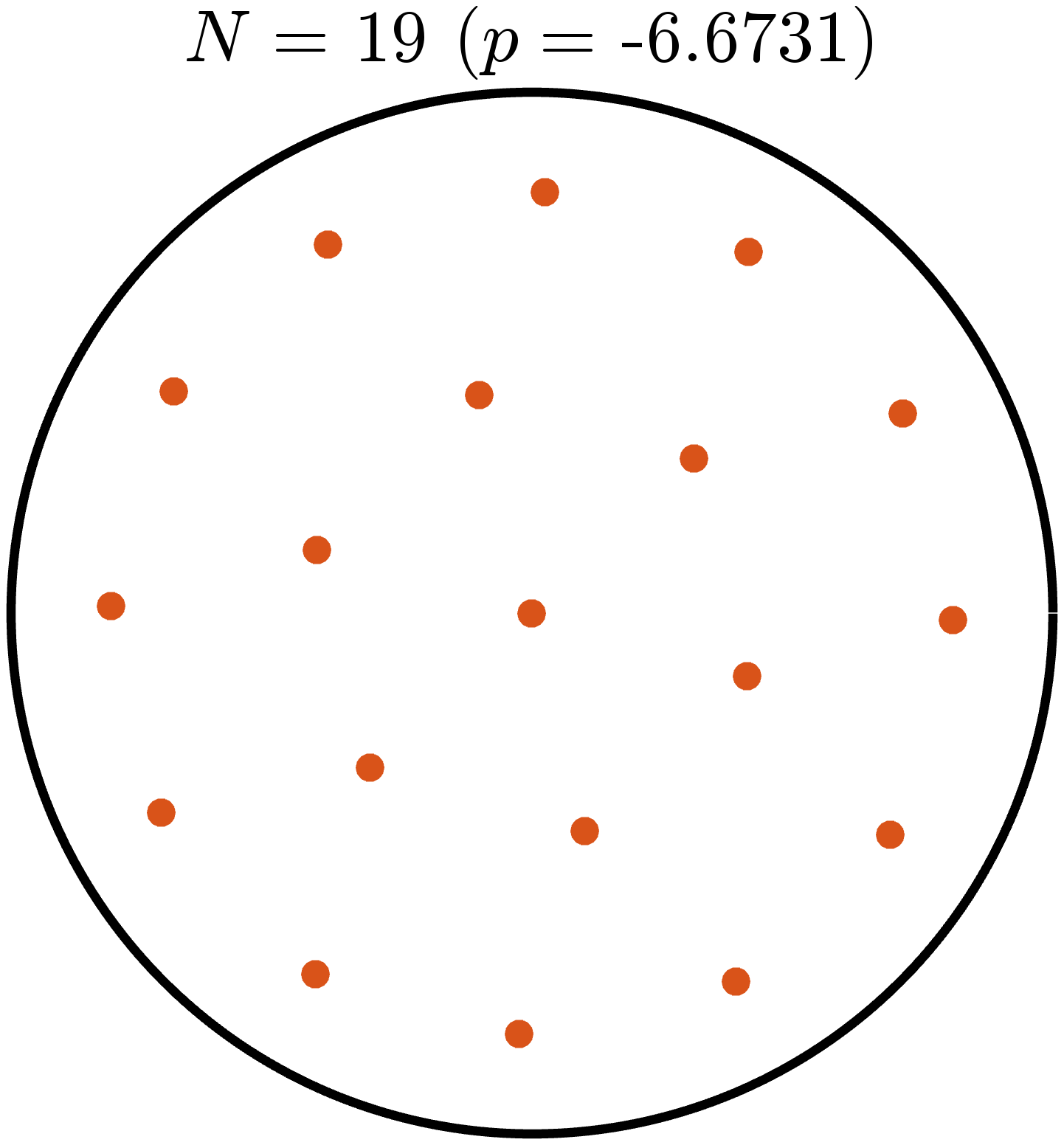}
\includegraphics[width=0.19\textwidth]{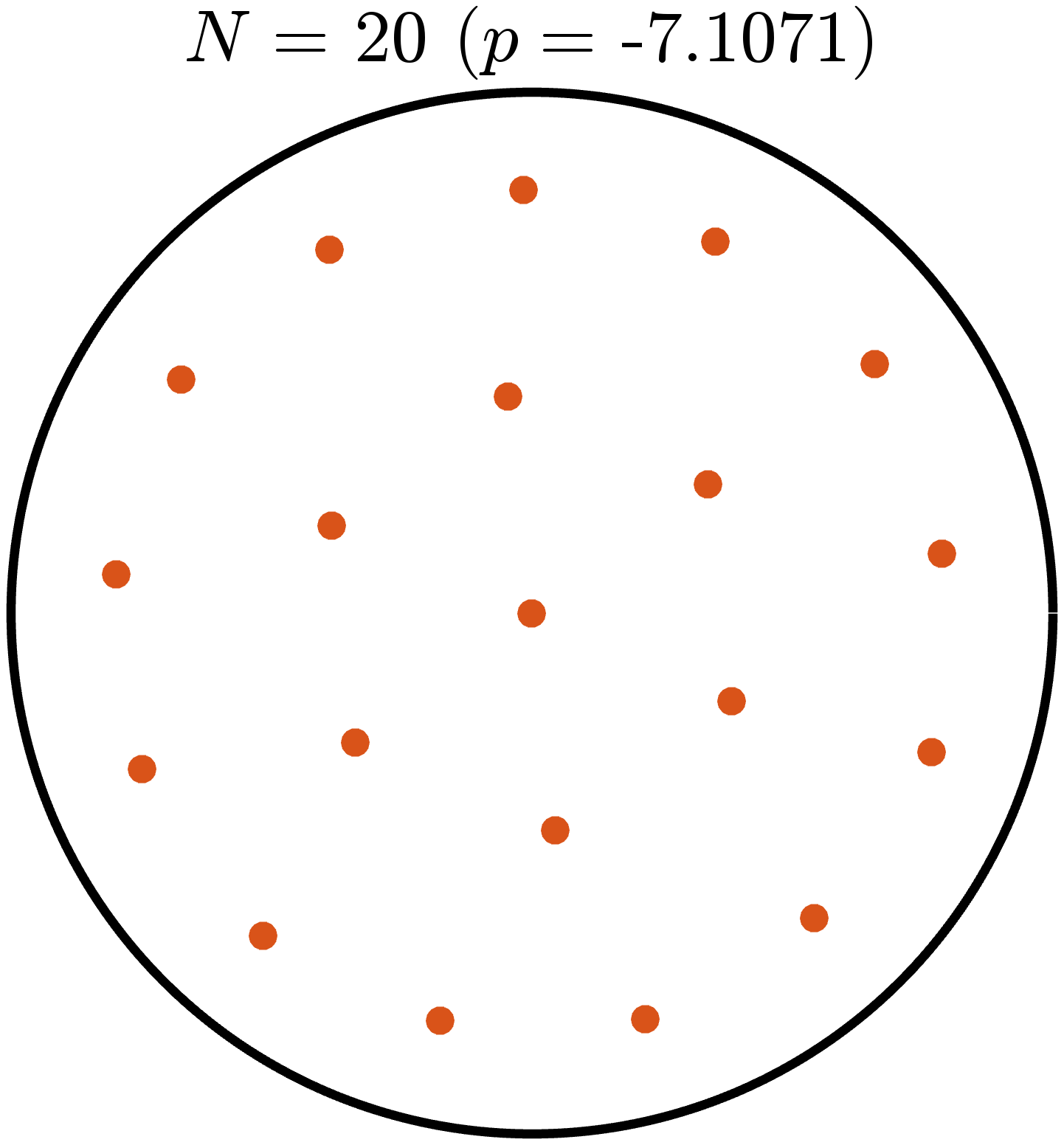}\\[5pt]
\includegraphics[width=0.19\textwidth]{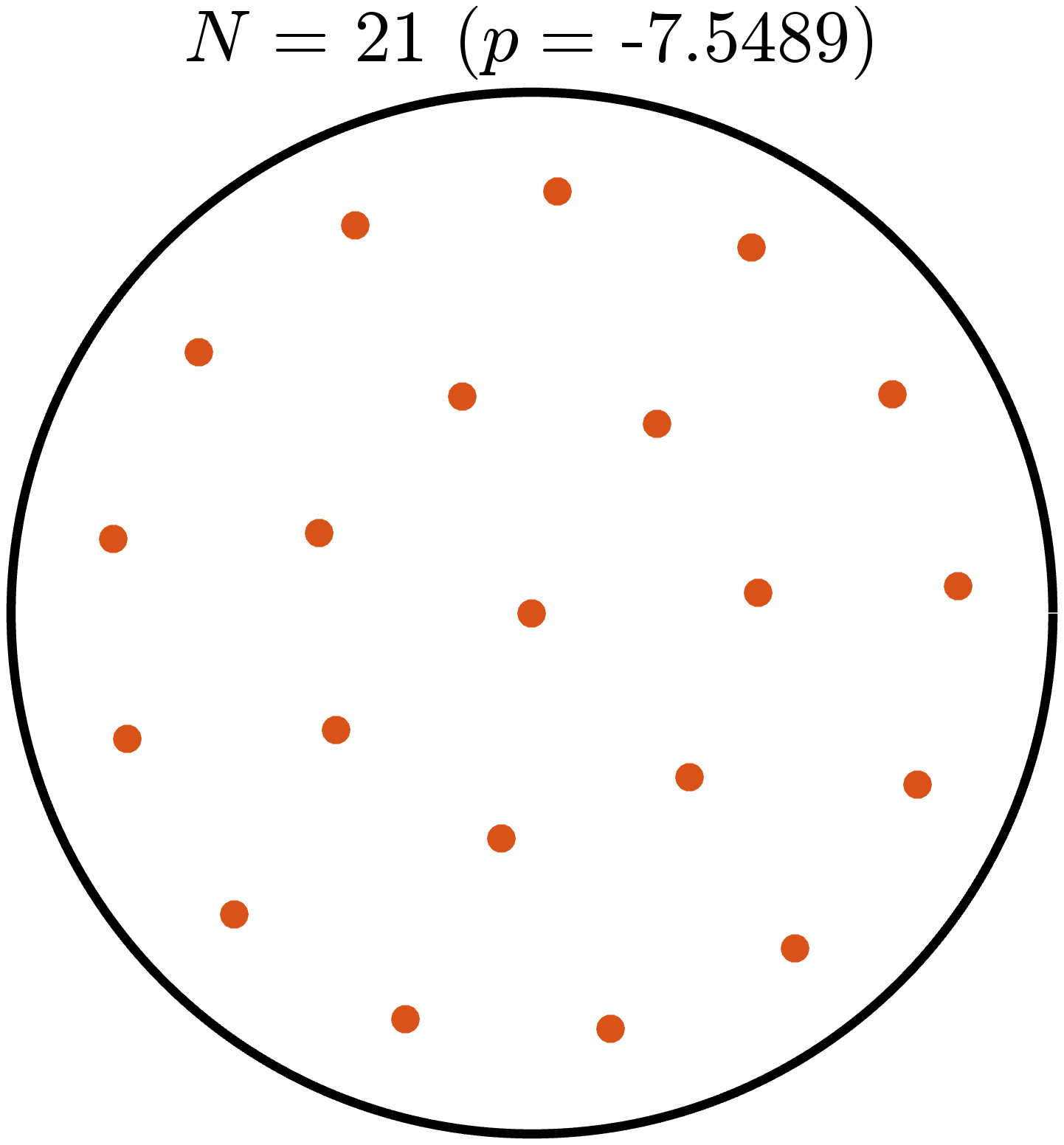}
\includegraphics[width=0.19\textwidth]{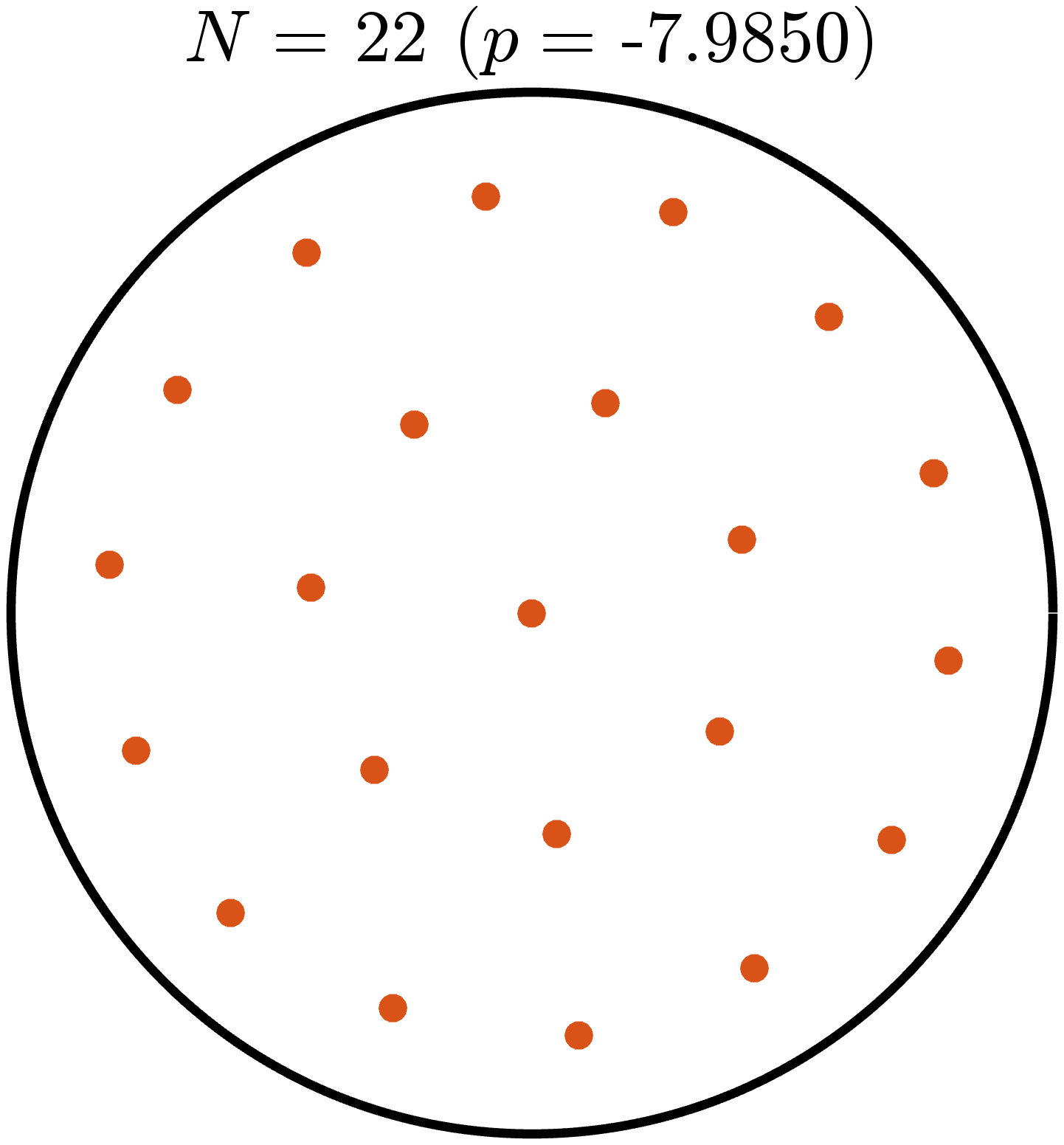}
\includegraphics[width=0.19\textwidth]{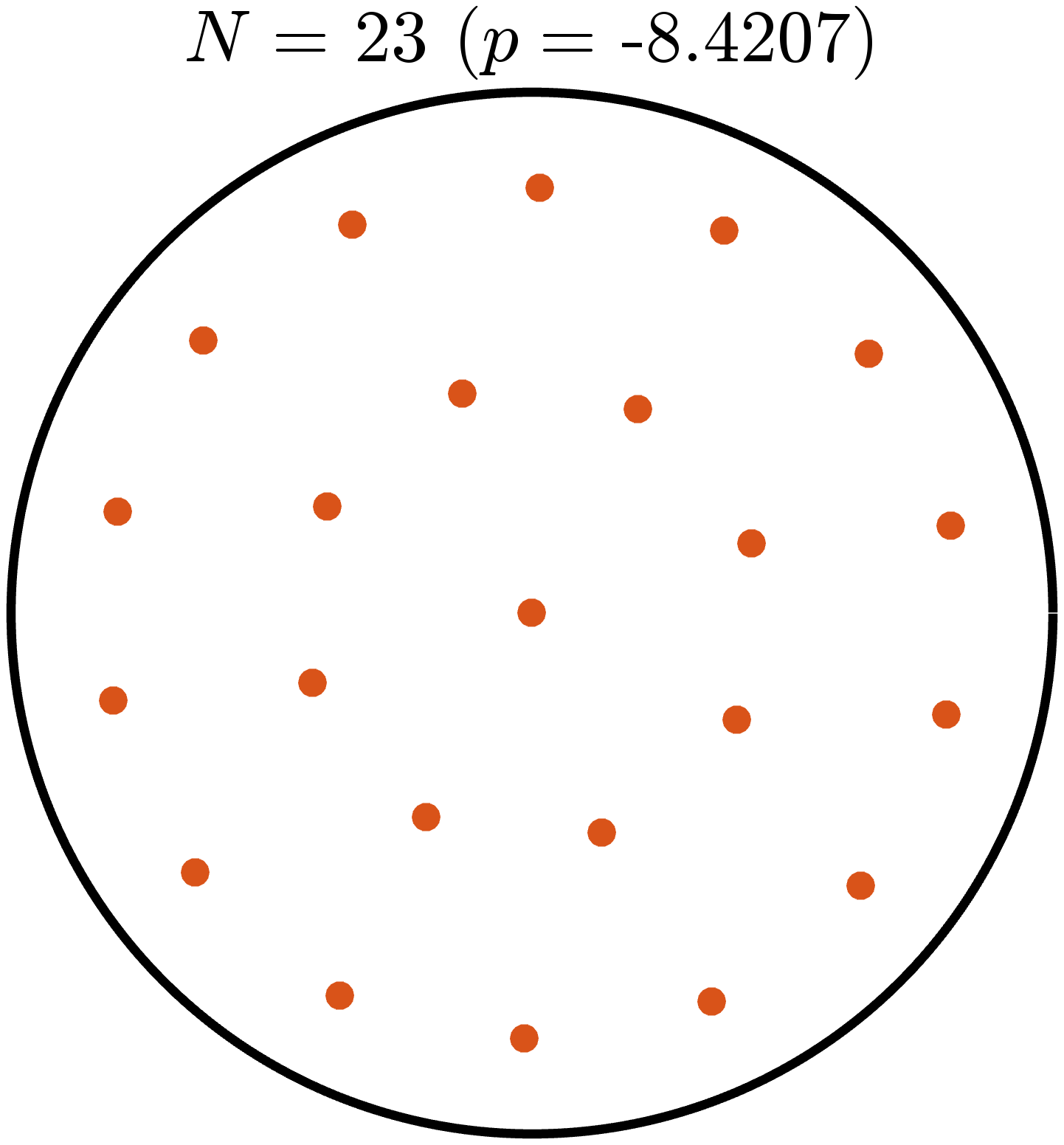}
\includegraphics[width=0.19\textwidth]{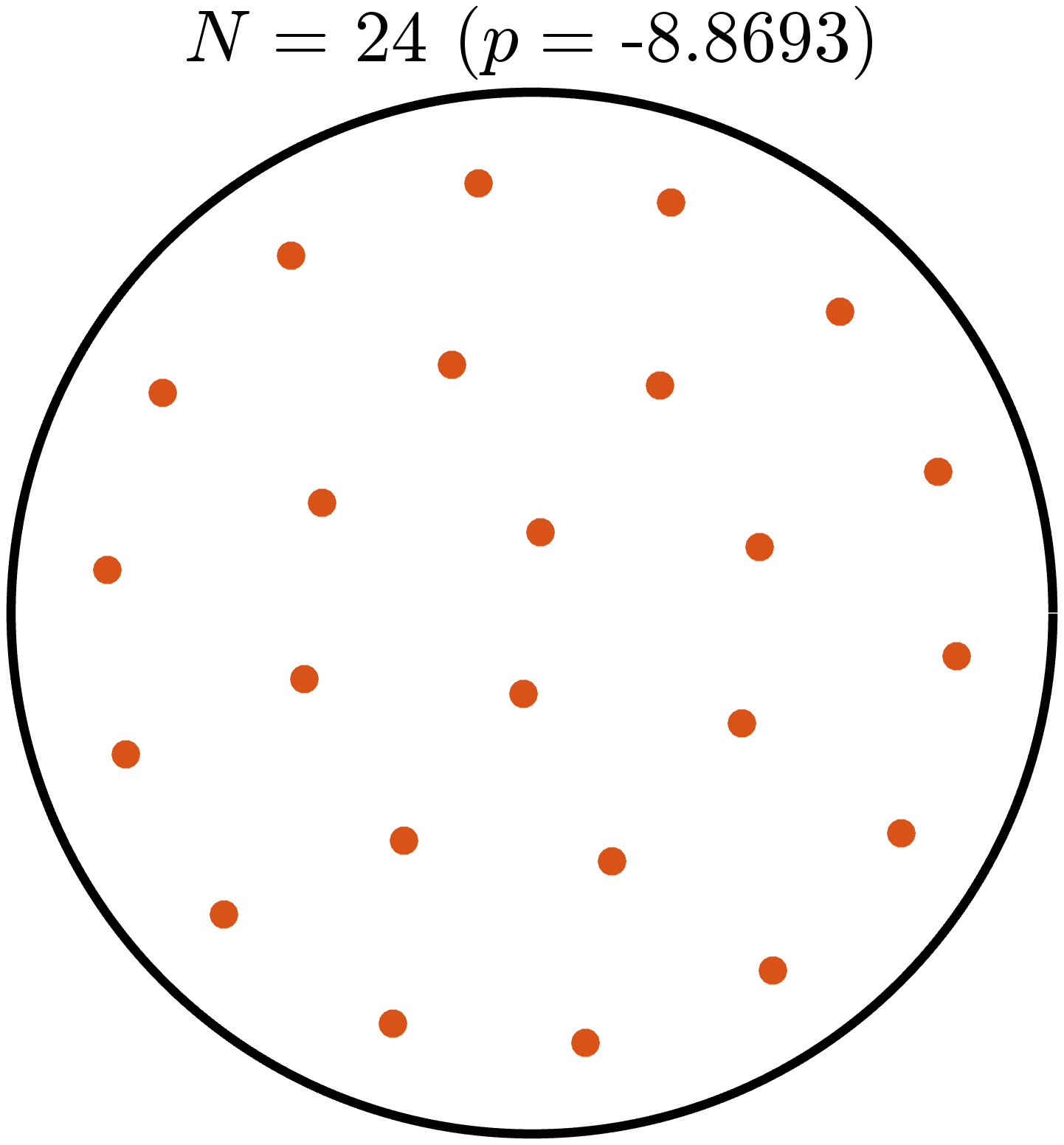}
\includegraphics[width=0.19\textwidth]{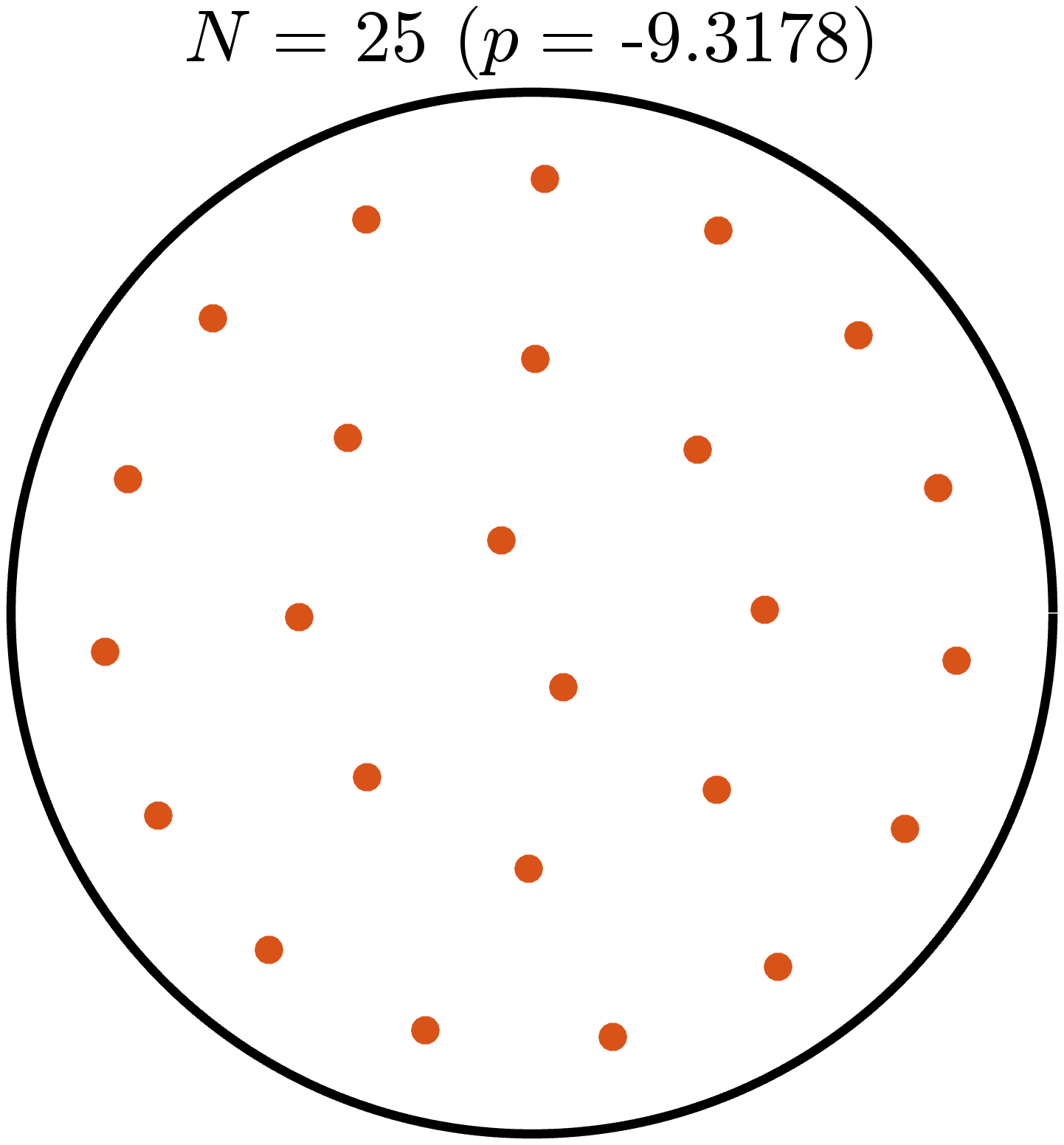}
    \caption{Locally optimal trap configurations for trap centers in the disk for $N= 6,\ldots,25$. This figure reproduces the results of \cite[Fig.~9]{KTW_2005} and provides a highly non-trivial validation of the integrated numerical method for $\Gint$ and optimization routines. As mentioned in \cite{KTW_2005}, some of these patterns are not global minima. \label{fig:disk_optimal}}
\end{figure}

\subsection{Capture at interior traps.}

In this section, we investigate the optimal placement of internal traps for the disk, ellipse, and a random domain.

In the disk case, the exact solution \eqref{eq:GreensDisk} for $\Rint_b$ has been used to find minimizing configurations. The disk is the only known case where an exact formula is available (as distinguished from a series solution) and minimizing configurations have been found in the form of constrained concentric rings, potentially with a trap at the origin. 

In Fig.~\ref{fig:disk_optimal}, we reproduce the result of \cite[Fig.~9]{KTW_2005} which used numerical optimization to determine the optimal configuration for $N=6,\ldots,25$ internal traps. These results, based on optimization over the $2N$ coordinate variables, provide slightly more optimal values than those constrained to lie on concentric circles.

\begin{figure}
    \centering
    \includegraphics[width=0.2425\textwidth]{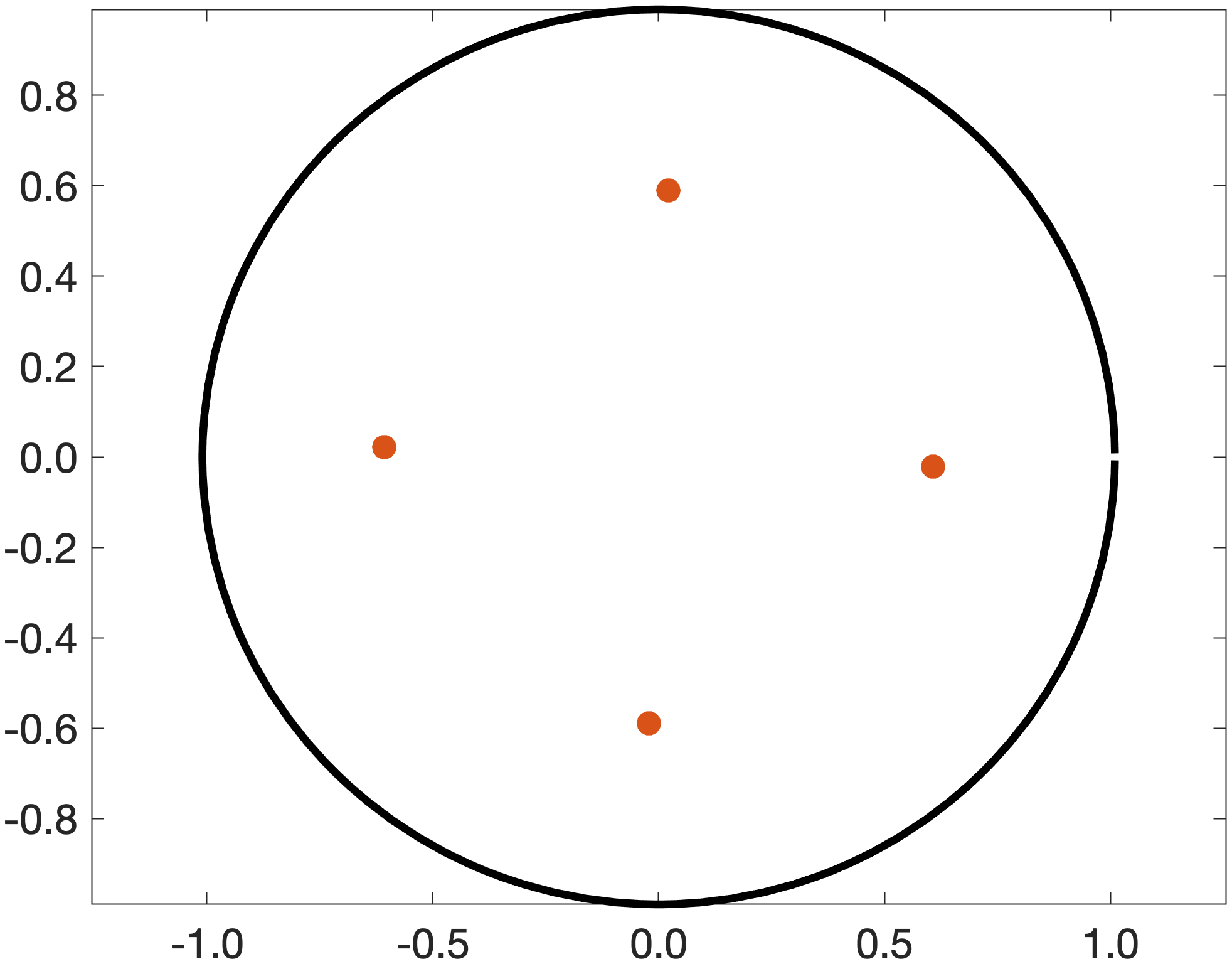}
    \includegraphics[width=0.242\textwidth]{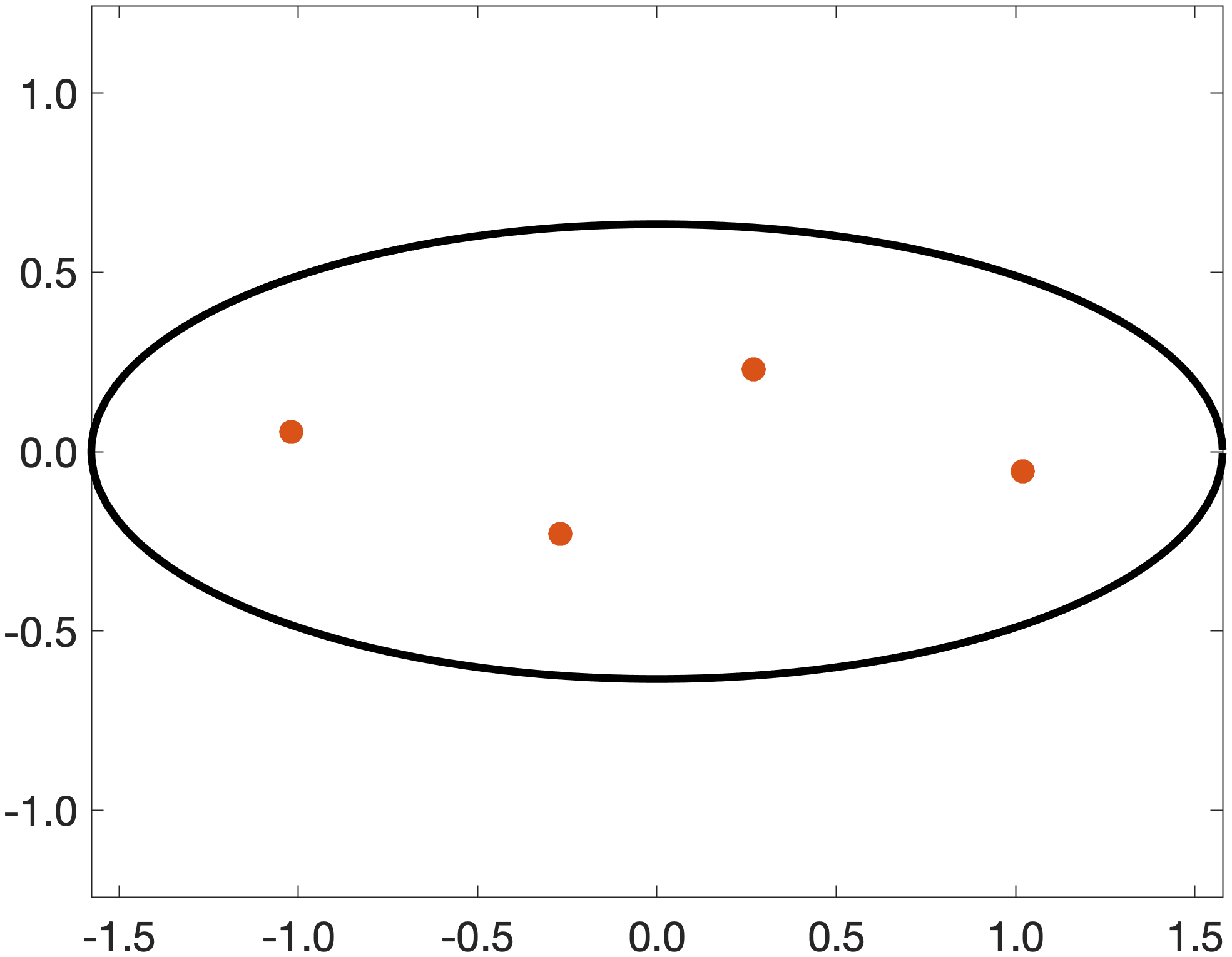}
     \includegraphics[width=0.2425\textwidth]{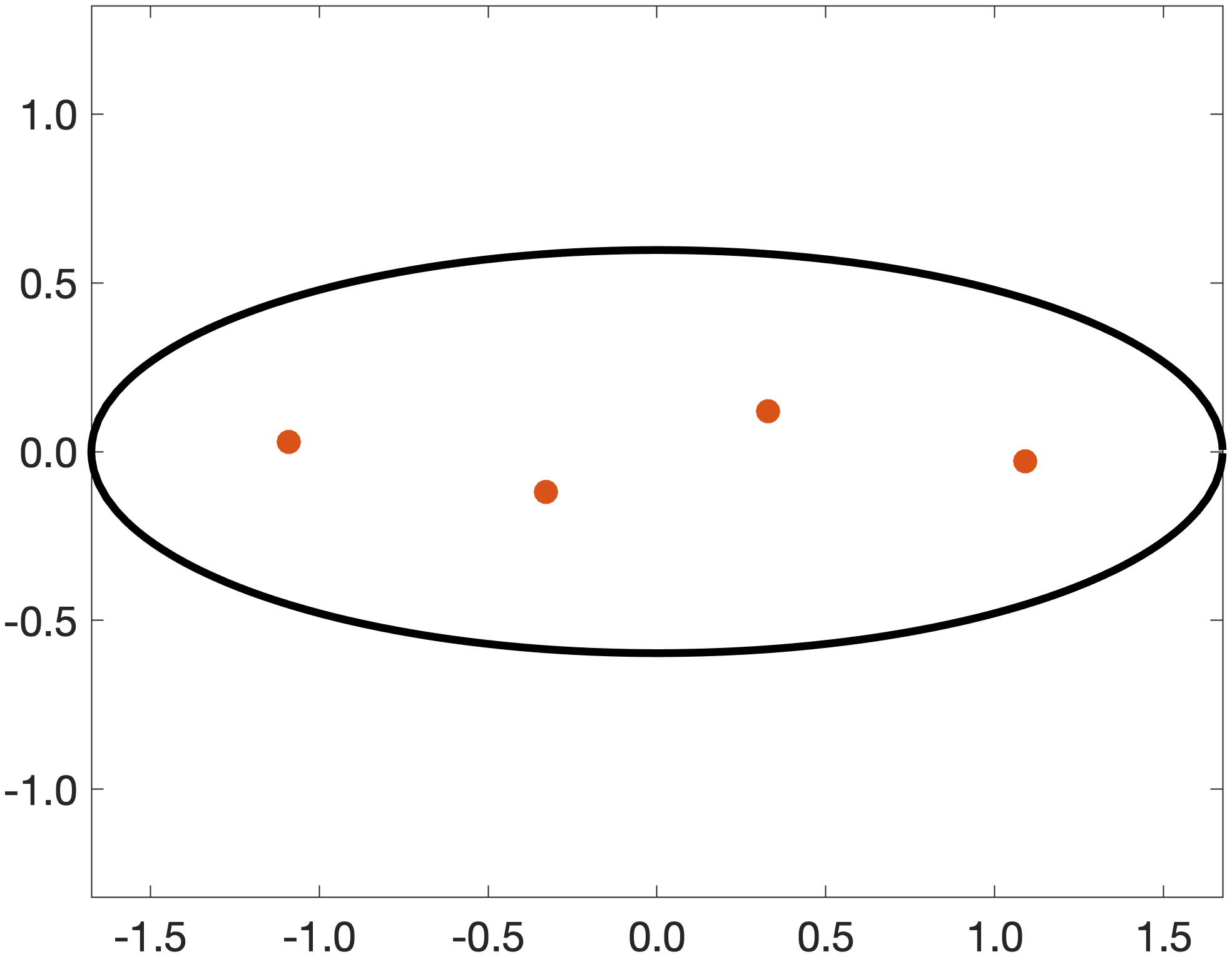}
    \includegraphics[width=0.2425\textwidth]{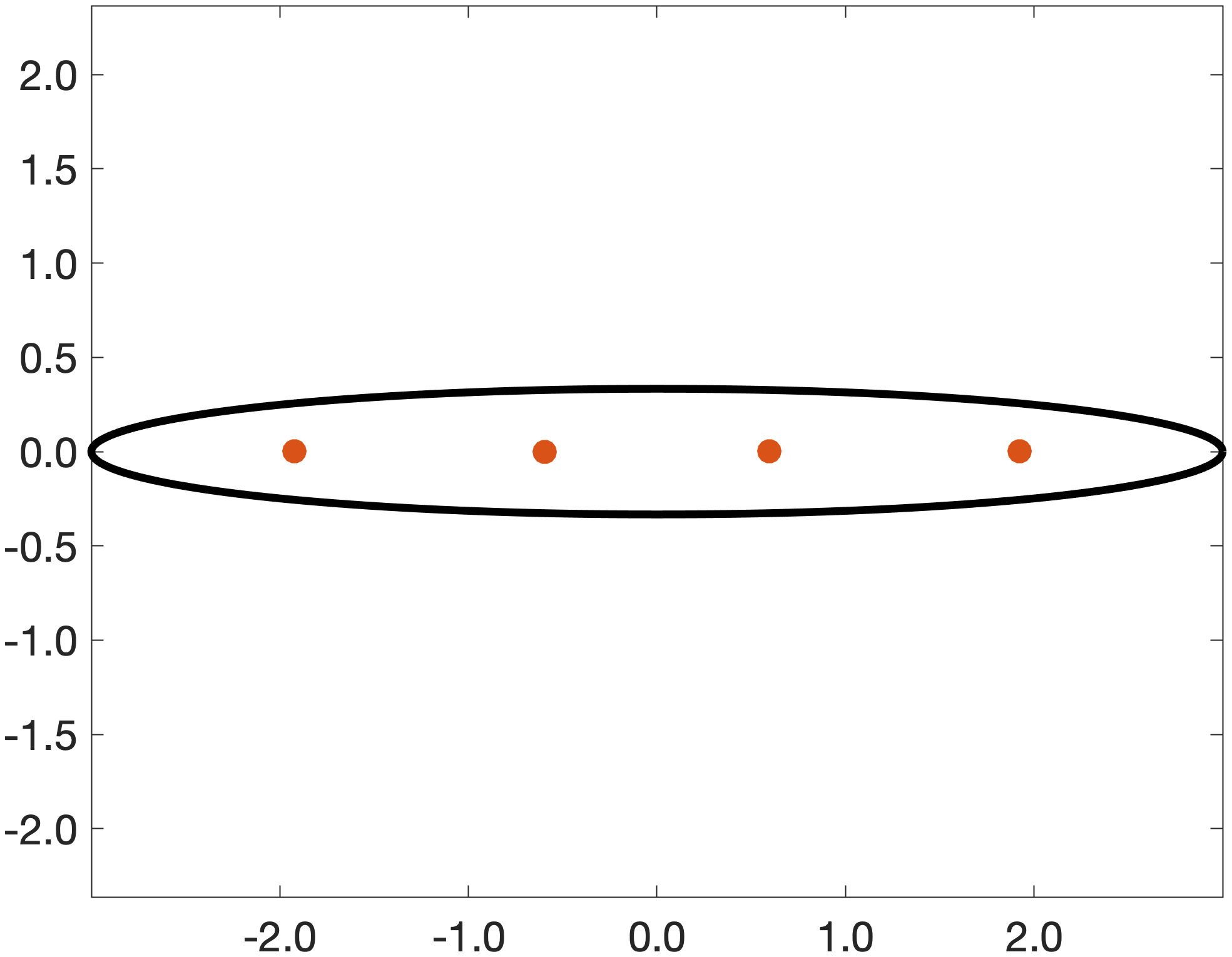}
    \caption{Optimal patterns for $N=4$ traps in an elliptical domain. As a partial confirmation of our process, we have reproduced Fig.~11 of \cite{Sarafa2021} which shows the optimizing configurations for $N=4$ traps in elliptical domains of area $\pi$ major axes  $a=1$ ($b=1$), $a=1.577$ ($b\approx0.634$), $a=1.675$ $ (b\approx 0.597)$, $a=3$ ($b\approx0.333$). \label{fig:EllipseN=4} }
\end{figure}

In \cite{Sarafa2021}, the explicit solution of \eqref{eqn_neumG} in elliptical domains was used to identify a bifurcation of optimal placement away from linear alignment. Specifically, the optimal arrangement of traps corresponds to equally spaced placement along the semi-major axis at small $N$ before bifurcation off the semi-major axis at a critical eccentricity.  In Fig.~\ref{fig:EllipseN=4}, we plot results for $N=4$ traps and a ellipse of area $|\Omega| = \pi ab = \pi$ for various semi-major axis $a$ and $b=1/a$.

\begin{figure}[t]
    \centering
\includegraphics[width=0.19\textwidth]{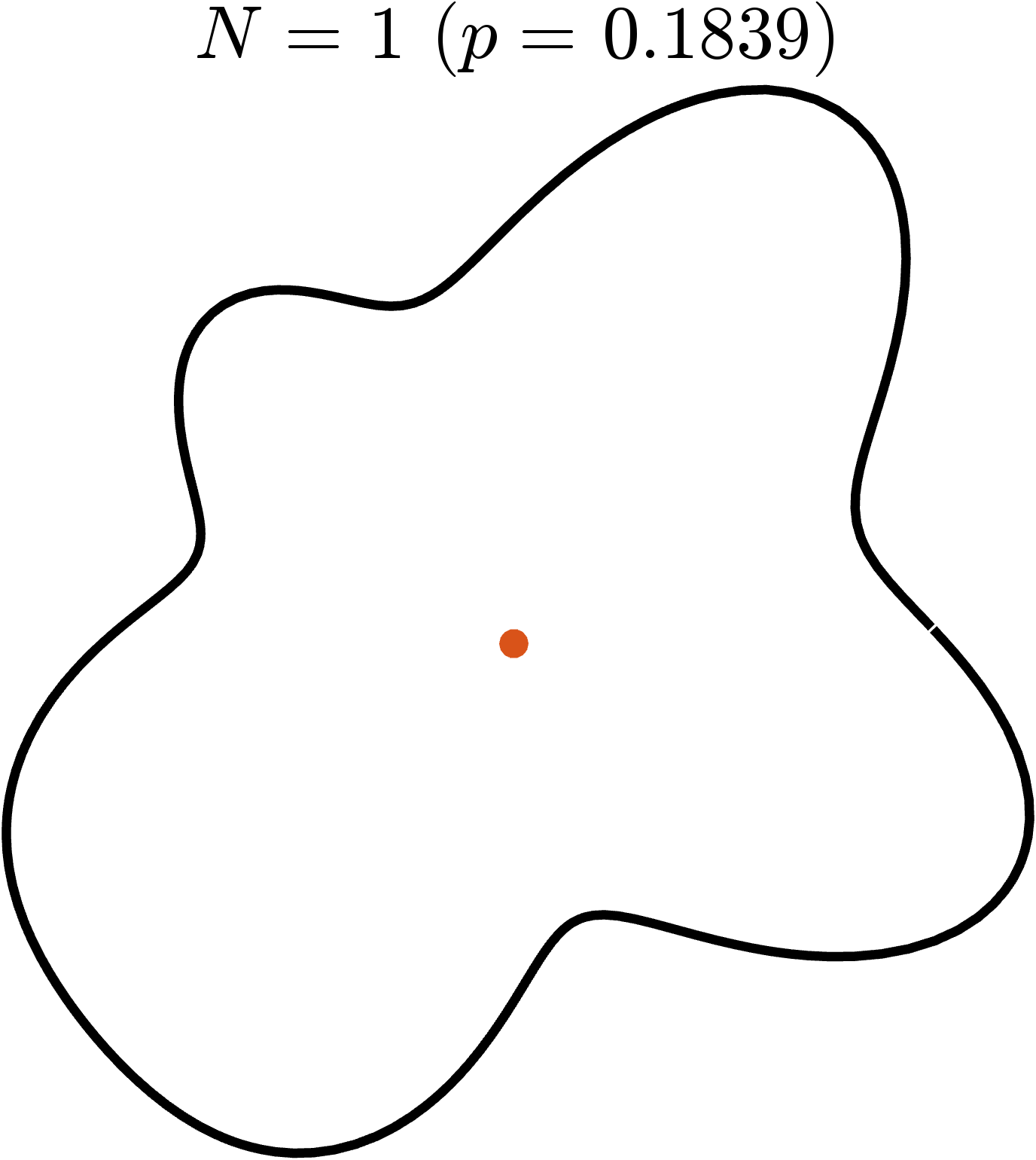}
\includegraphics[width=0.19\textwidth]{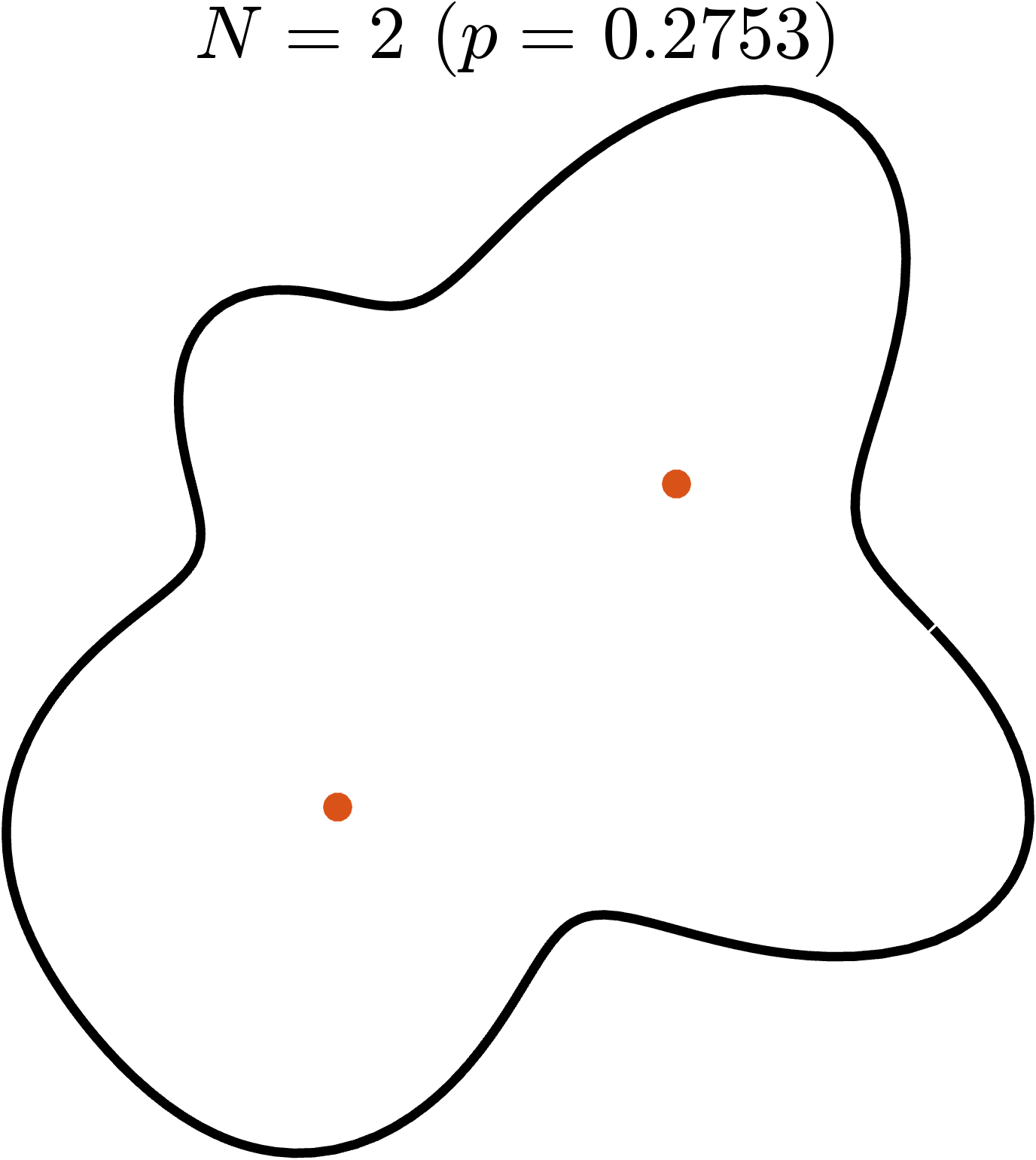}
\includegraphics[width=0.19\textwidth]{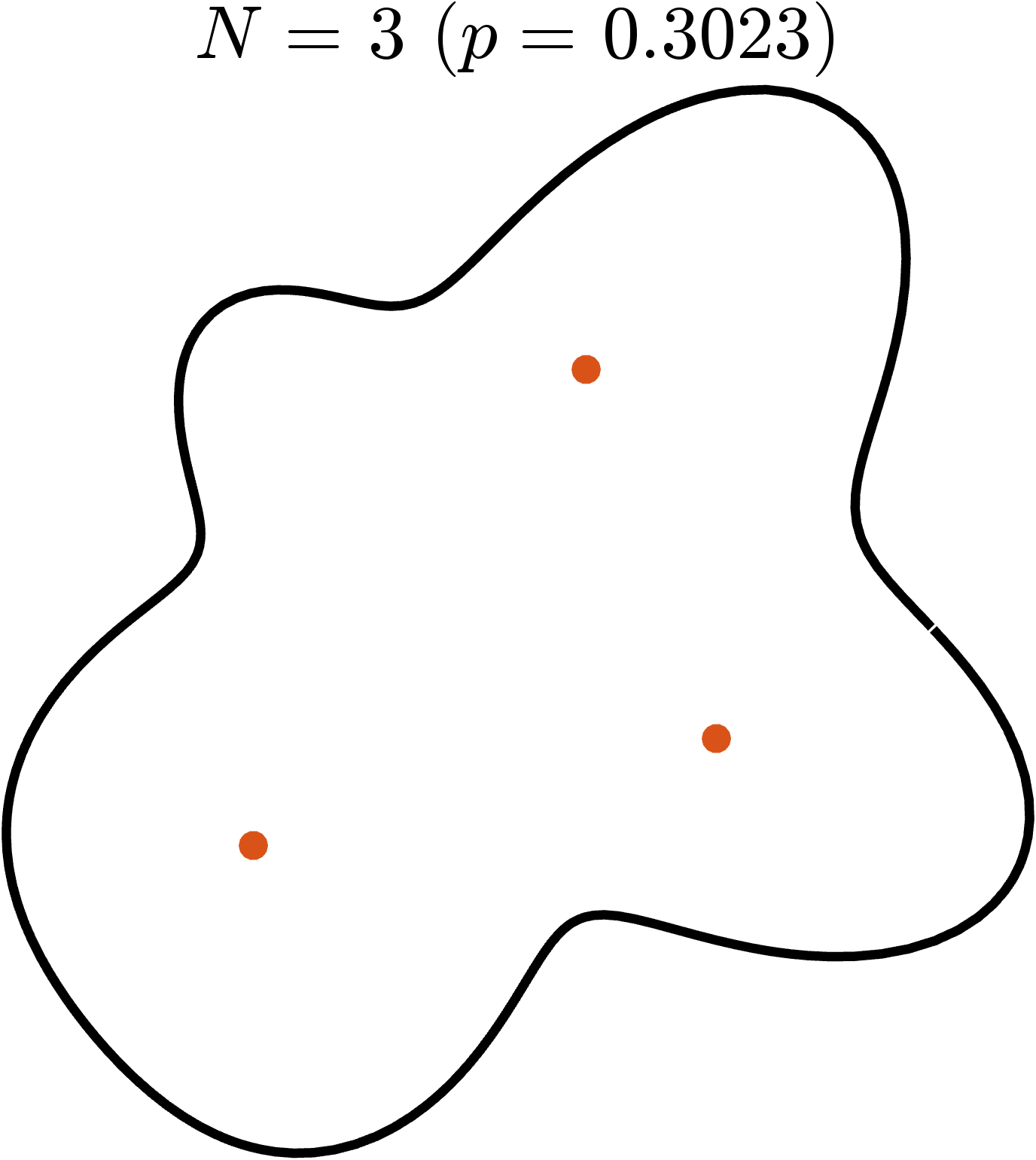}
\includegraphics[width=0.19\textwidth]{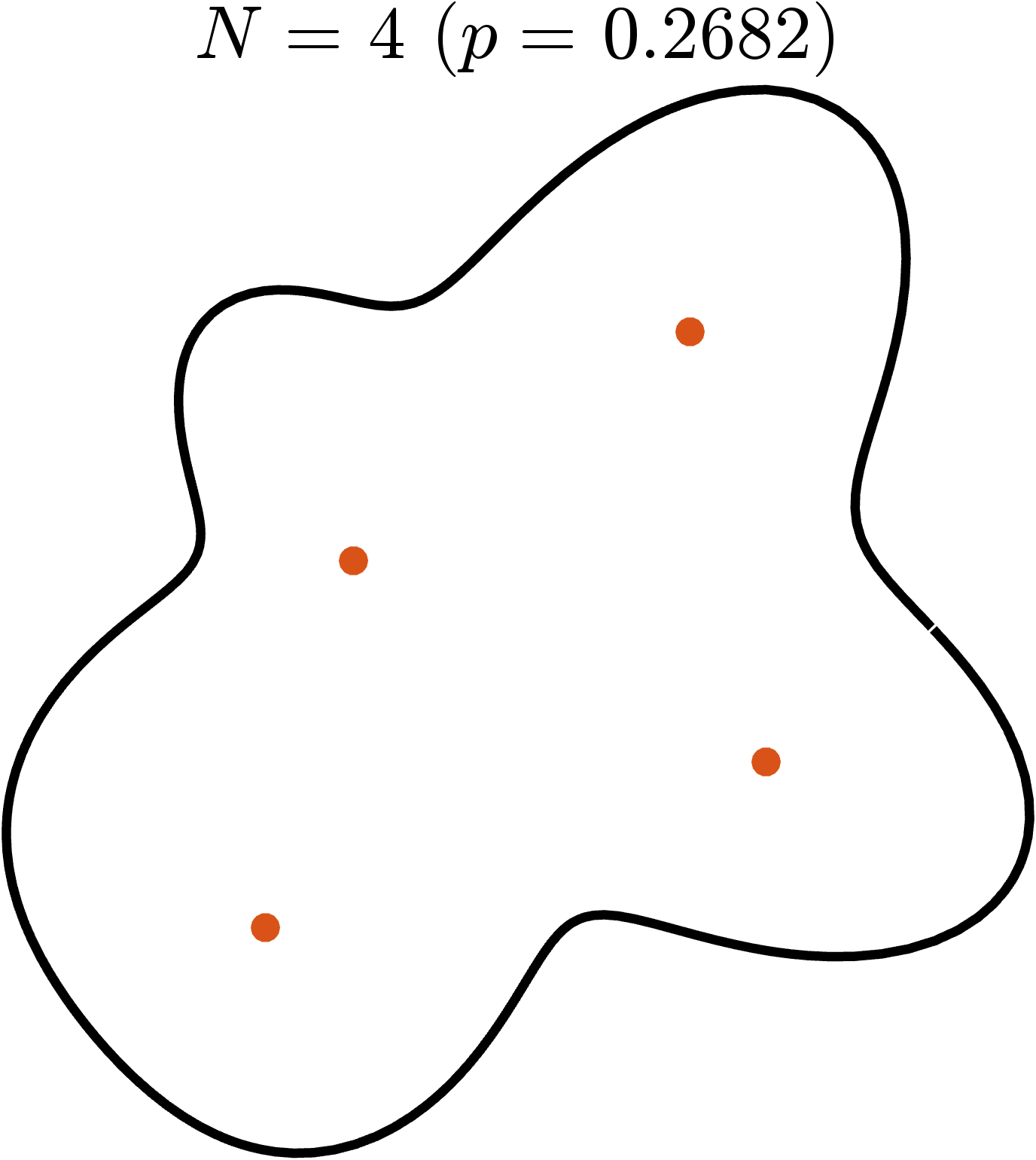}
\includegraphics[width=0.19\textwidth]{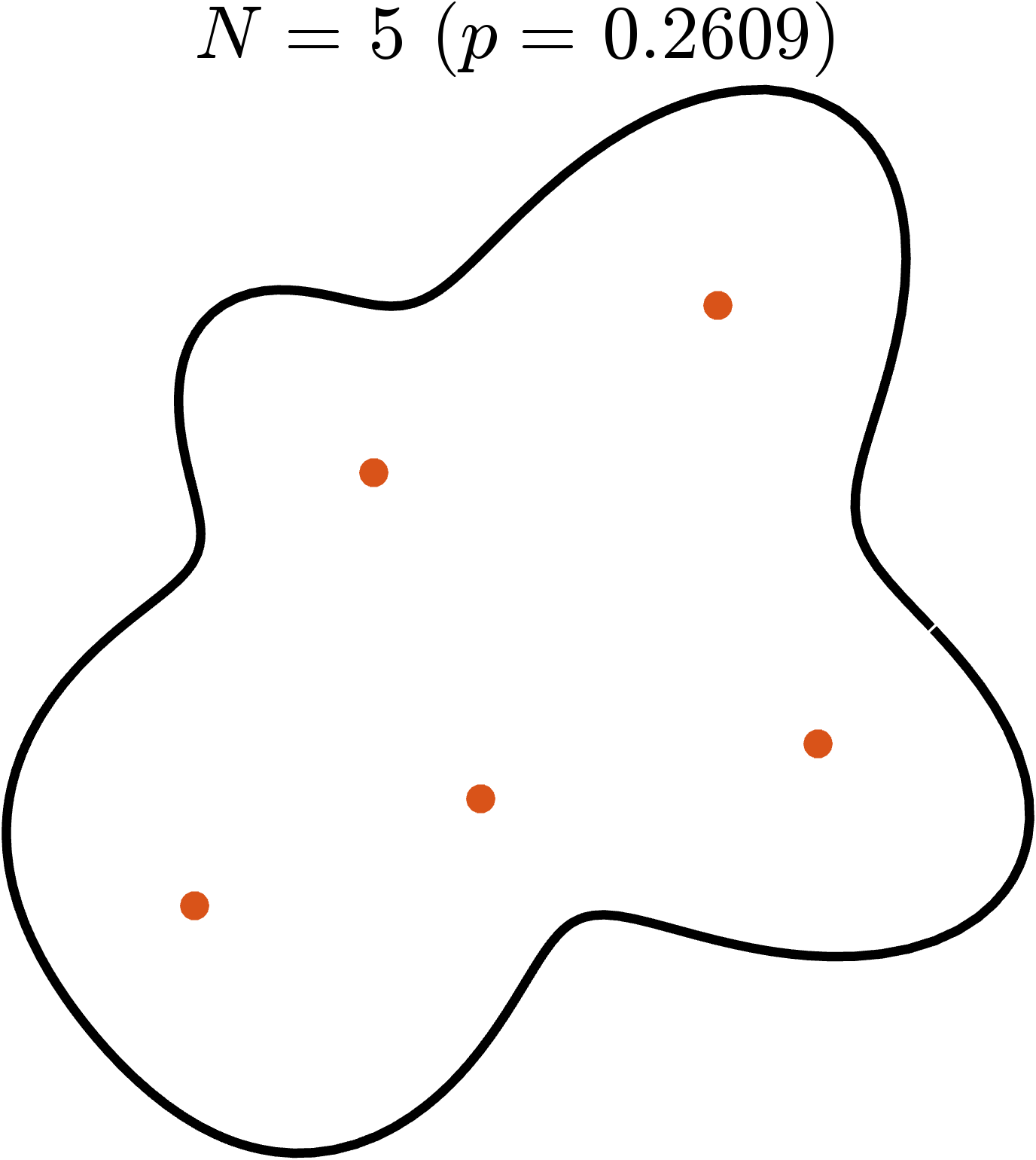}\\[5pt]
\includegraphics[width=0.19\textwidth]{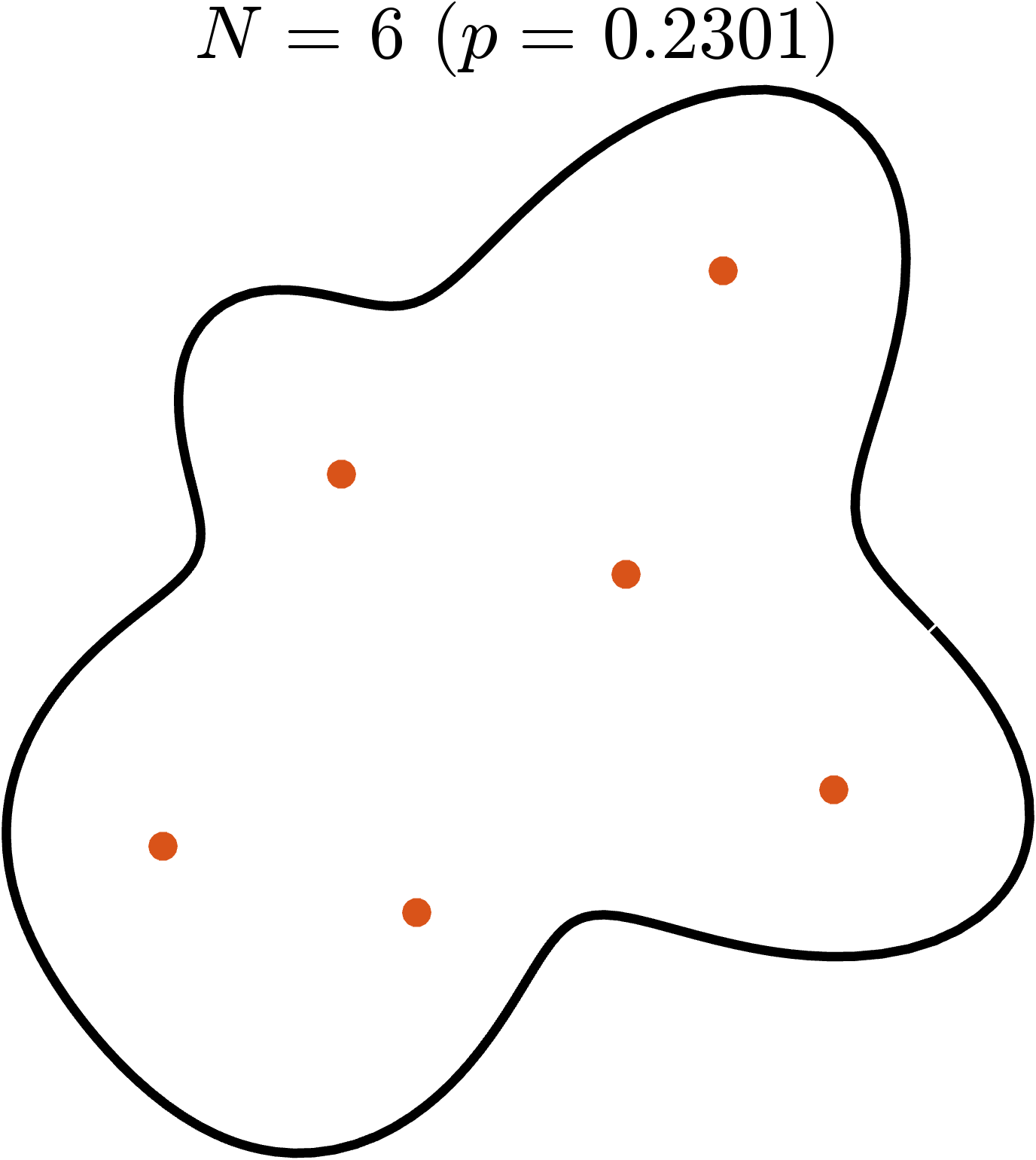}
\includegraphics[width=0.19\textwidth]{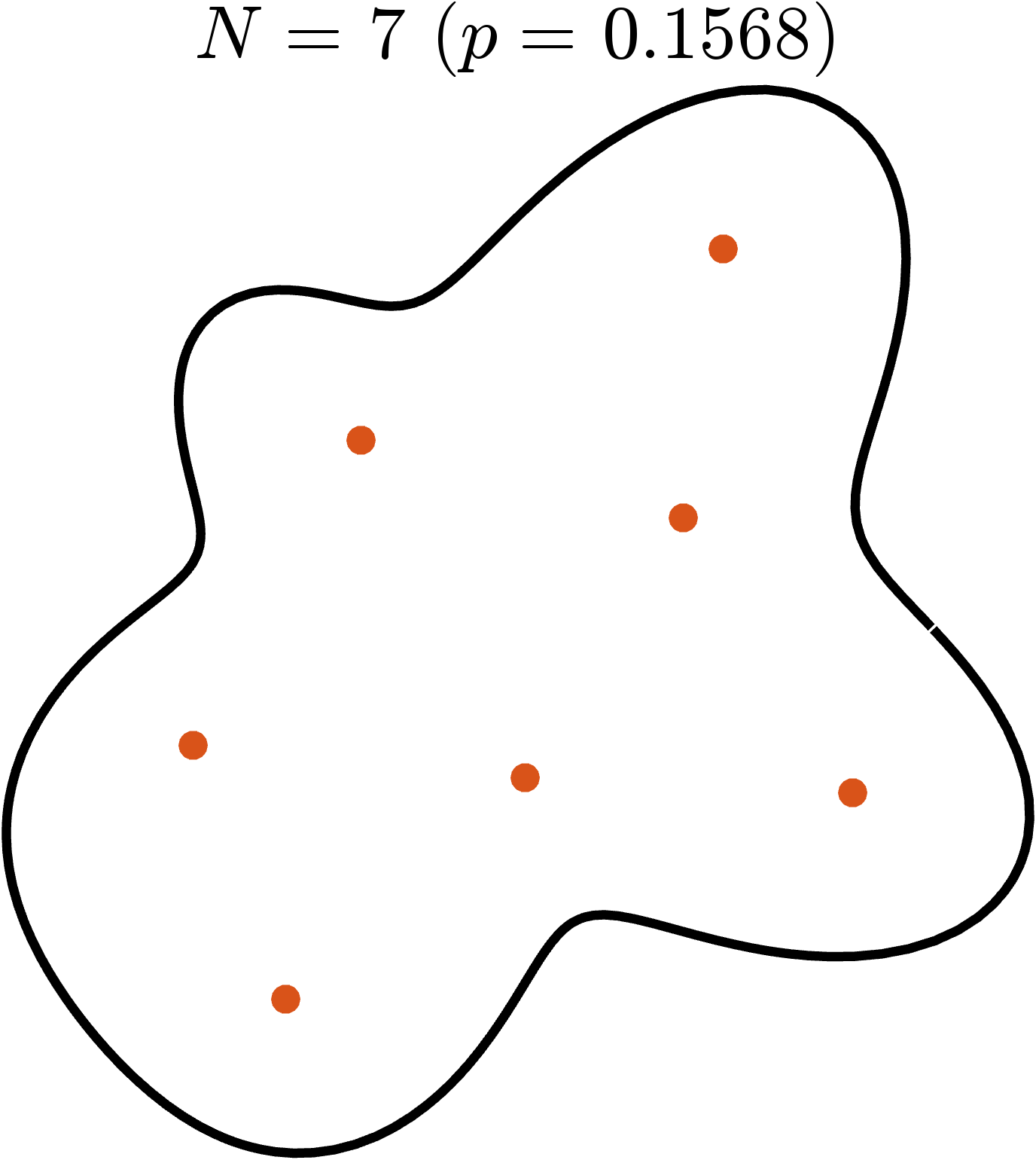}
\includegraphics[width=0.19\textwidth]{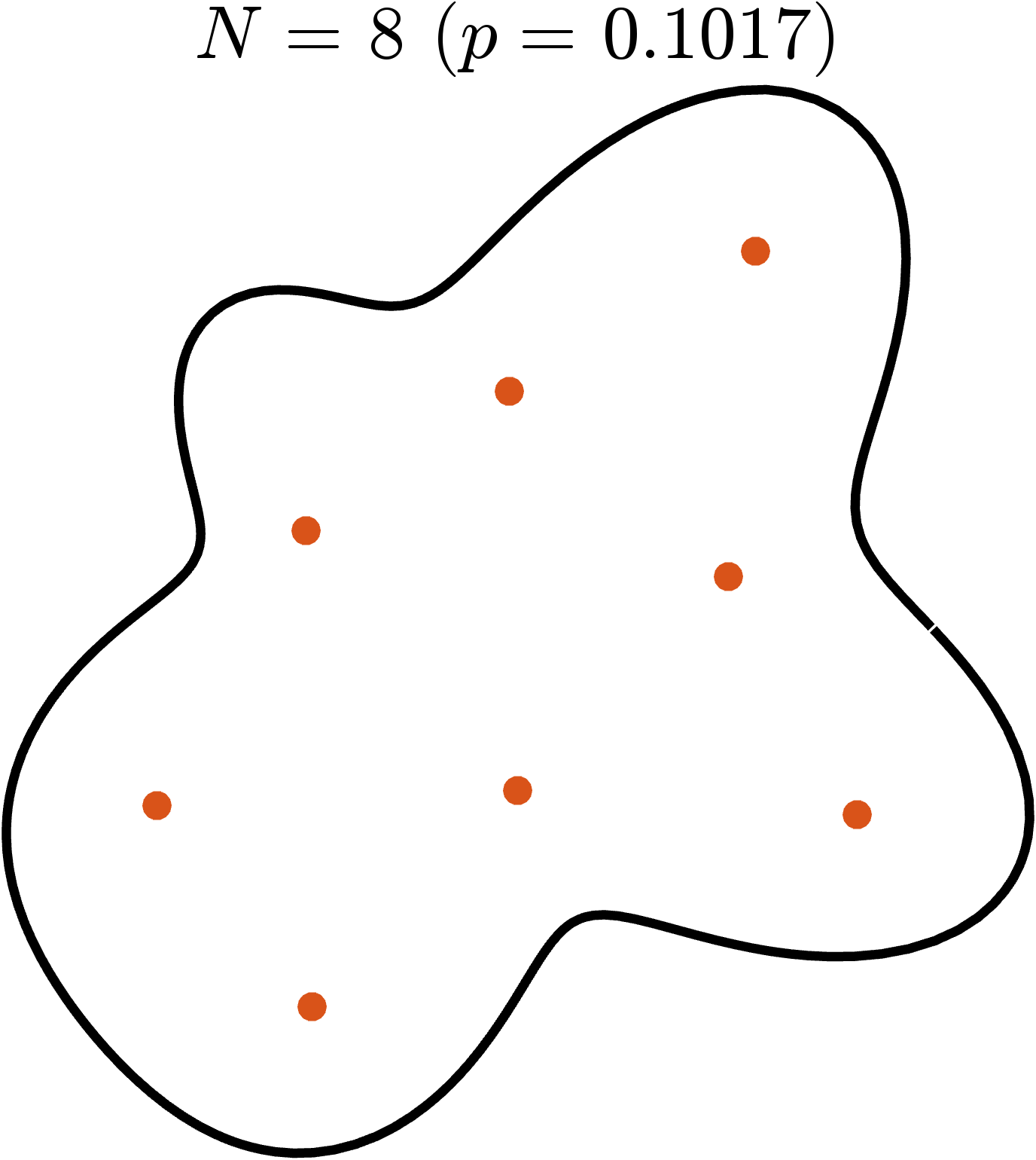}
\includegraphics[width=0.19\textwidth]{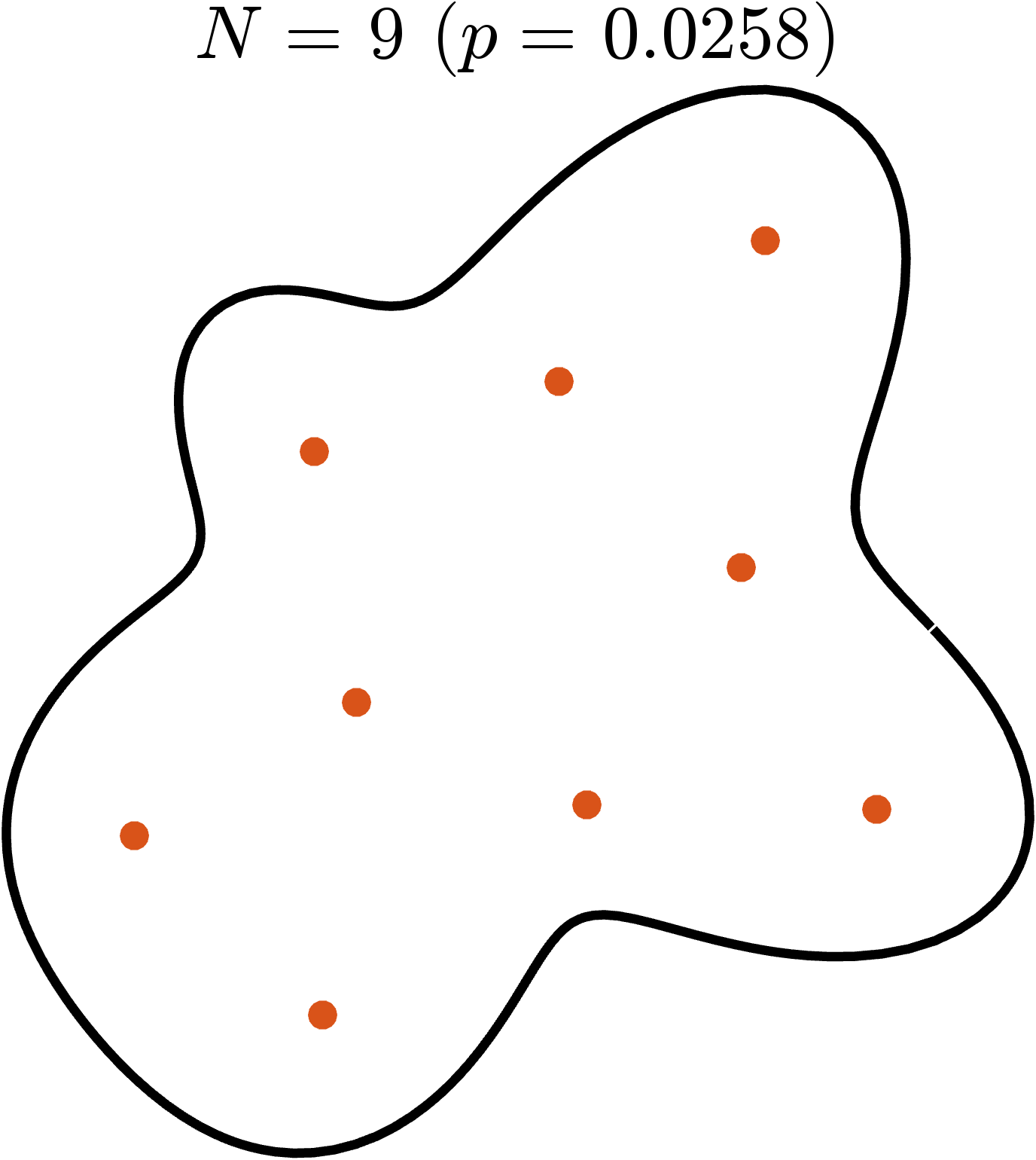}
\includegraphics[width=0.19\textwidth]{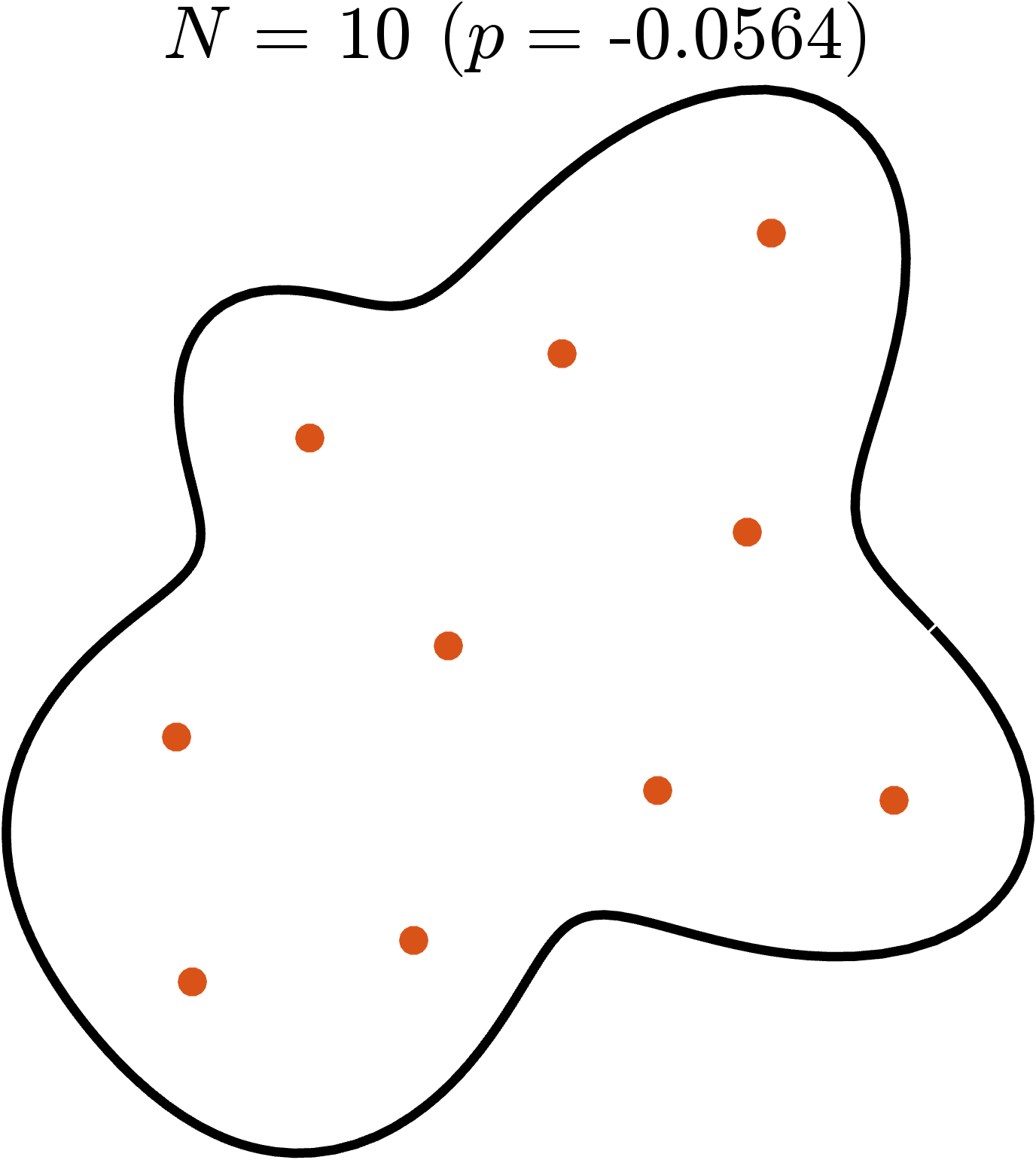}\\[5pt]
\includegraphics[width=0.19\textwidth]{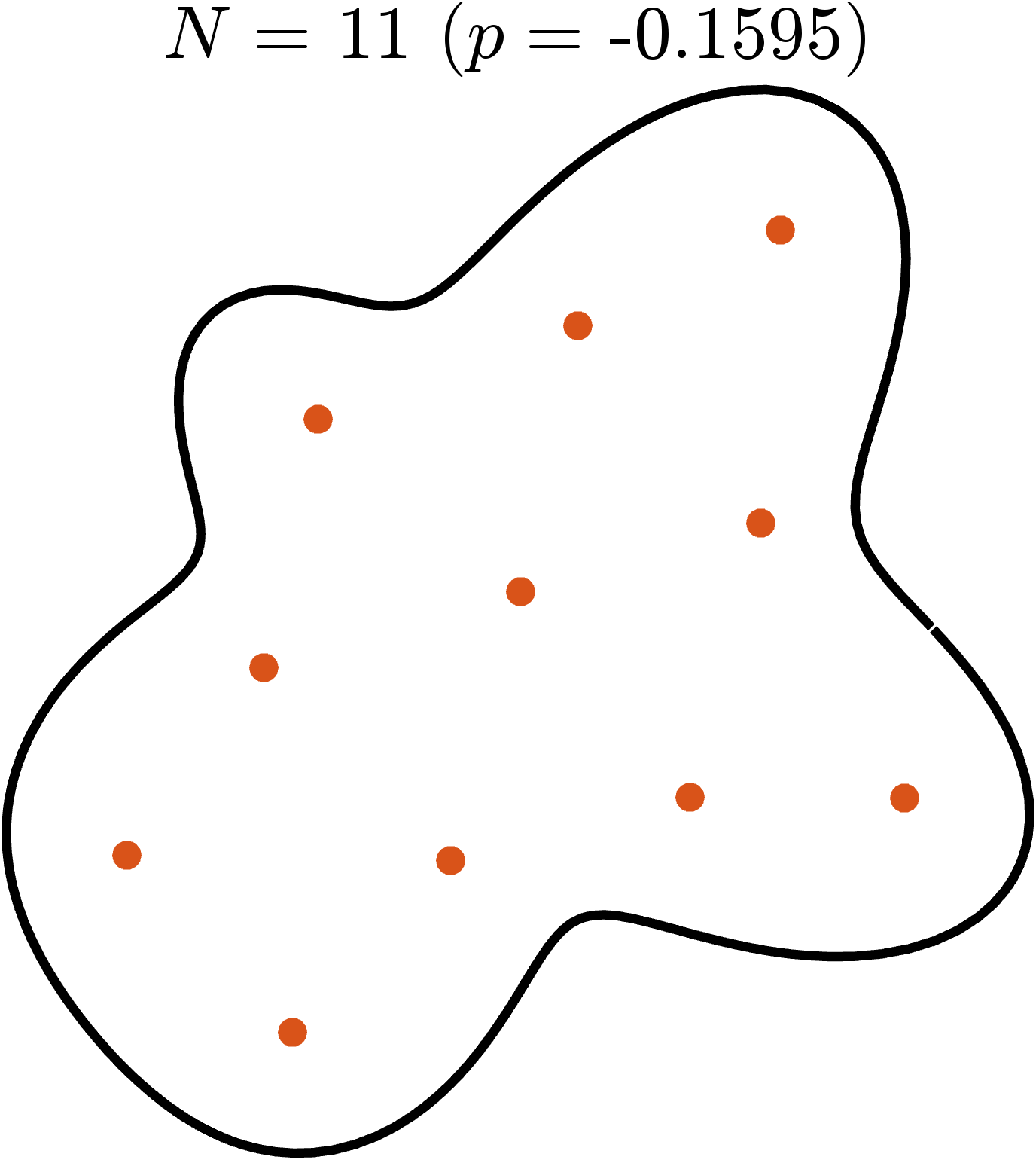}
 \includegraphics[width=0.19\textwidth]{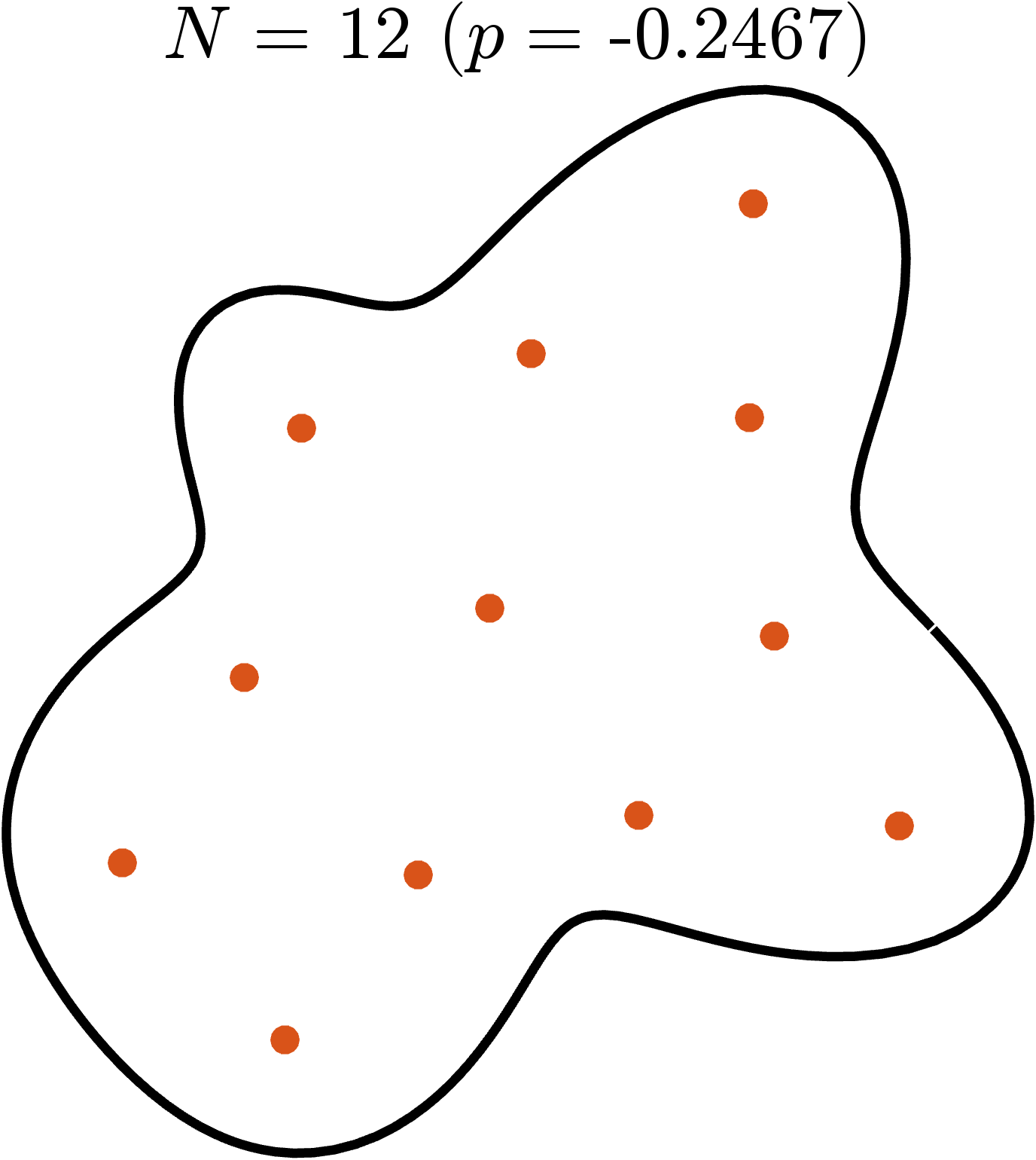}
\includegraphics[width=0.19\textwidth]{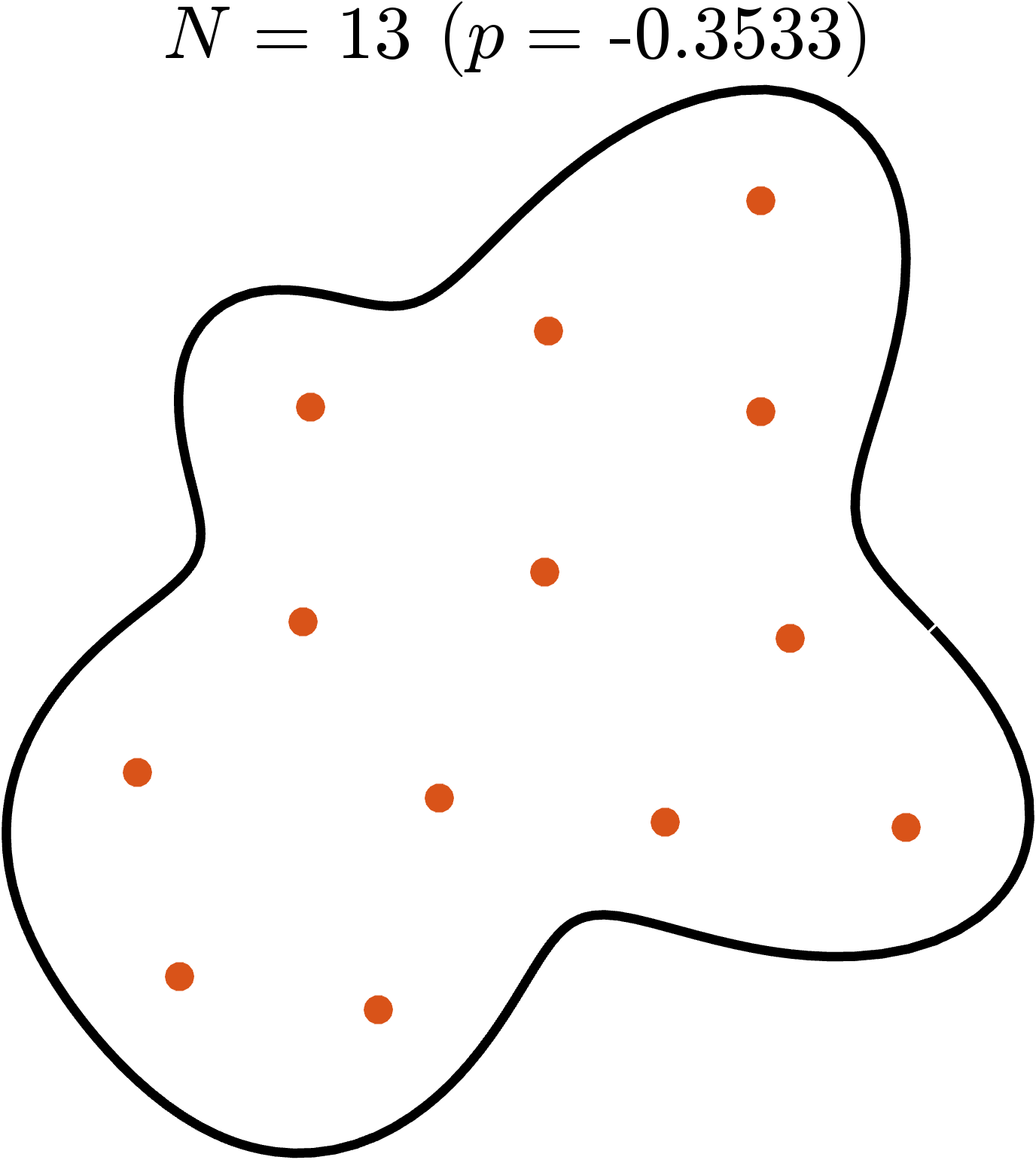}
\includegraphics[width=0.19\textwidth]{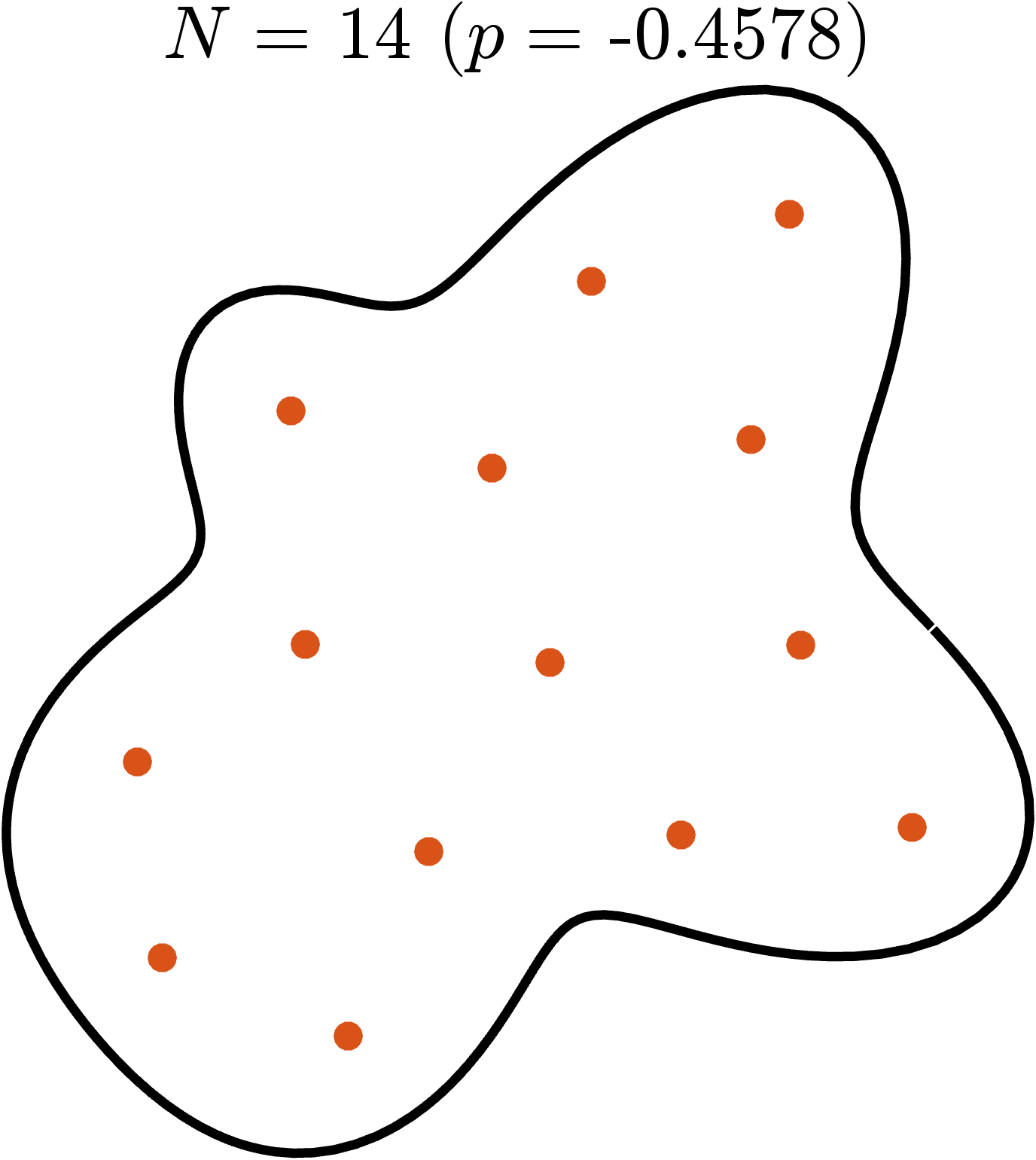}
\includegraphics[width=0.19\textwidth]{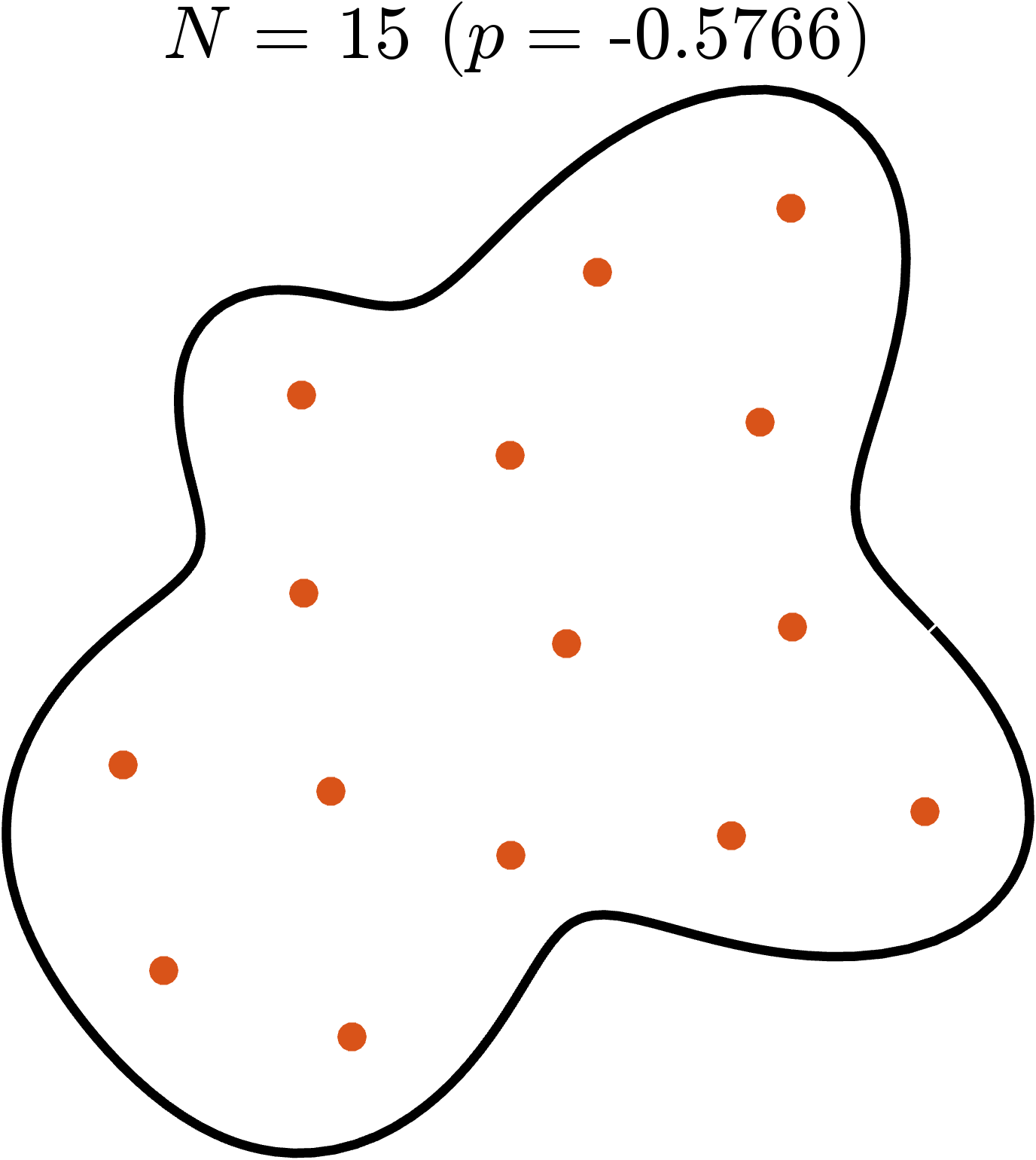}\\[5pt]
\includegraphics[width=0.19\textwidth]{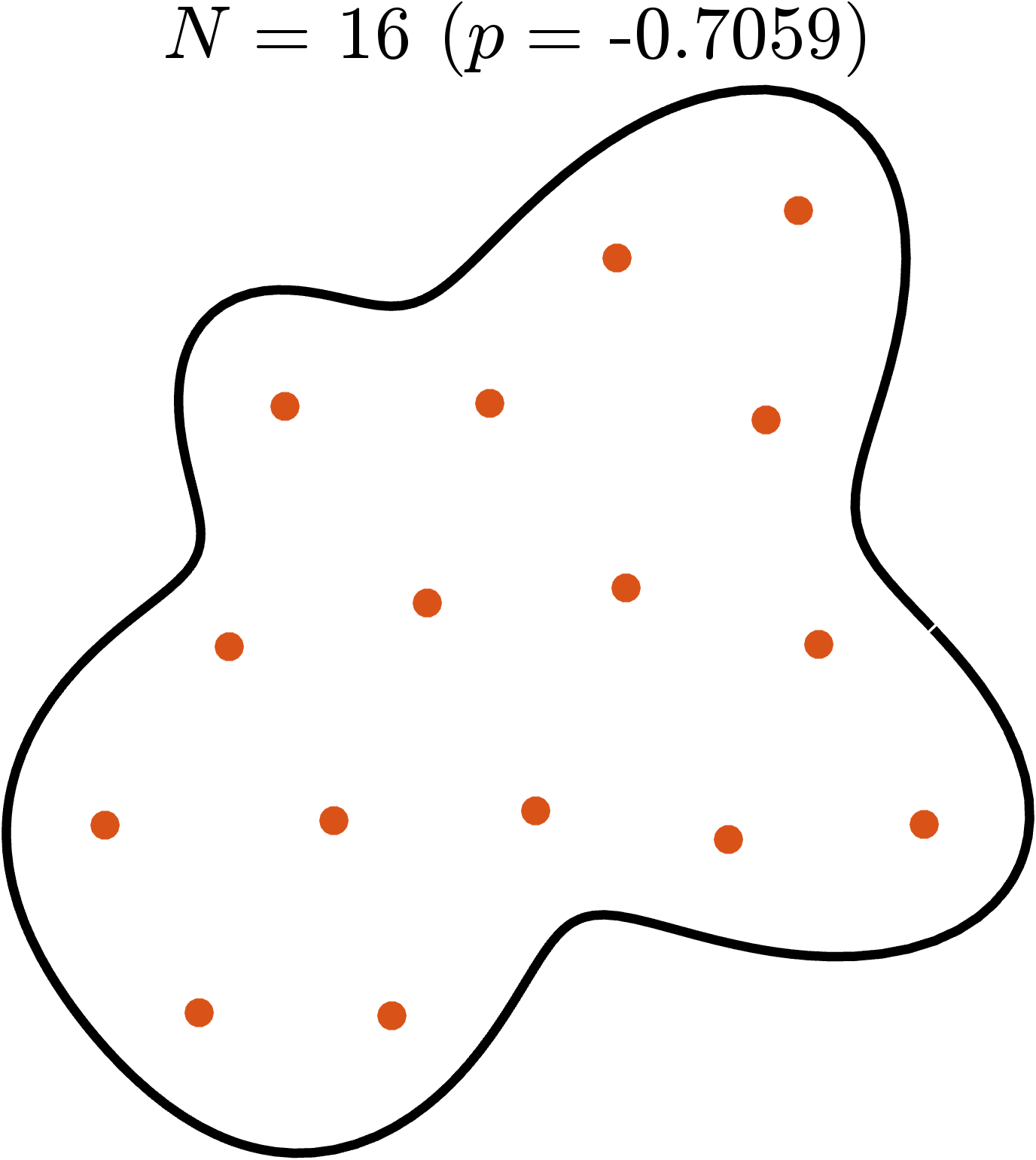}
\includegraphics[width=0.19\textwidth]{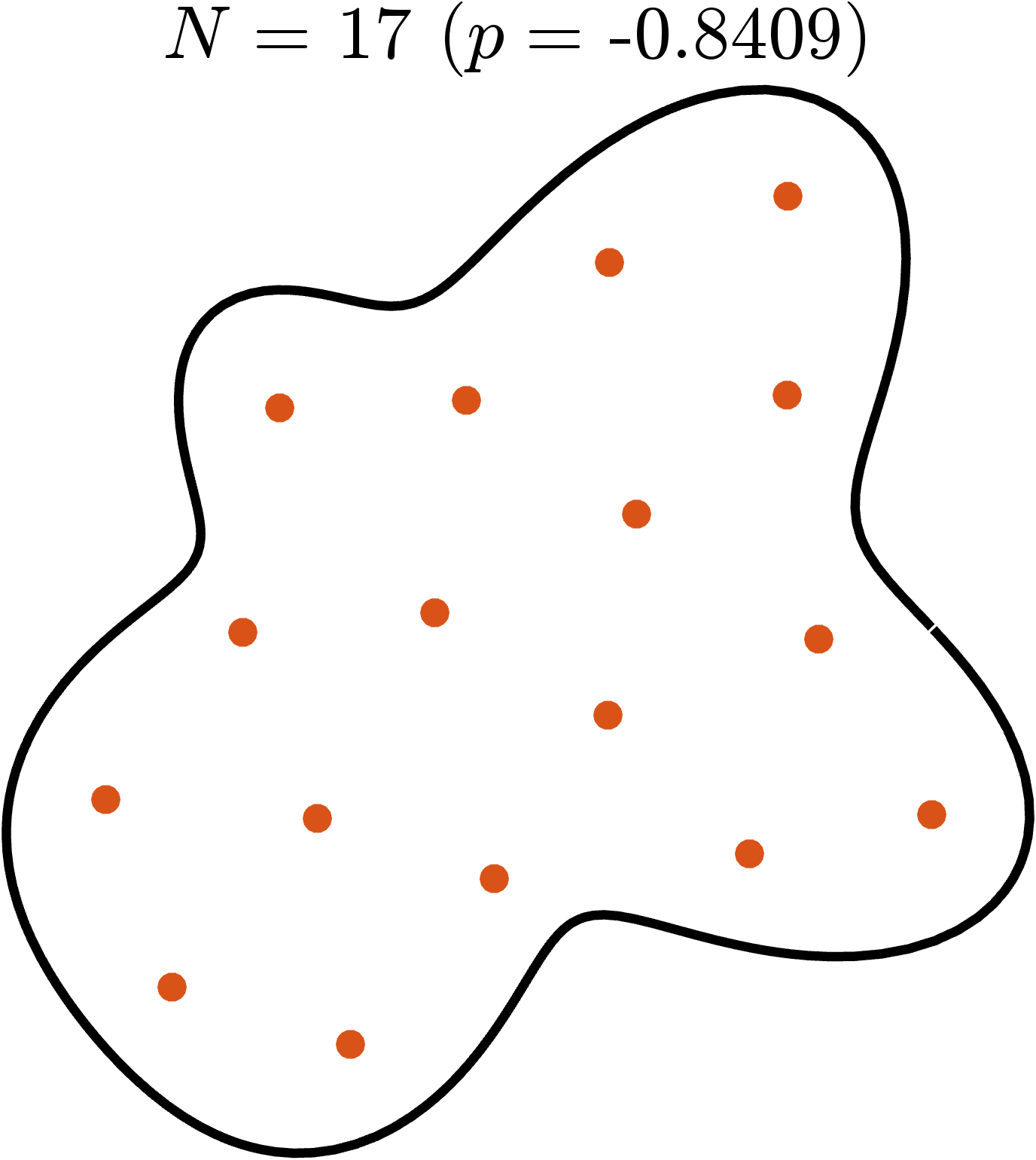}
\includegraphics[width=0.19\textwidth]{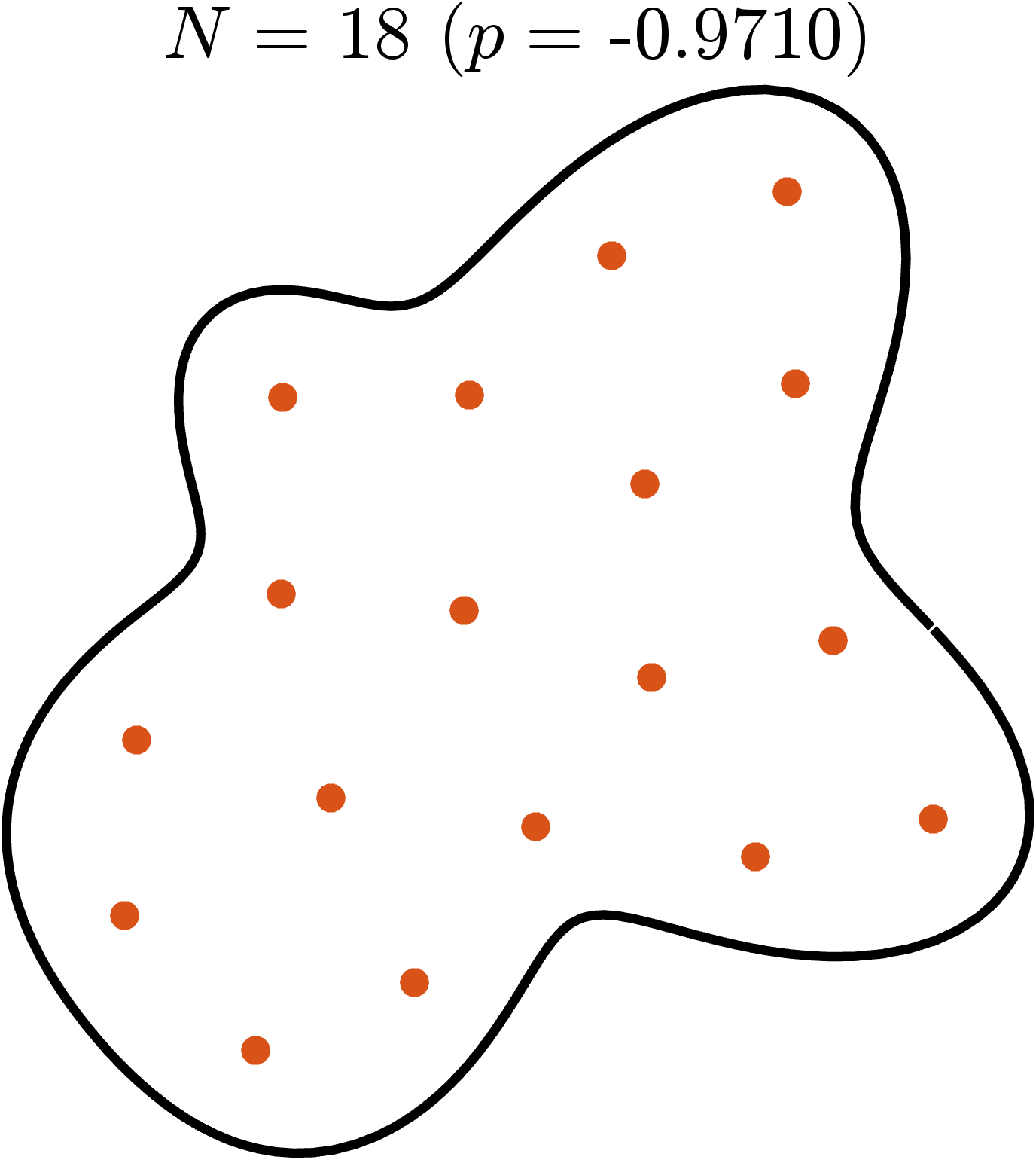}
\includegraphics[width=0.19\textwidth]{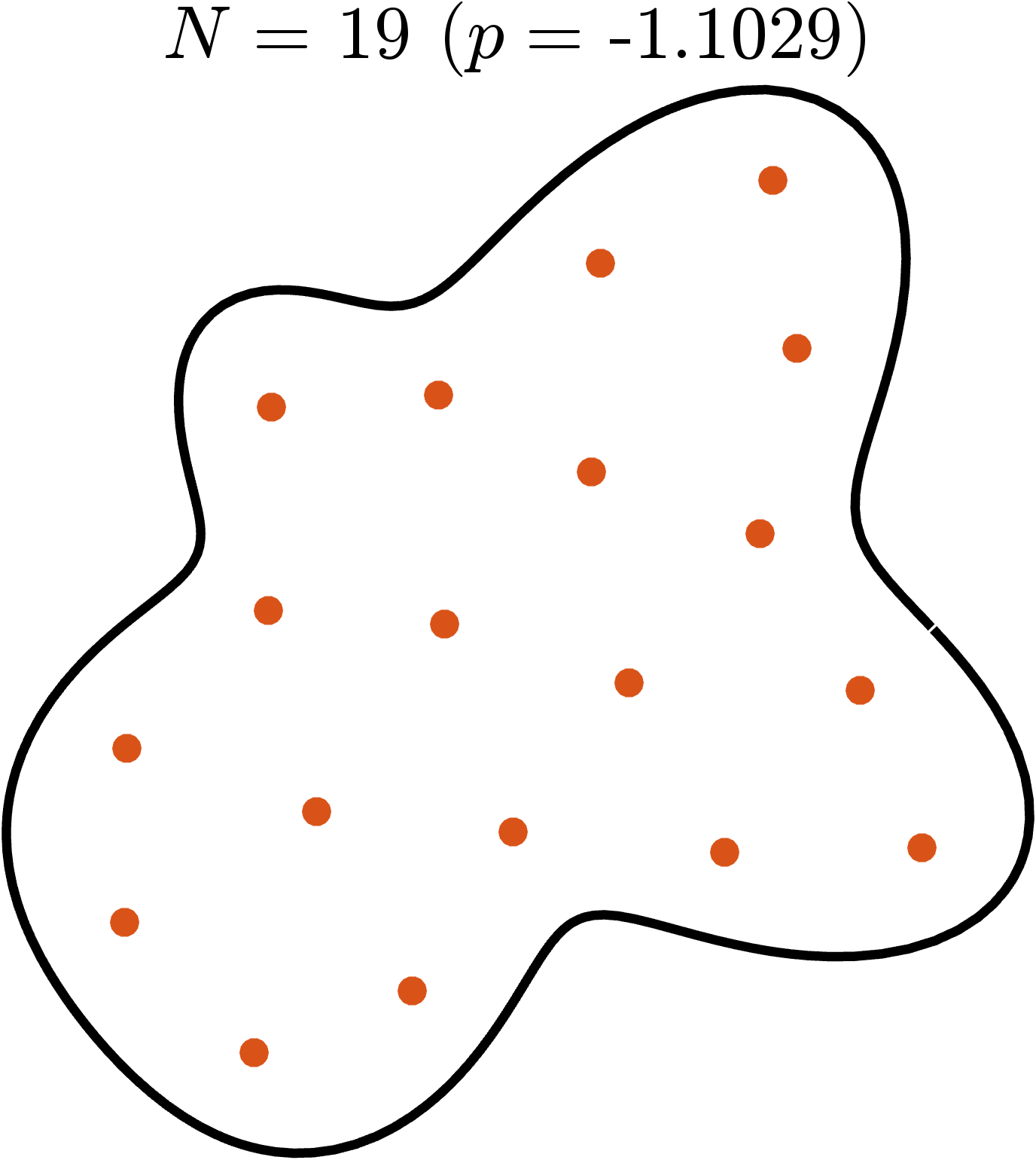}
\includegraphics[width=0.19\textwidth]{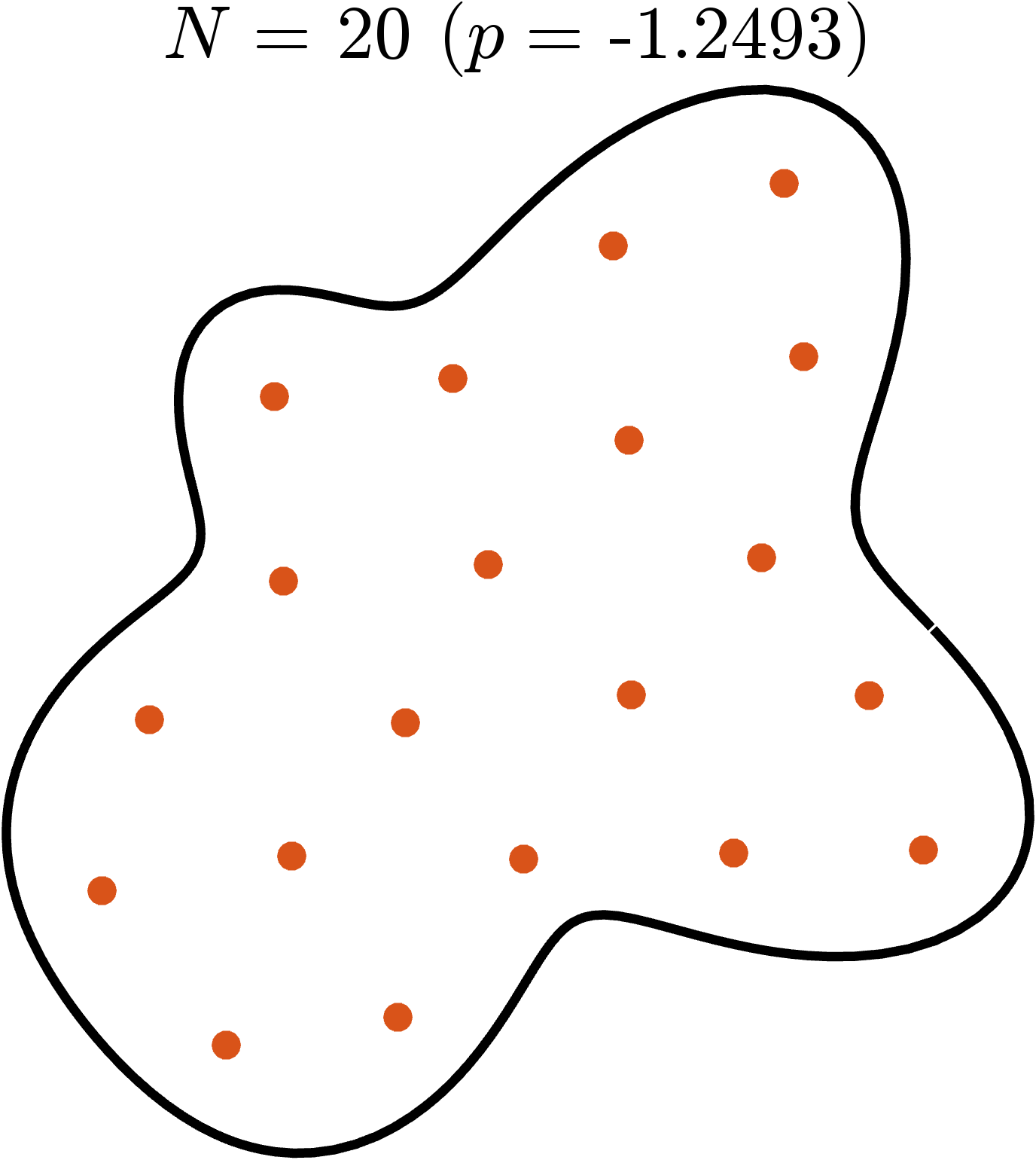}
\caption{Example of local minima of $p = p(\emph{\bx}_1,\ldots,\emph{\bx}_N)$ defined in \eqref{eq_Discrete} for a \emph{random} domain \eqref{eq:RandomDomain} with $N=1,\ldots,20$ interior traps.\label{fig:Random_Domain}}
\end{figure}

To demonstrate the method on more general geometries, we consider a family of random geometries generated by the parametric representation
\begin{equation}\label{eq:RandomDomain}
    \bx(\theta) = r(\theta)(\cos\theta,\sin\theta), \qquad r(\theta) = a_0 + \sum_{k=1}^M a_k \cos k\theta + b_k \sin k \theta.
\end{equation}
Our process is to sample standard normal random variables for $\{a_k\}_{k=1}^M$ and $\{b_k\}_{k=1}^M$ then set $a_0 = 1.1\times\sum_{k=1}^N \big(|a_k| + |b_k|\big)$. In Fig.~\ref{fig:Random_Domain}, we show the result of numerical optimization of \eqref{eq_Discrete} for trap numbers $N=1,\ldots,20$. At smaller $N$ we observe a tendency of traps to be located near the center and the geometric lobes of the region. For large numbers, we see that the traps become almost equally spaced throughout the domain.

\begin{figure}
    \centering
    \includegraphics[width = 0.24\textwidth]{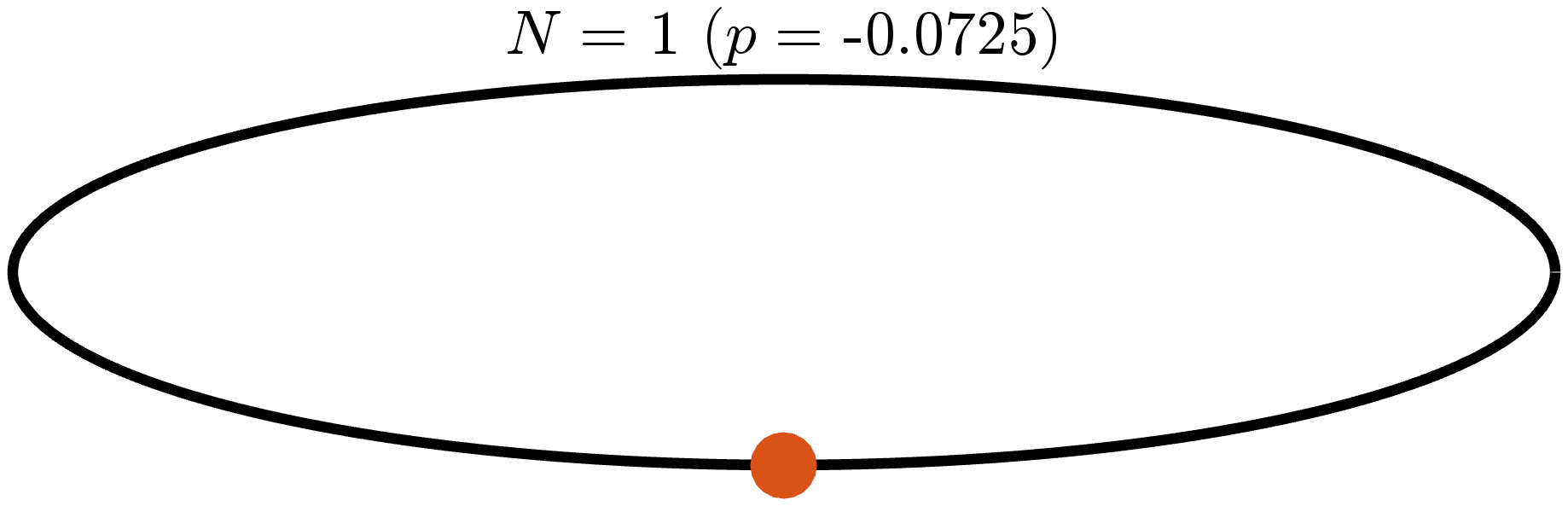}
    \includegraphics[width = 0.24\textwidth]{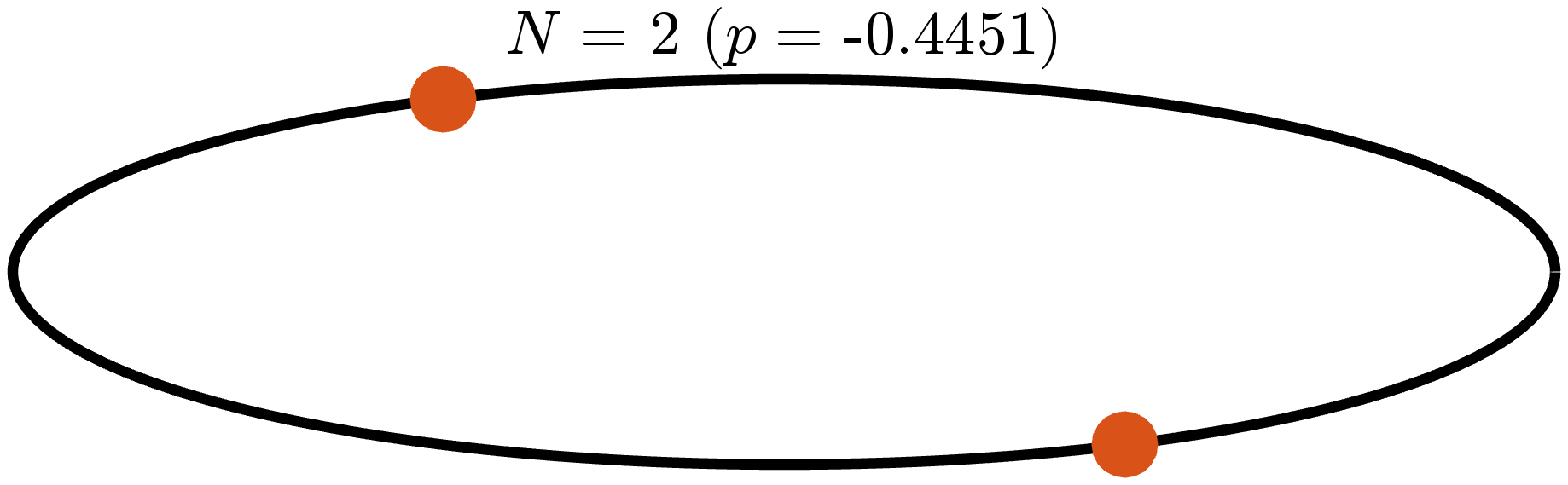}
    \includegraphics[width = 0.24\textwidth]{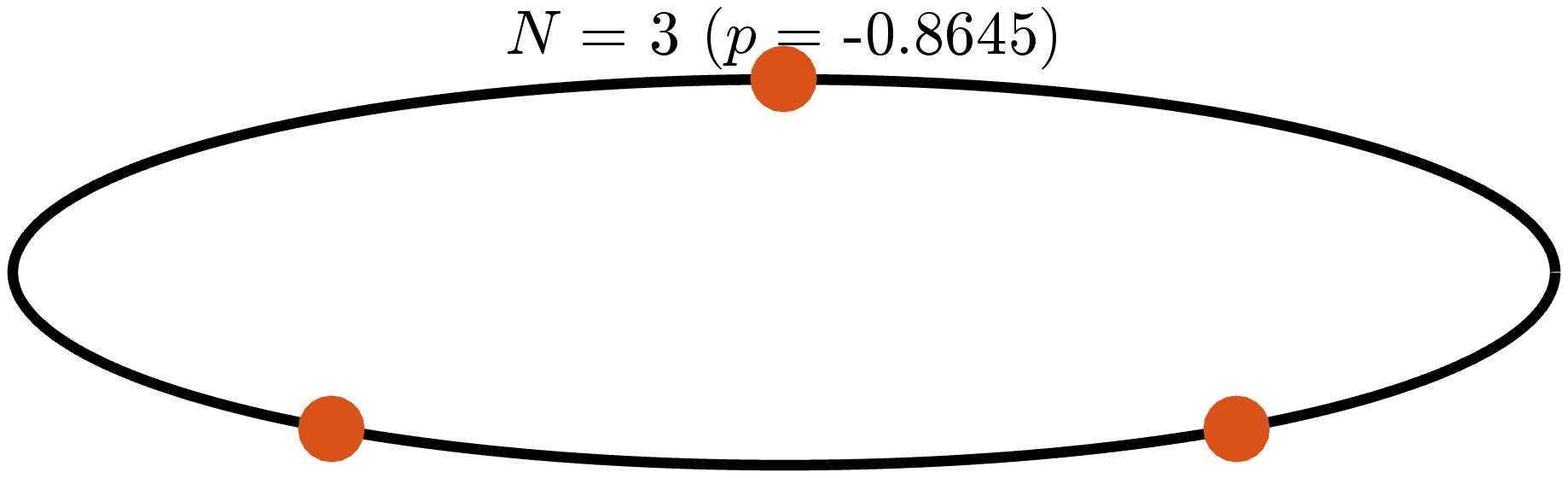}
    \includegraphics[width = 0.24\textwidth]{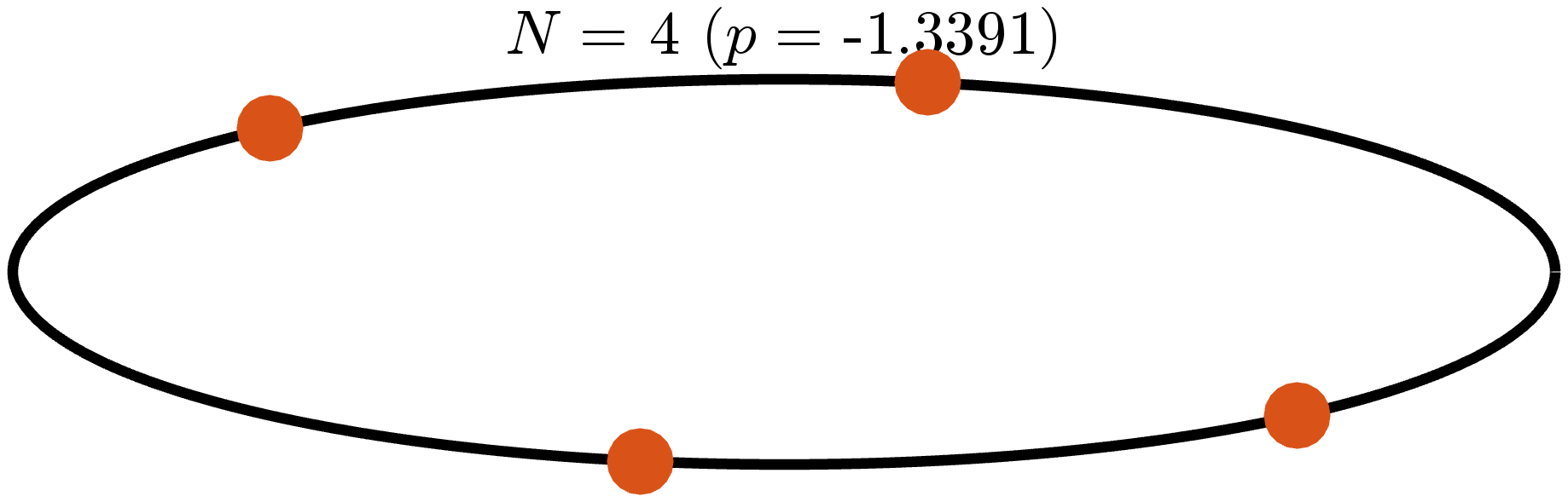}\\[5pt]
    \includegraphics[width = 0.24\textwidth]{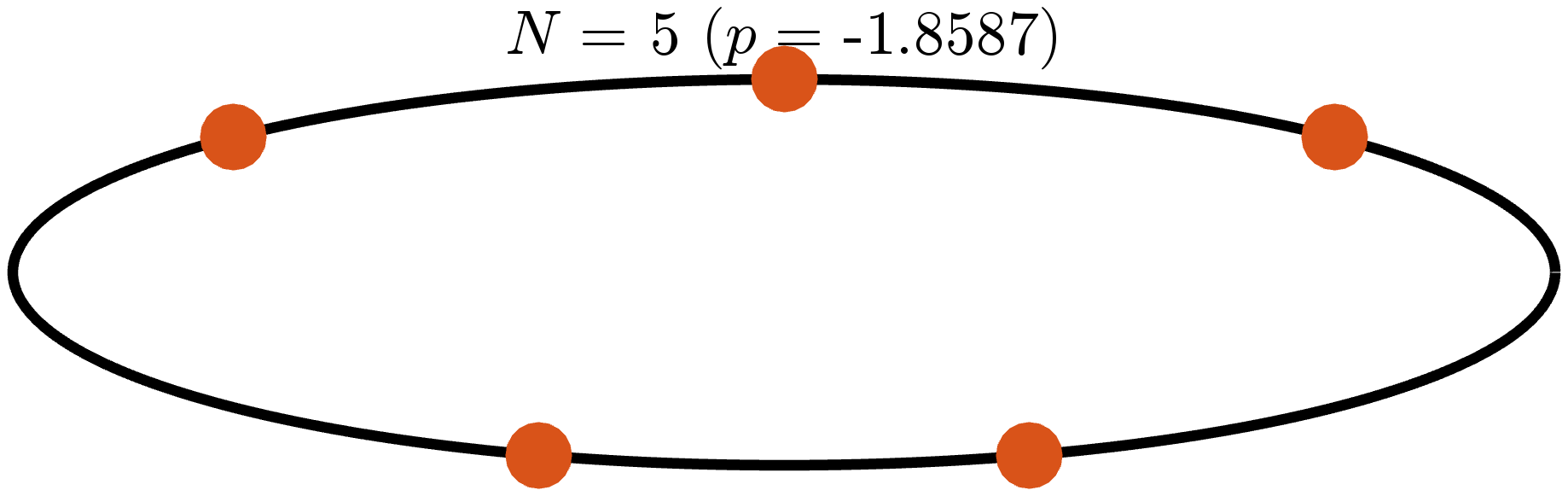}
    \includegraphics[width = 0.24\textwidth]{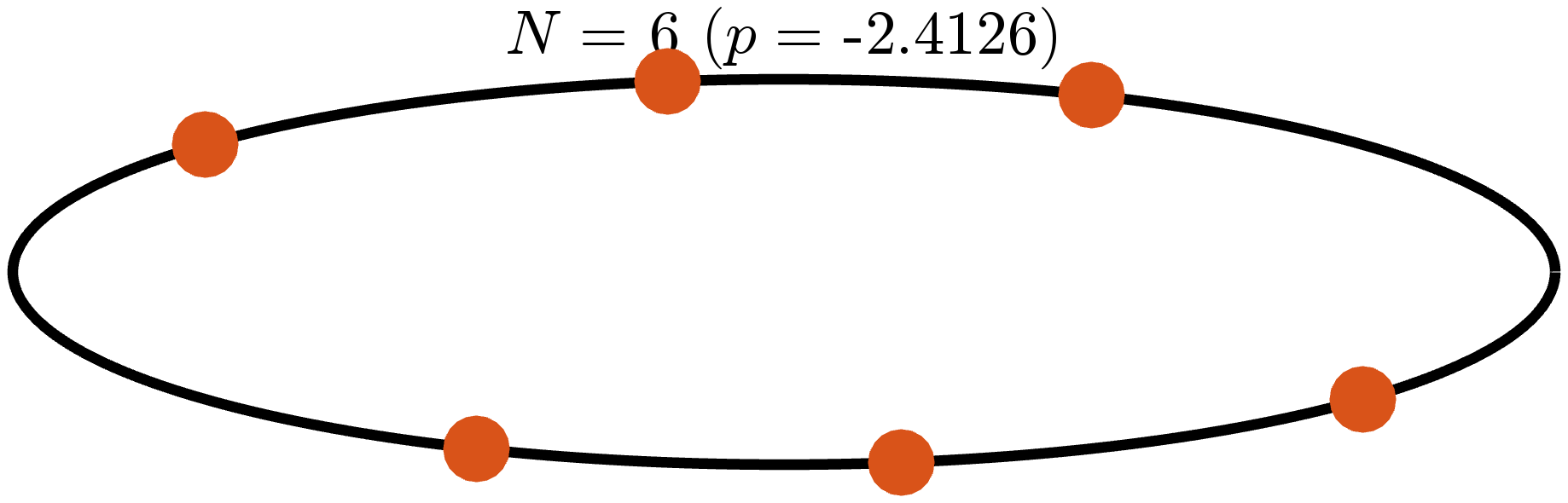}
    \includegraphics[width = 0.24\textwidth]{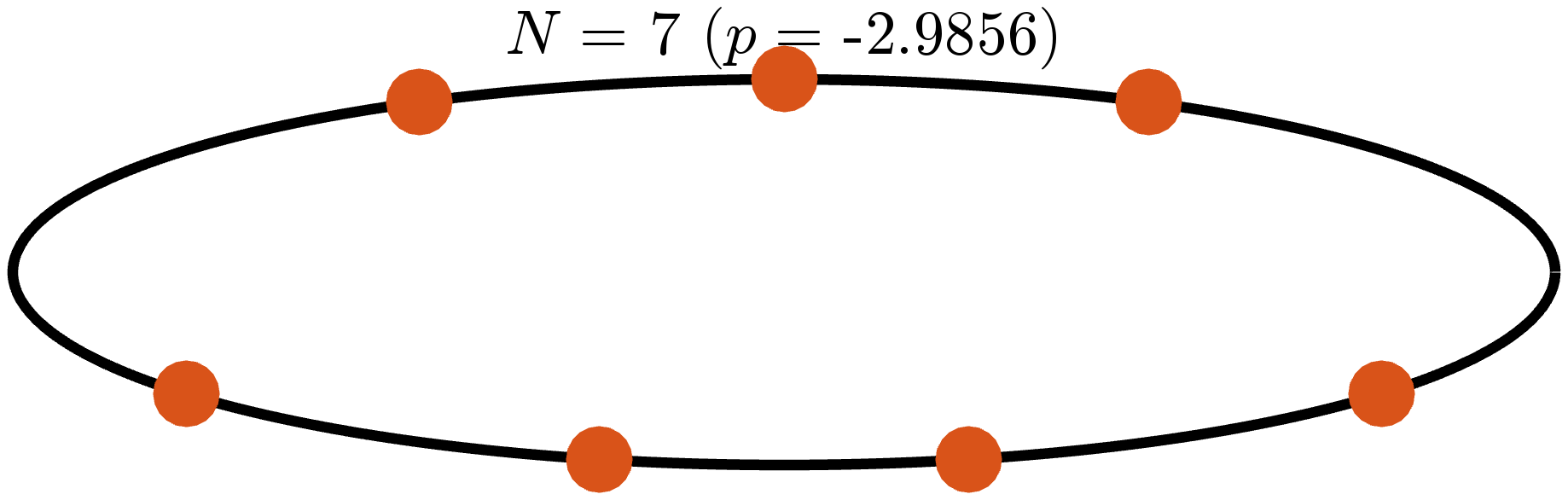}
    \includegraphics[width = 0.24\textwidth]{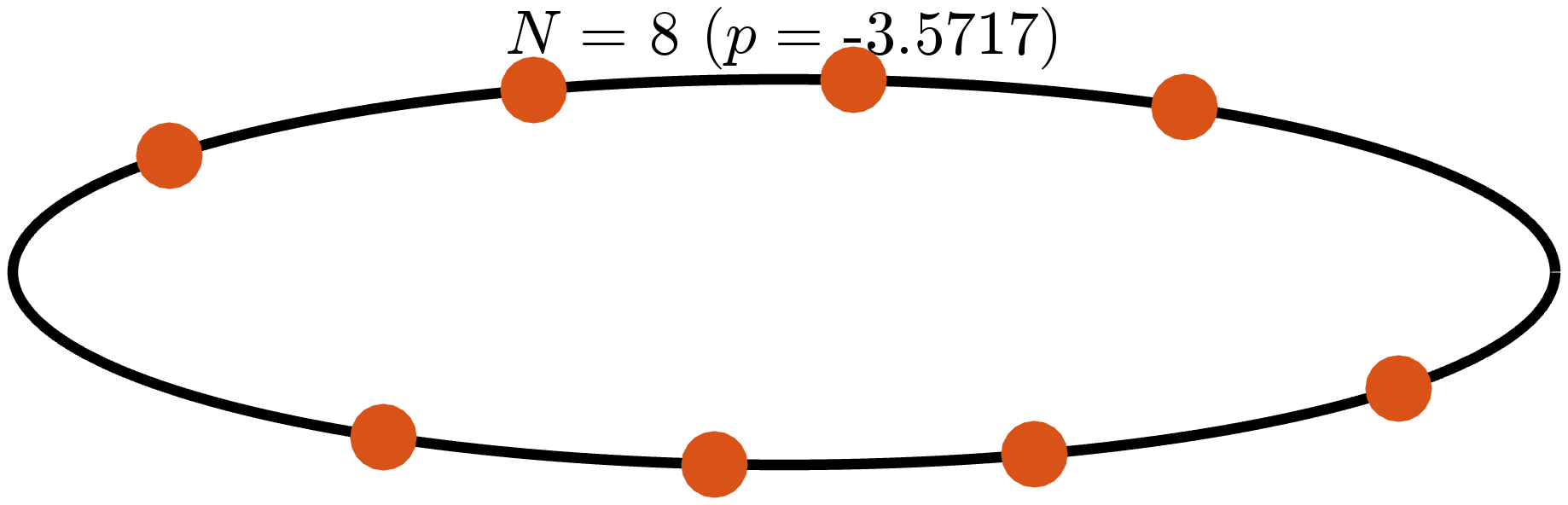}\\[5pt]
    \includegraphics[width = 0.24\textwidth]{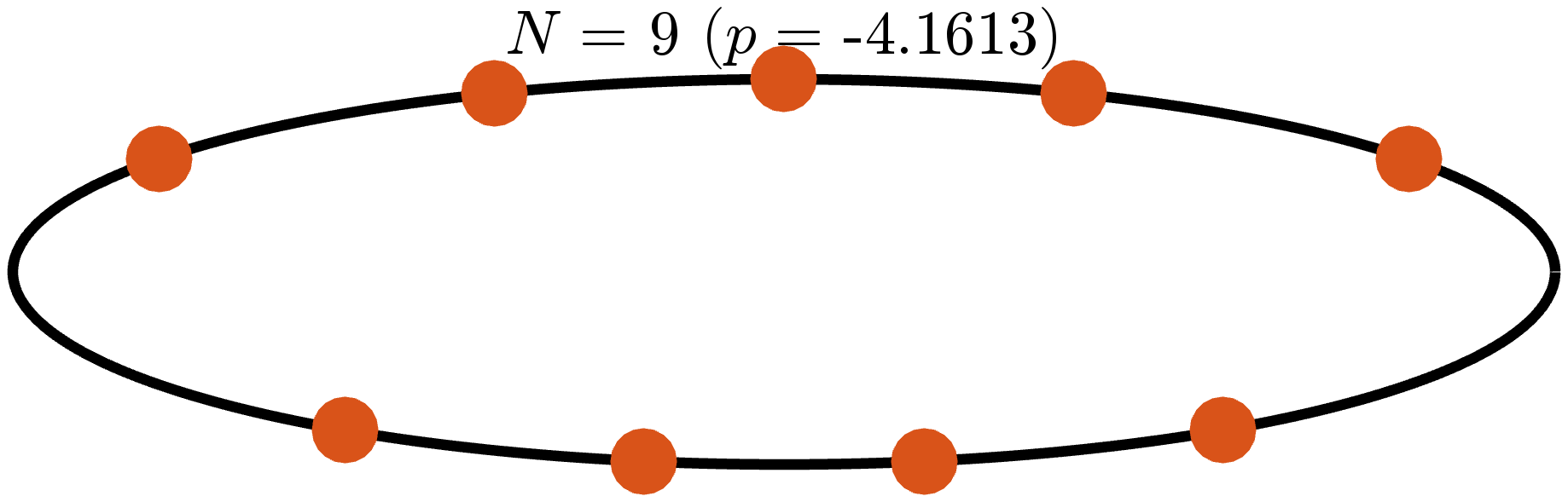}
    \includegraphics[width = 0.24\textwidth]{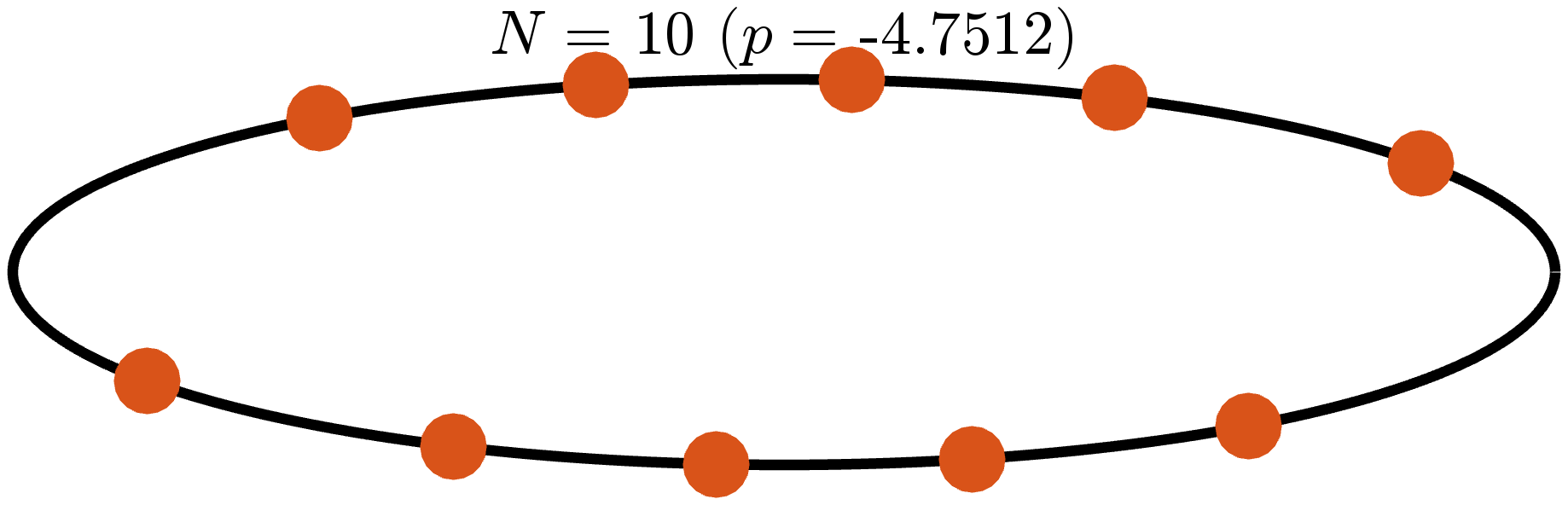}
    \includegraphics[width = 0.24\textwidth]{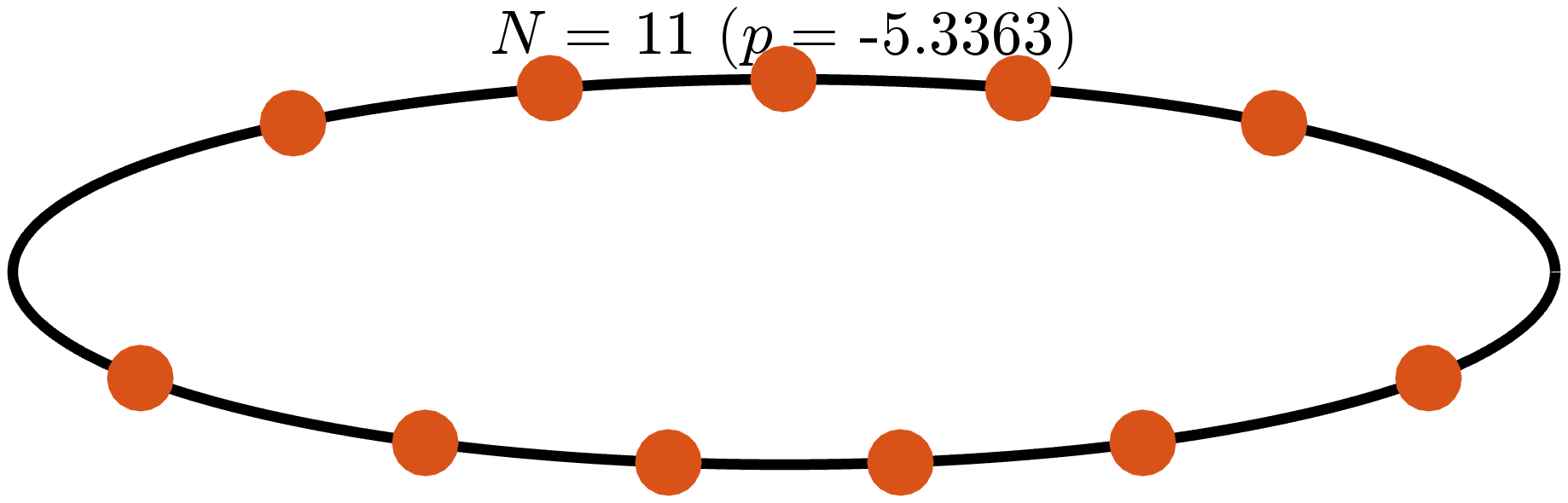}
    \includegraphics[width = 0.24\textwidth]{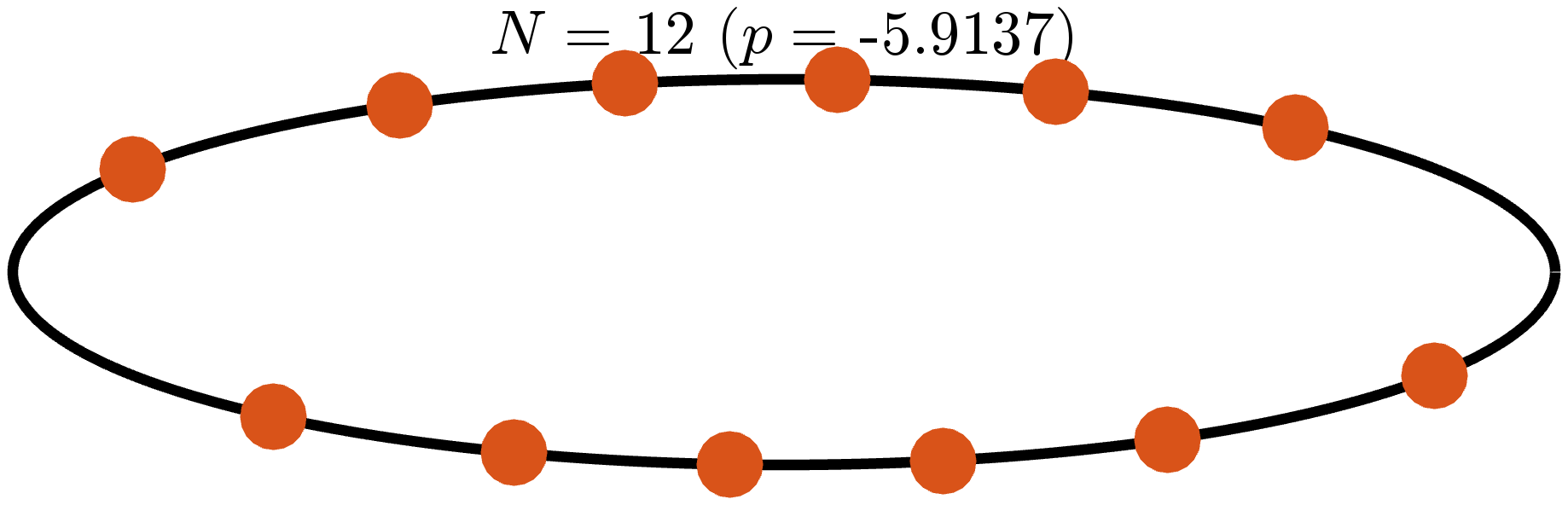}
    \caption{Minimizing locations of the discrete energy $p(\emph{\bx}_1,\ldots,\emph{\bx}_N)$ for boundary windows located on an ellipse of area $|\Omega| = \pi$ and dimensions $a= 2$ $b = 1/2$. From top left to bottom right, the patterns are shown for $N=1$ to $N=12$ windows respectively.}
\end{figure}

\begin{figure}
    \centering
    \includegraphics[width = 0.159\textwidth]{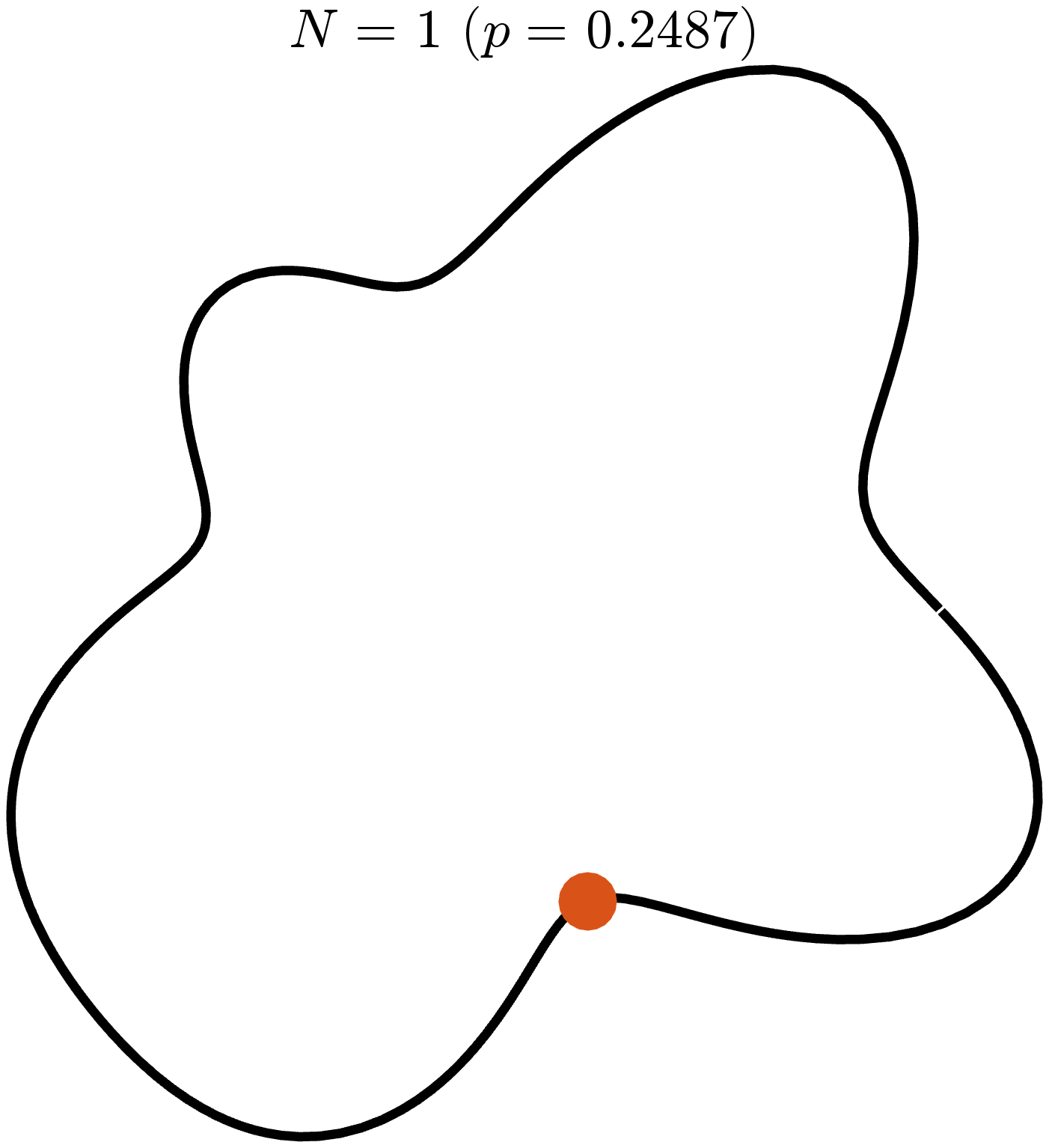}
    \includegraphics[width = 0.159\textwidth]{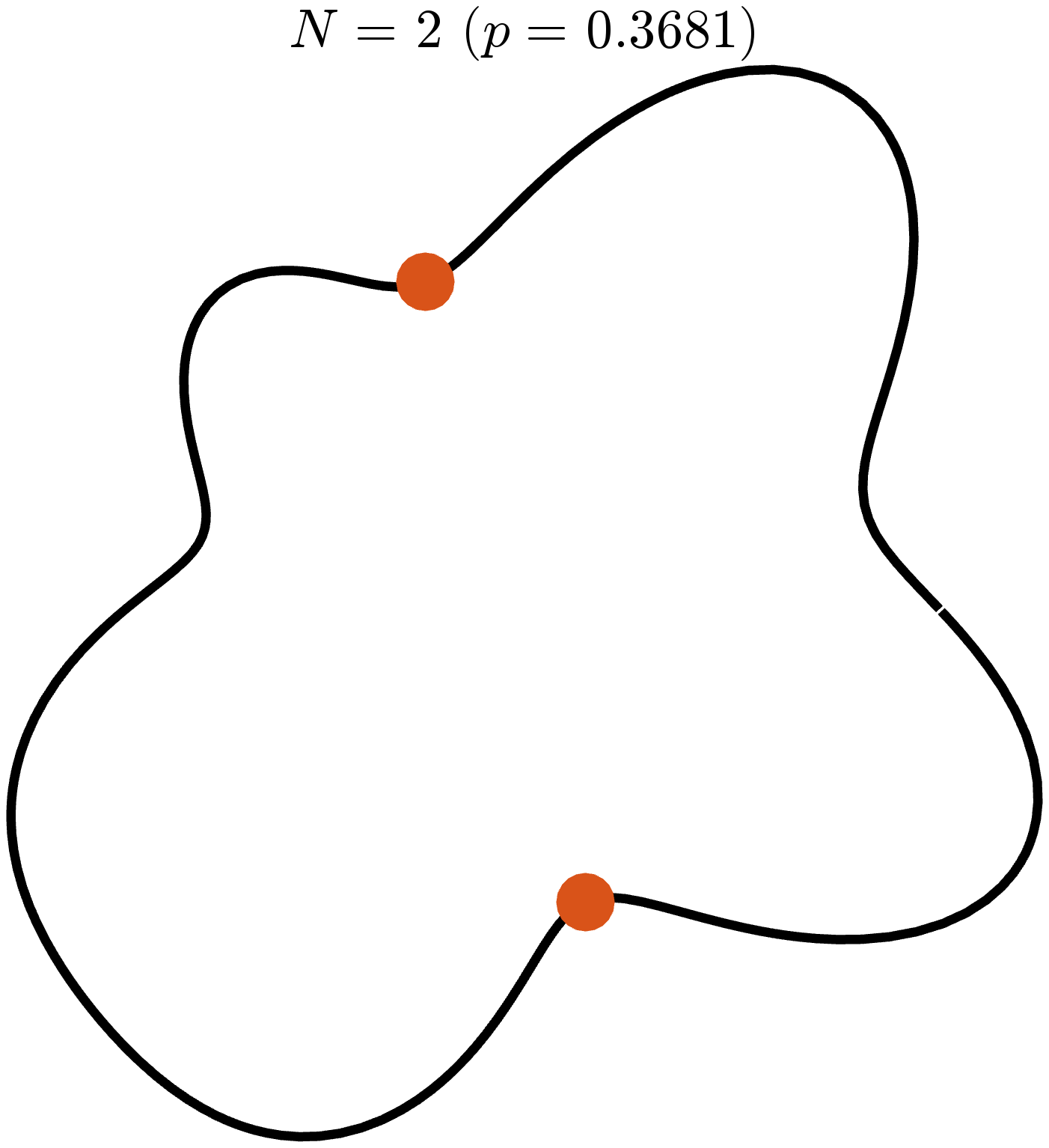}
    \includegraphics[width = 0.159\textwidth]{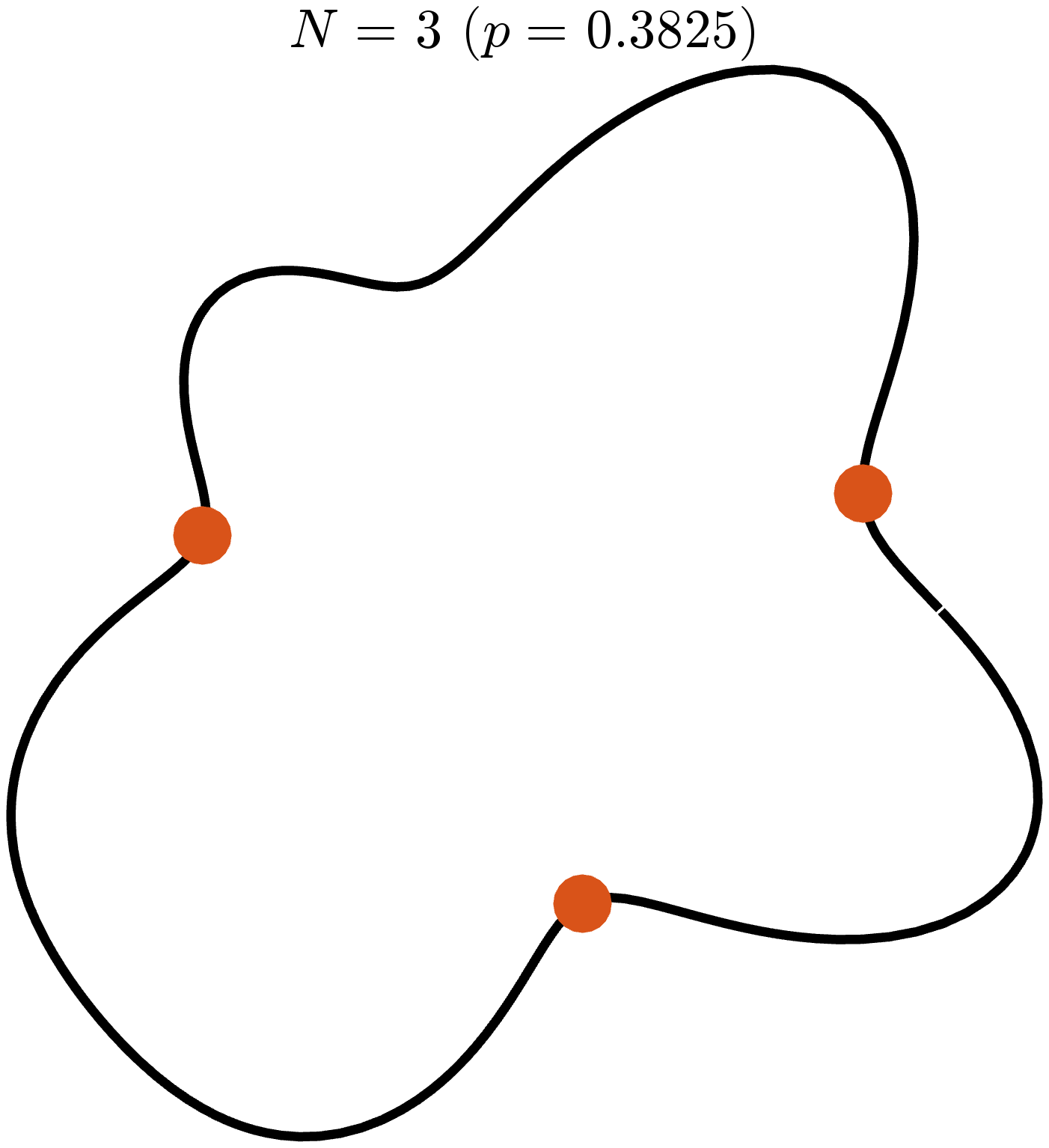}
    \includegraphics[width = 0.159\textwidth]{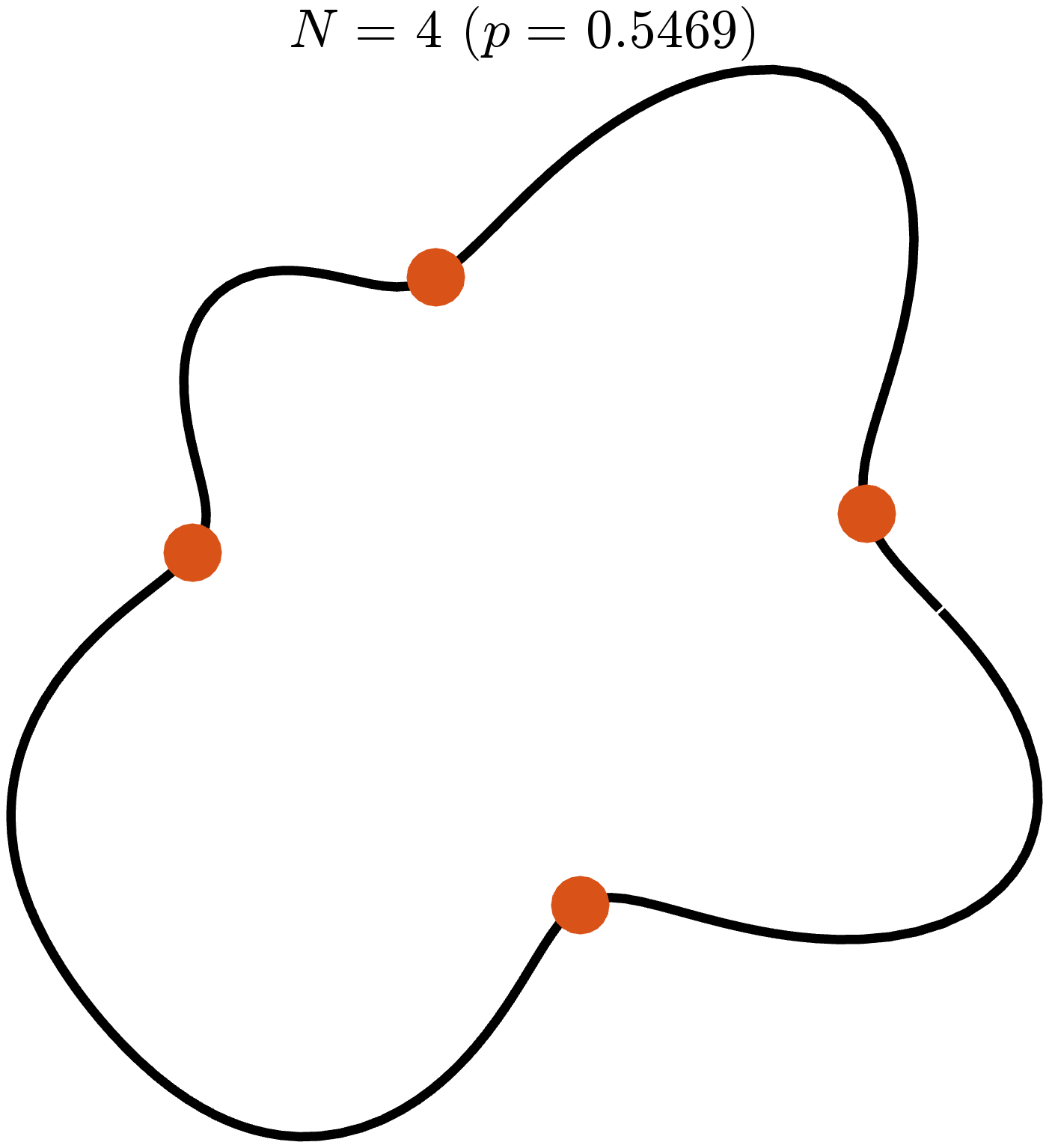}
    \includegraphics[width = 0.159\textwidth]{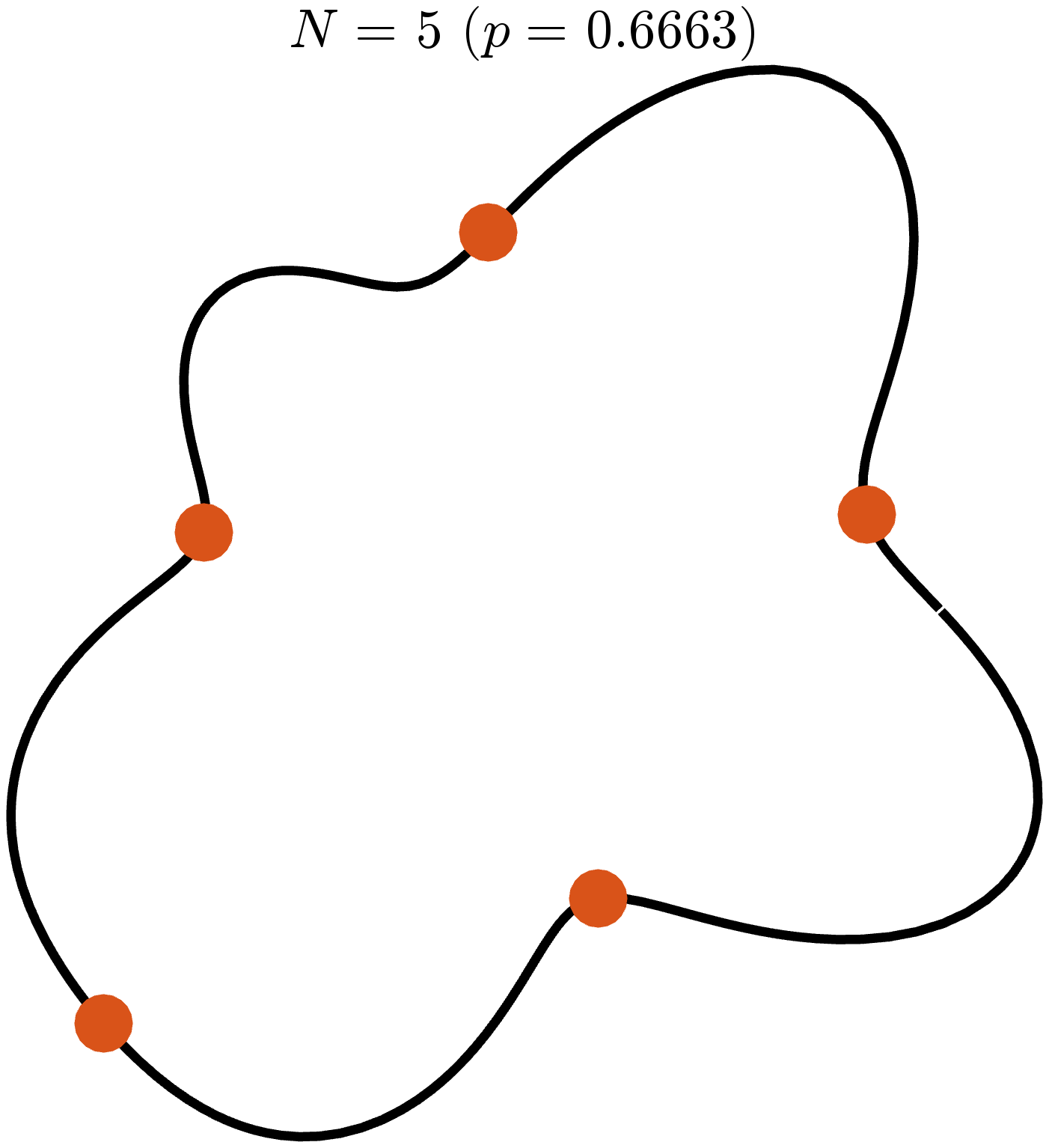}
    \includegraphics[width = 0.159\textwidth]{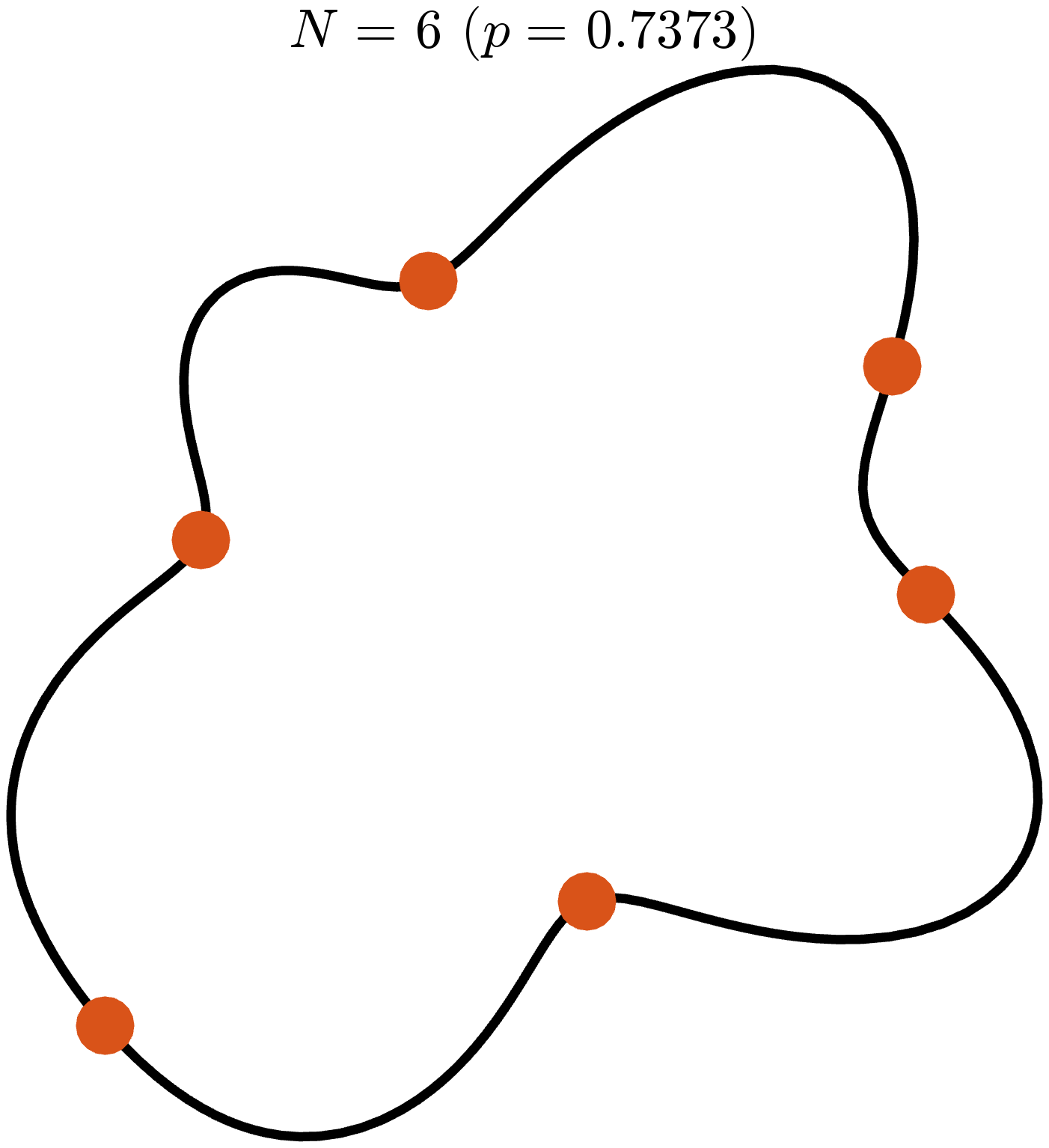}\\[5pt]
    \includegraphics[width = 0.159\textwidth]{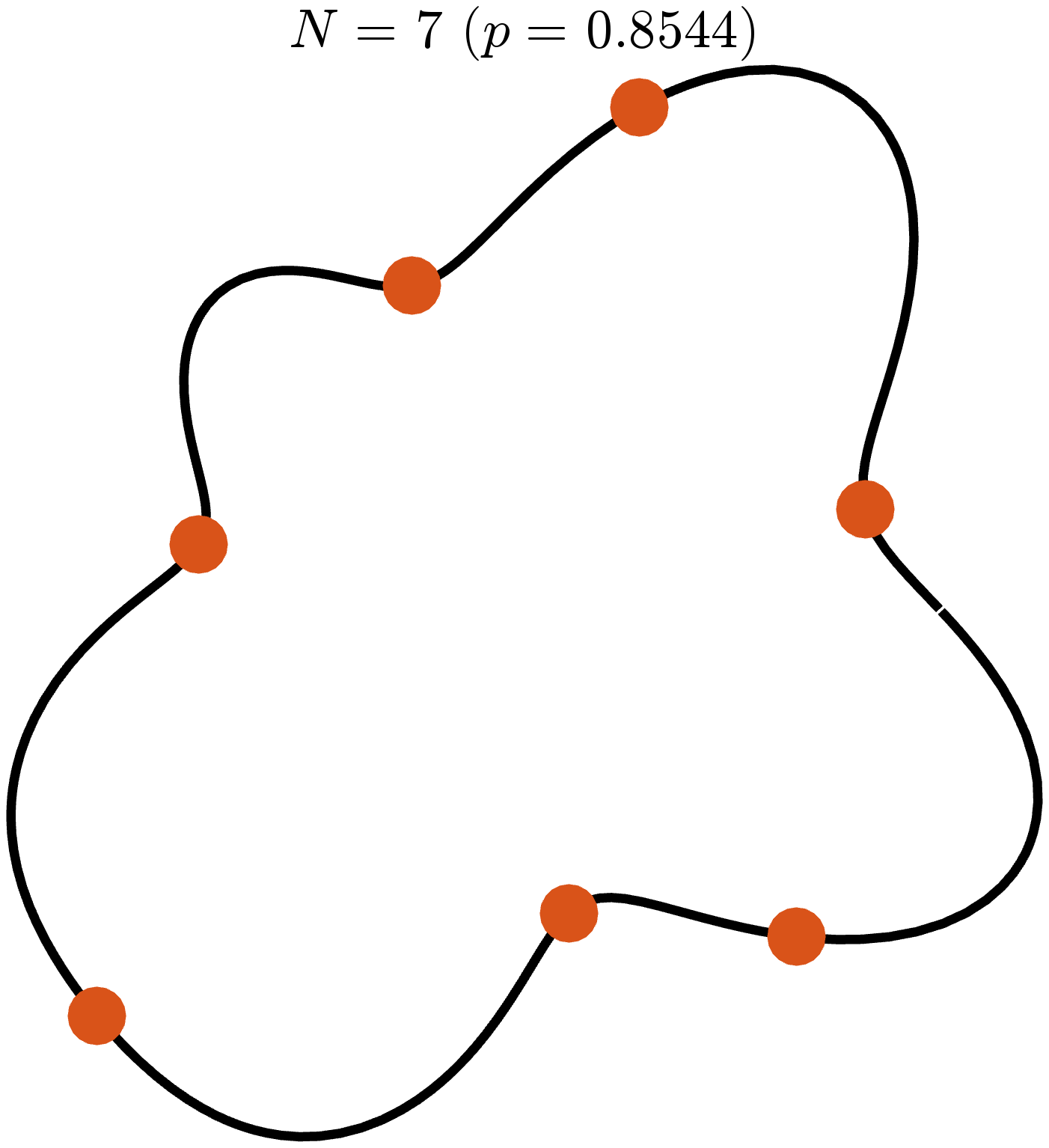}
    \includegraphics[width = 0.159\textwidth]{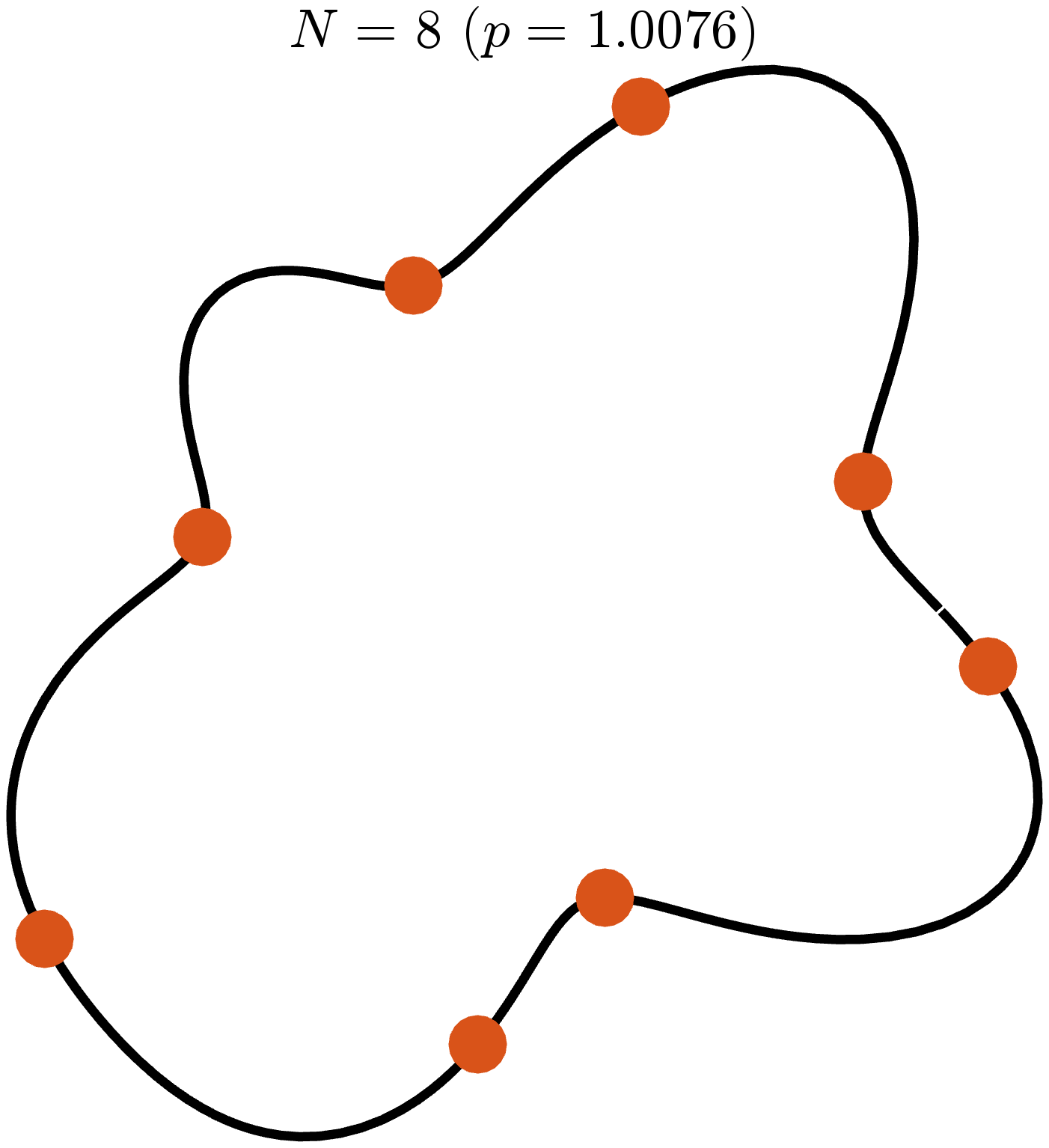}
    \includegraphics[width = 0.159\textwidth]{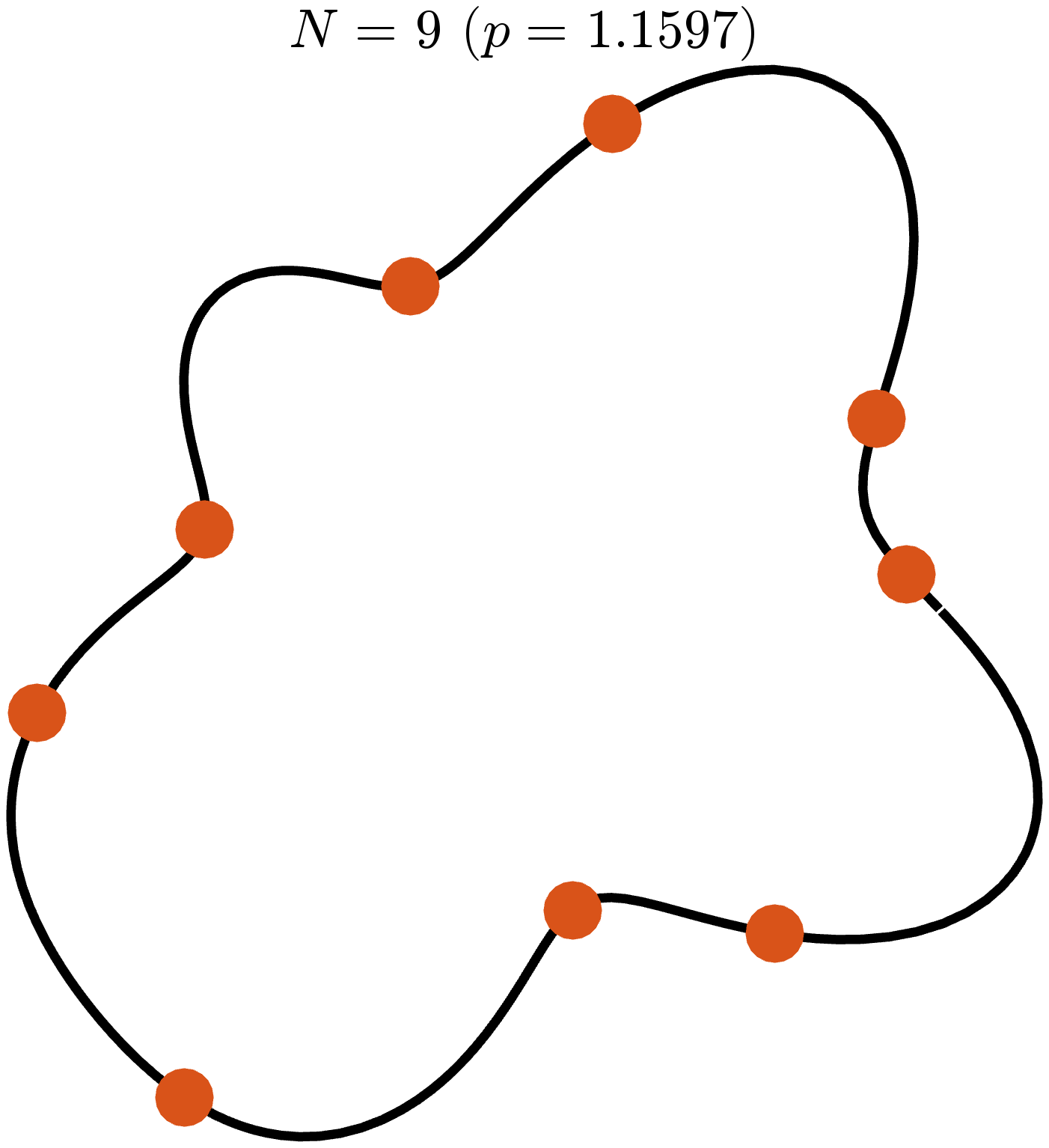}
    \includegraphics[width = 0.159\textwidth]{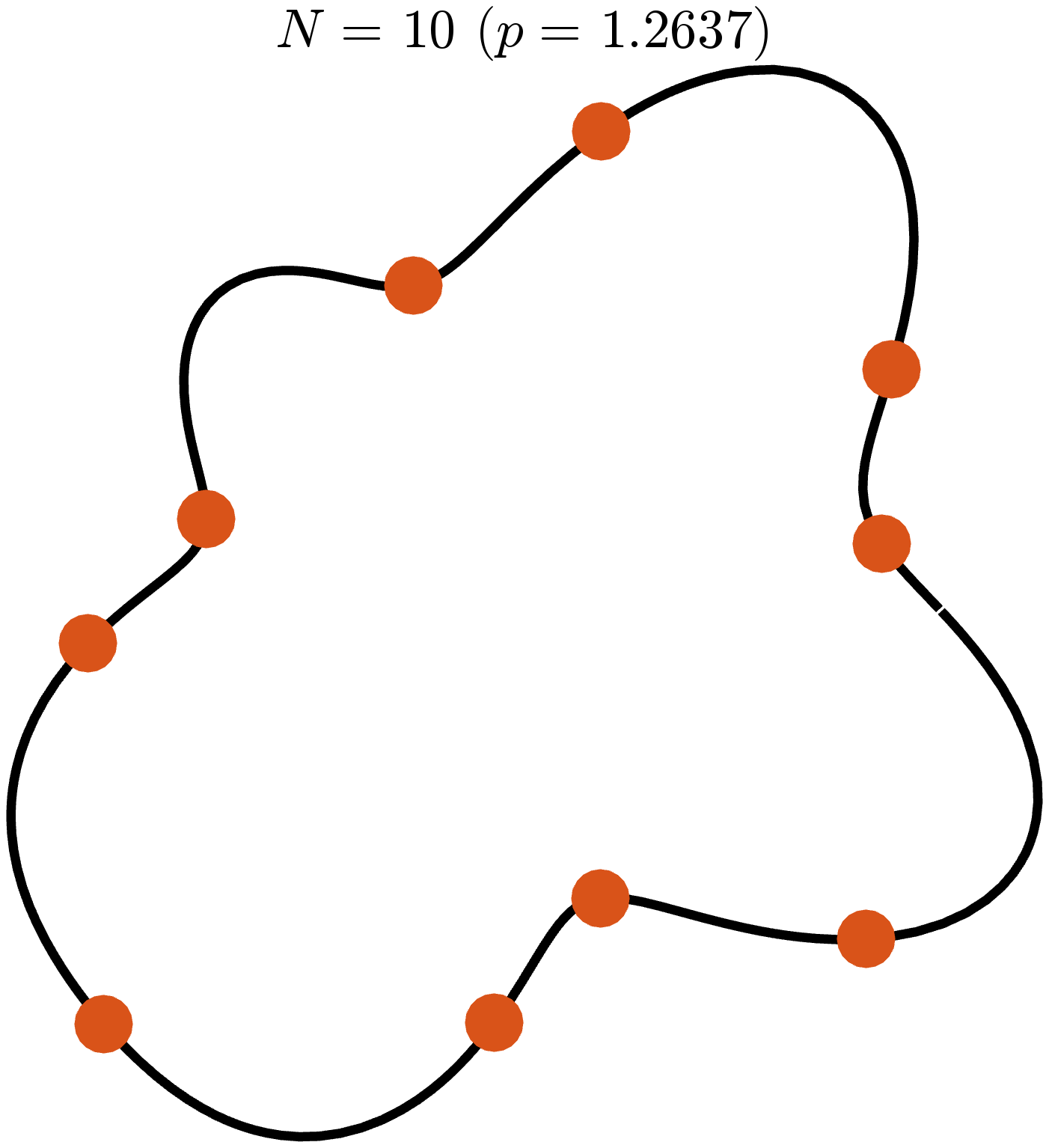}
    \includegraphics[width = 0.159\textwidth]{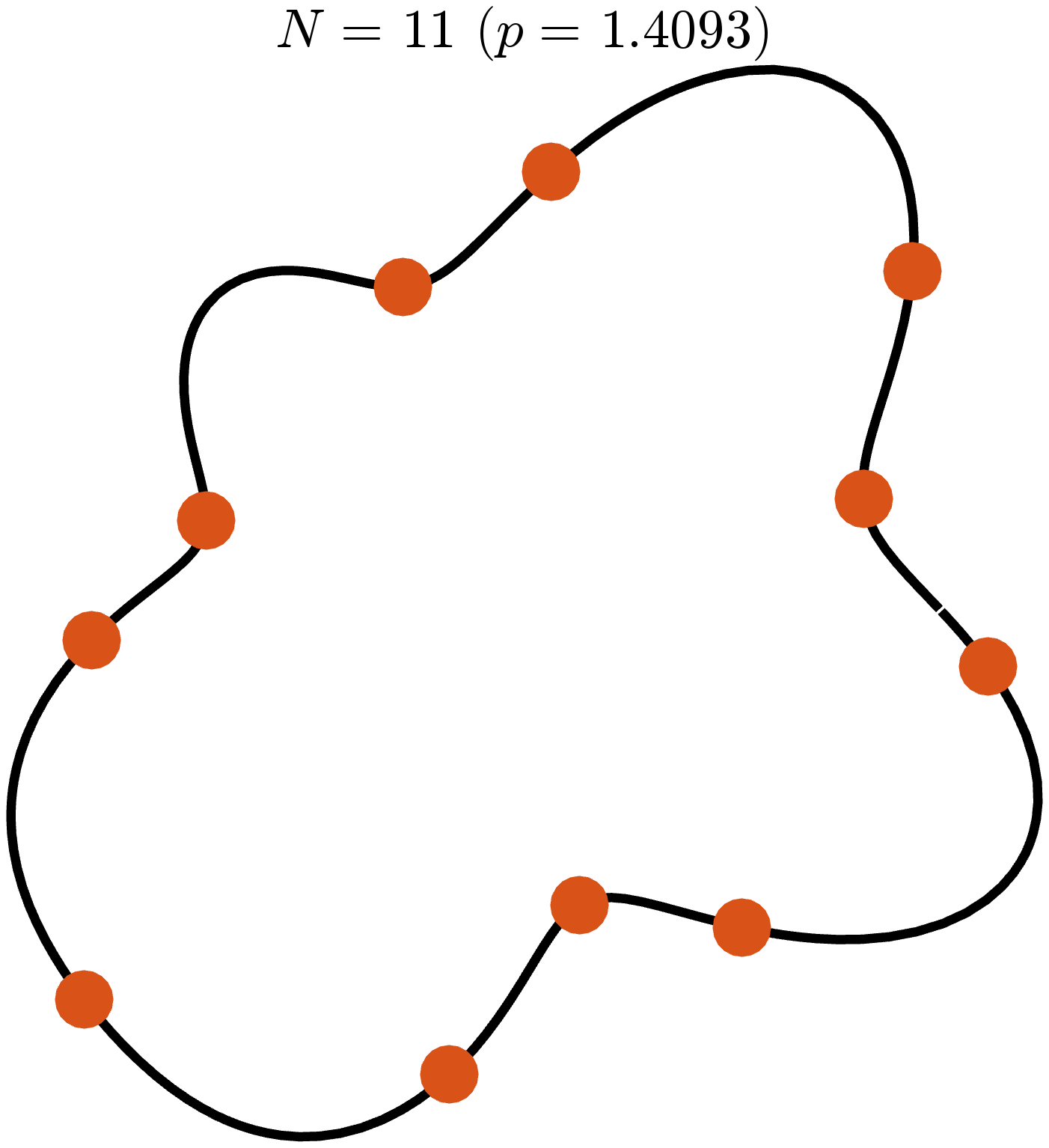}
    \includegraphics[width = 0.159\textwidth]{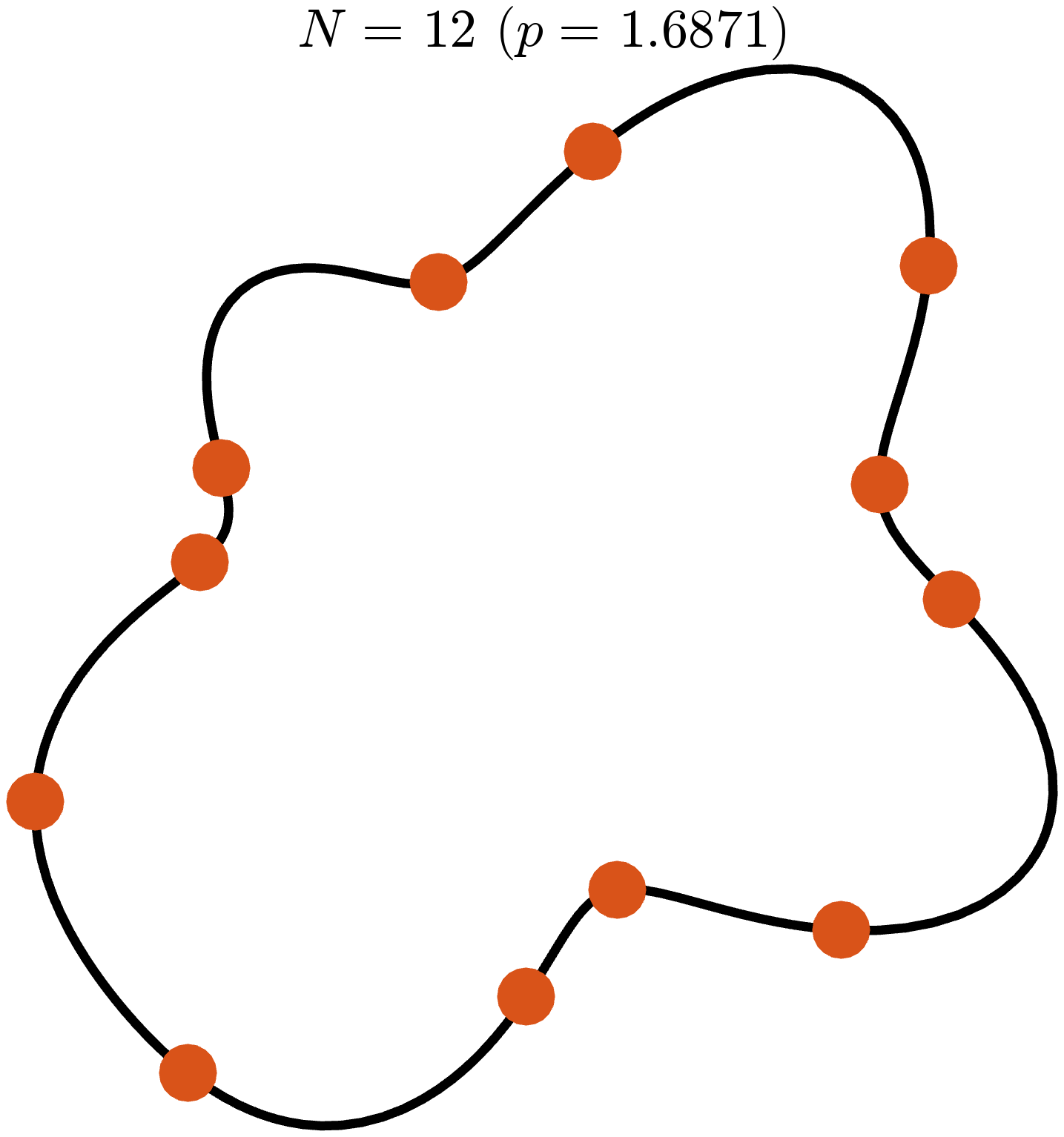}
    \caption{Minimizing locations of the discrete energy $p(\emph{\bx}_1,\ldots,\emph{\bx}_N)$ for boundary windows located on a random domain \eqref{eq:RandomDomain}. From top left to bottom right, the patterns are shown for $N=1$ to $N=12$ windows respectively.\label{fig:randomSurf}}
\end{figure}

\subsection{Escape through boundary windows.}\label{sec:Results_NEP}

\begin{figure}[htbp]
\centering
\subfigure[Ellipse $a=2$, $b=1$.]{\includegraphics[width=0.45\textwidth]{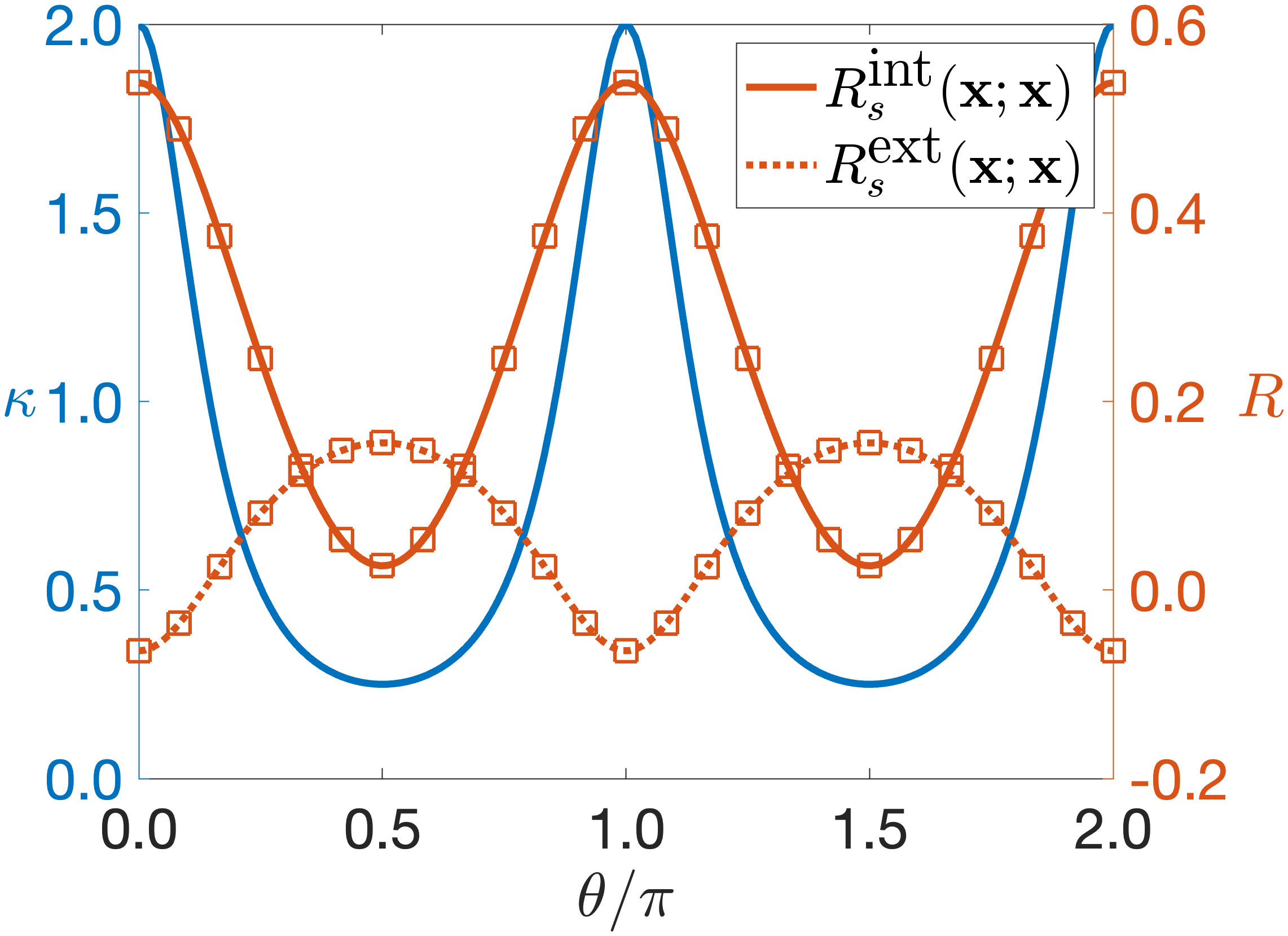}
\label{fig:curvature_vs_R_a}}
\quad
\subfigure[The random domain.]{\includegraphics[width=0.45\textwidth]{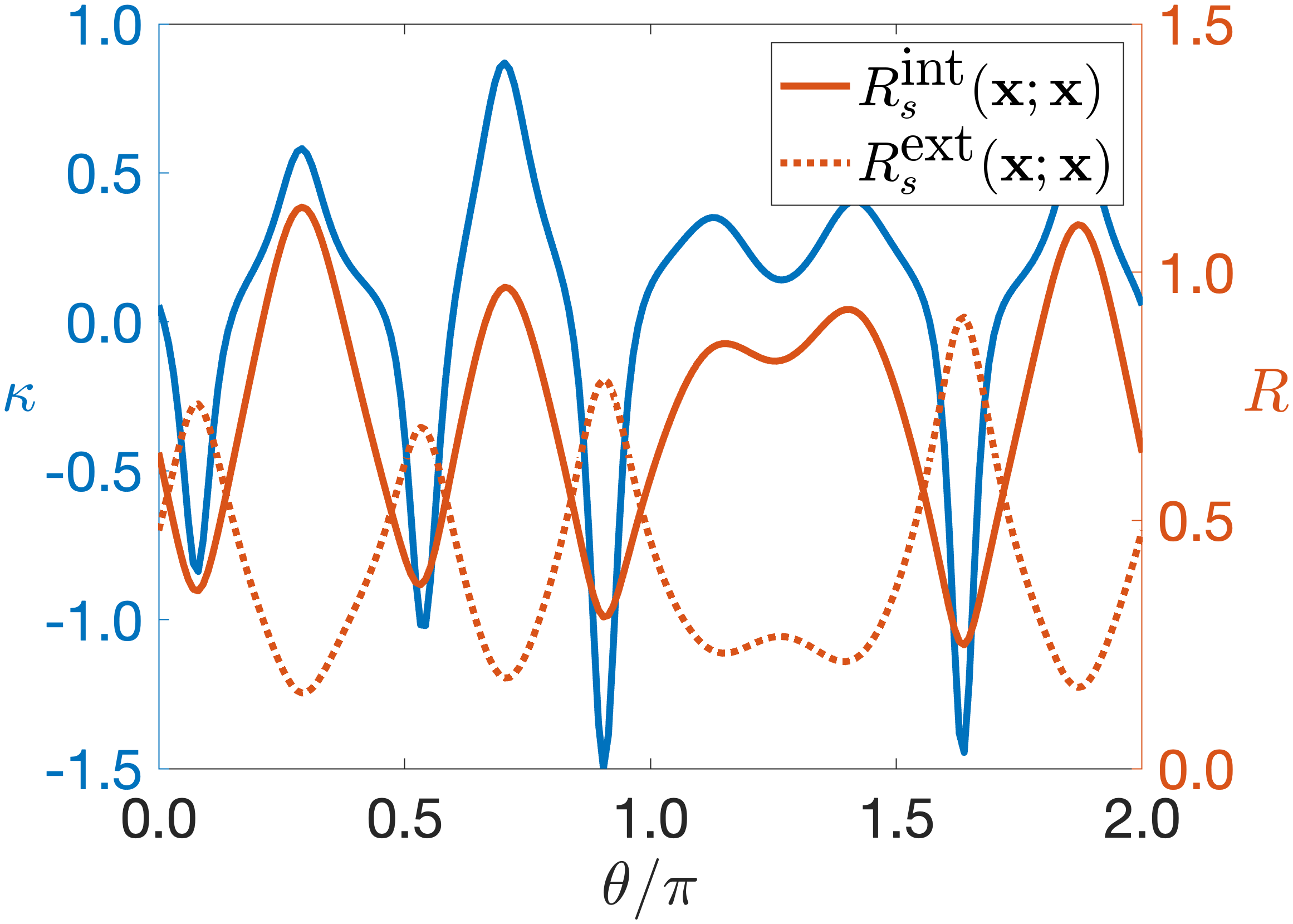}\label{fig:curvature_vs_R_b}}
    \caption{The values of $\Rint_{s}(\emph{\bx}(\theta))$, $\Rext_{s}(\emph{\bx}(\theta))$ and $\kappa(\theta)$ for the elliptical and random domain \eqref{eq:RandomDomain} for boundary parameterization $\partial\Omega = \{ \emph{\bx}(\theta) \ | \ \theta\in(0,2\pi)\}$. We remark that the left figure partially reproduces Fig.~7(a) of \cite{ward2010a} and displays results both from our numerical method (lines) and the closed form expressions (squares) given in \eqref{GreensEllipse_Surf}.\label{fig:curvature_vs_R}}
\end{figure}

It is natural to investigate whether the optimizing locations of the boundary windows are related in some way to the local geometric features of the boundary. In Fig.~\ref{fig:curvature_vs_R} we plot the regular part $\Rint_s(\bx)$ for $\bx(\theta)=(x_1(\theta),x_2(\theta))$ where $\theta\in(0,2\pi)$ is parameterization of the boundary. Together with the regular part, we plot the signed boundary curvature
\[
\kappa(\theta) = \frac{x_1' x_2'' - x_2' x_1''}{(x_1'^2 +x_2'^2)^{\frac{3}{2}} },
\]
as it varies along $\partial\Omega$. In the case of the elliptical domain, we observe in Fig.~\ref{fig:curvature_vs_R_a} a clear correspondence between the maxima and minima of the curvature and those of the regular part. This suggests that for a single boundary window on an elliptical domain, the placement should be at the boundary point of lowest curvature in order to minimize the GMFPT \eqref{eqn:tau_0} of exit from the domain.

In the case of the random domain \eqref{eq:RandomDomain}, we notice in Fig.~\ref{fig:curvature_vs_R_b} more nuance in the relationship between the maxima and minima of $\kappa(\theta)$ and $\Rint_s(\bx(\theta))$. First, a clear proximity between critical points of $\Rint_s(\bx(\theta))$ and critical points of $\kappa(\theta)$ is observed. Second, we notice that neither the global maximum nor global minimum of $R_s(\bx(\theta))$ corresponds to the global maximum or global minimum of $\kappa(\theta)$. This reflects the fact that the regular part $\Rint_s$ is a global quantity that depends both on local features at $\bx(\theta)\in\partial\Omega$ and the entire domain $\Omega$. These examples do, however, suggest a general principle that the critical points of $R_s(\bx(\theta))$ and $\kappa(\theta)$ are \emph{close} to one another. However, as seen in the carefully constructed counter example of \cite[$\S 4.1$]{LindsayWard2010}, the nature of these critical points may not correspond. It remains an open problem to fully describe the relationship between $R_s(\bx)$ and the local and global properties of $\partial\Omega$. 

As the number $N$ of boundary traps increases, the optimization problem becomes influenced by interaction between traps and in particular their pairwise distance through the logarithmic free space Green's function. Interestingly, for $N=2$ we see in the ellipse case $a= 2$ that the optimal locations are not the points $(0,\pm b)$ as would be suggested by a preference towards points of lowest boundary curvature. Instead, the optimal configuration moves away from regions of preferred curvature in order to increase the pairwise distance. However, we do observe in this elliptical example that points have a preference to cluster around these two points of most negative curvature. 

This trend is reinforced when we consider the optimal placement for the random domain (cf.~Fig.~\ref{fig:randomSurf}). In Fig.~\ref{fig:curvature_vs_R_b} we observe that $\Rint_s$ has four distinct local minima. The optimizing location for traps $N= 1,\ldots,4$ corresponds to a placement at each of these four points. As the number of placed points increases, we observe a pattern where clustering occurs around these points of minimizing curvature.

\subsection{Geometric shielding effects in source detection.}

Ratiometric sensing is a proposed mechanism for inferring the directionality of a signaling external sources through a comparison of receptor activity over the cellular surface. The underlying concept is that heightened receptor activity in a particular section of the cellular surface indicates the directionality of the source \cite{Lew2019,Ismael2016,Lakhani2017,Bumsoo2019}. In this example, we explore how cellular geometry can sharpen or modulate this mechanism. We consider the class of \emph{cassini} domains given parametrically for $\bx = (x_1,x_2)$ by
\begin{equation}\label{eq:cassini_domain}
    \big((x_1-a)^2 + x_2^2\big) \big((x_1+a)^2 + x_2^2\big) = b^4.
\end{equation}
This choice of domain, which varies over the two parameters $(a,b)$, allows us to consider a range of shapes which can reflect non-circularity observed in cellular shapes from ovals to rods like shapes and also dumbbell geometries. 

\begin{figure}[h!]
    \centering
    {\includegraphics[width=0.5\textwidth]{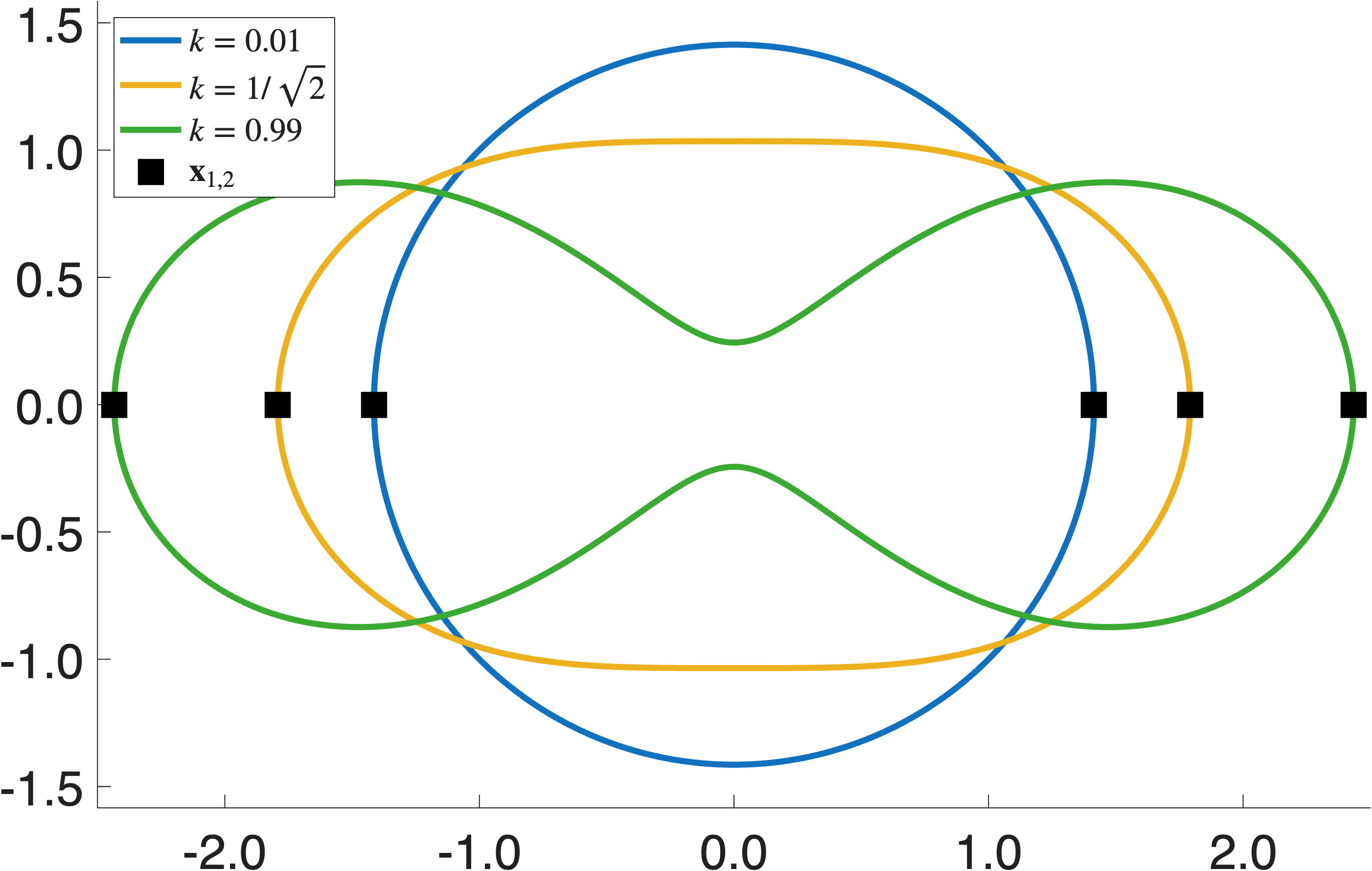}}
    \caption{The cassini family of domains \eqref{eq:cassini_domain} for fixed area $|\Omega|=\pi$. The parameter $k=\frac{a}{b}\in(0,1)$ describes oval geometries for $k\in(0,1/\sqrt{2})$, a \lq\lq pill\rq\rq\ shaped domain for $k=1/\sqrt{2}$ and dumbbell geometries as $k\to1$. There are two receptors of extent $\eps$ located at the rightmost ($\emph{\bx}_1$) and leftmost ($\emph{\bx}_2$) edges of the domain.\label{fig:cassin_domains}}
\end{figure}

If we fix the area $|\Omega|=b^2E(a^4/b^4)$ of these curves, we can reduce to a one-parameter family of curves for $k= a/b\in(0,1)$ where
\begin{equation}
  a = k b , \qquad b^2 = \frac{|\Omega|}{E(k^4)}.
\end{equation}
Here $E(m) = \int_0^{\frac{\pi}{2}}[1- m \sin^2 u]^{\frac12}du$ is the complete elliptic integral of the second kind. As $k\to0$, the curve becomes circular, while $k\to1$ represents a dumbbell shaped domain whose neck pinches at $k=1$. The case $k=1/\sqrt{2}$ represents a \lq\lq pill\rq\rq\ shaped geometry typical of bacteria (see Fig.~\ref{fig:cassin_domains}).

In order to explore how geometry can influence the mechanism of ratiometric sensing, we place a source at position $\bx=[R,0]$ and consider the \emph{differential splitting probabilities} defined (see \cite{BL2025}) as 
\begin{equation}\label{ex:splitting_cassini}
\Xi[\bx] = \phi^{\ast}_1(\bx) - \phi^{\ast}_2(\bx) \in(-1,1),
\end{equation}
where $\phi^{\ast}_1(\bx)$ is the splitting probability to the right receptor and $\phi^{\ast}_2(\bx)$ the splitting probability to the left receptor from $\bx\in\mathbb{R}^2\setminus\Omega$. The splitting probabilities are described in Sec.~\ref{sec:splitting} and a matched asymptotic solution in the small receptor limit in given by \eqref{eq:splitting_asy}. In this setting, values of $\Xi\approx 1$ are associated with a strong directional signal to the right, $\Xi\approx -1$ with a strong directional signal to the left and when $\Xi\approx 0$, there is little directional information to be inferred from the splitting probabilities.

\begin{figure}[htbp]
\centering
\subfigure[$\Xi$ against $k$ for $R=5$.]{\includegraphics[width=0.45\textwidth]{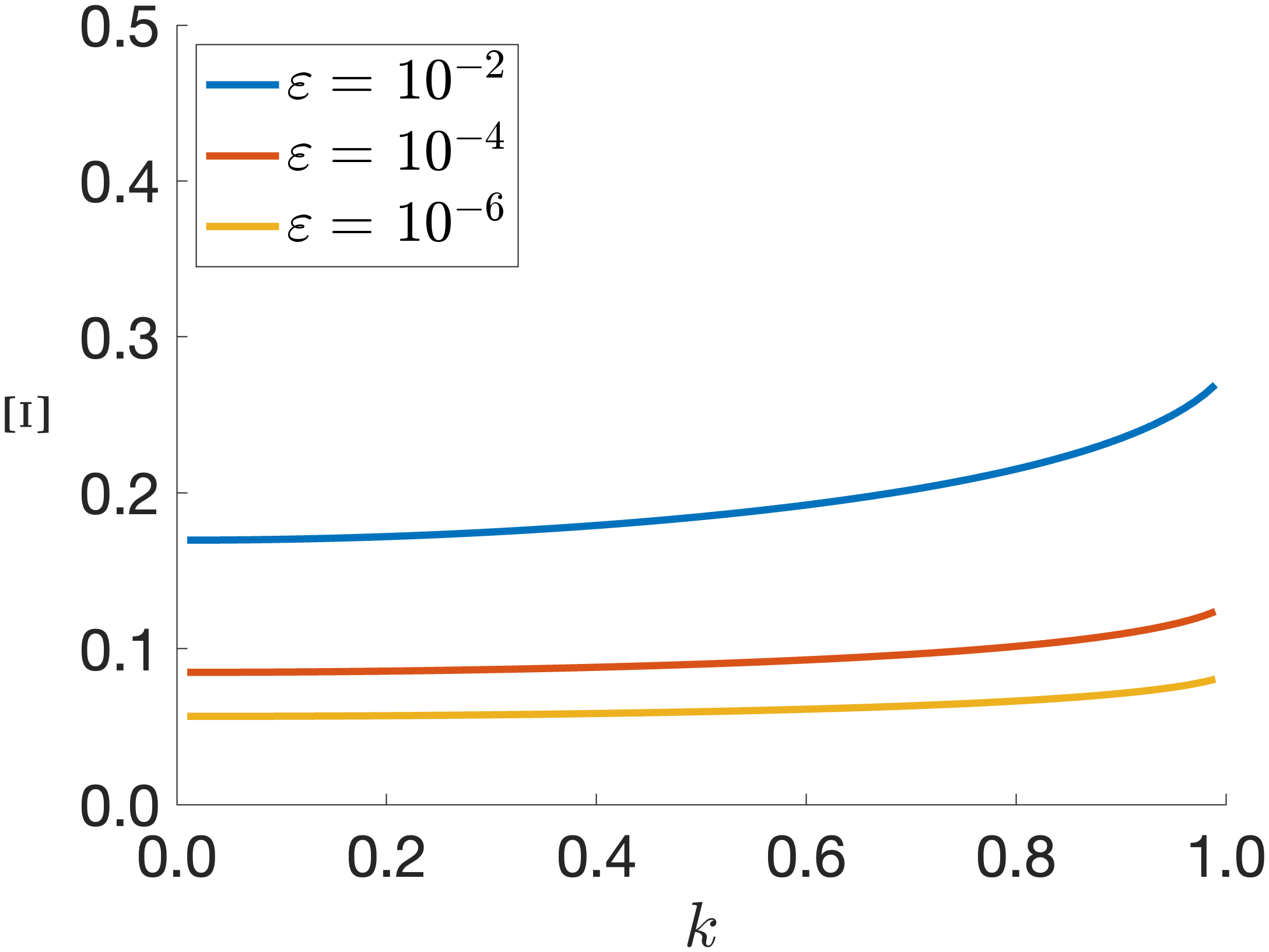}\label{fig:cassini_split_ex_a}}\qquad
\subfigure[$\Xi$ against $k$ for $\eps=10^{-4}$.]{\includegraphics[width=0.45\textwidth]{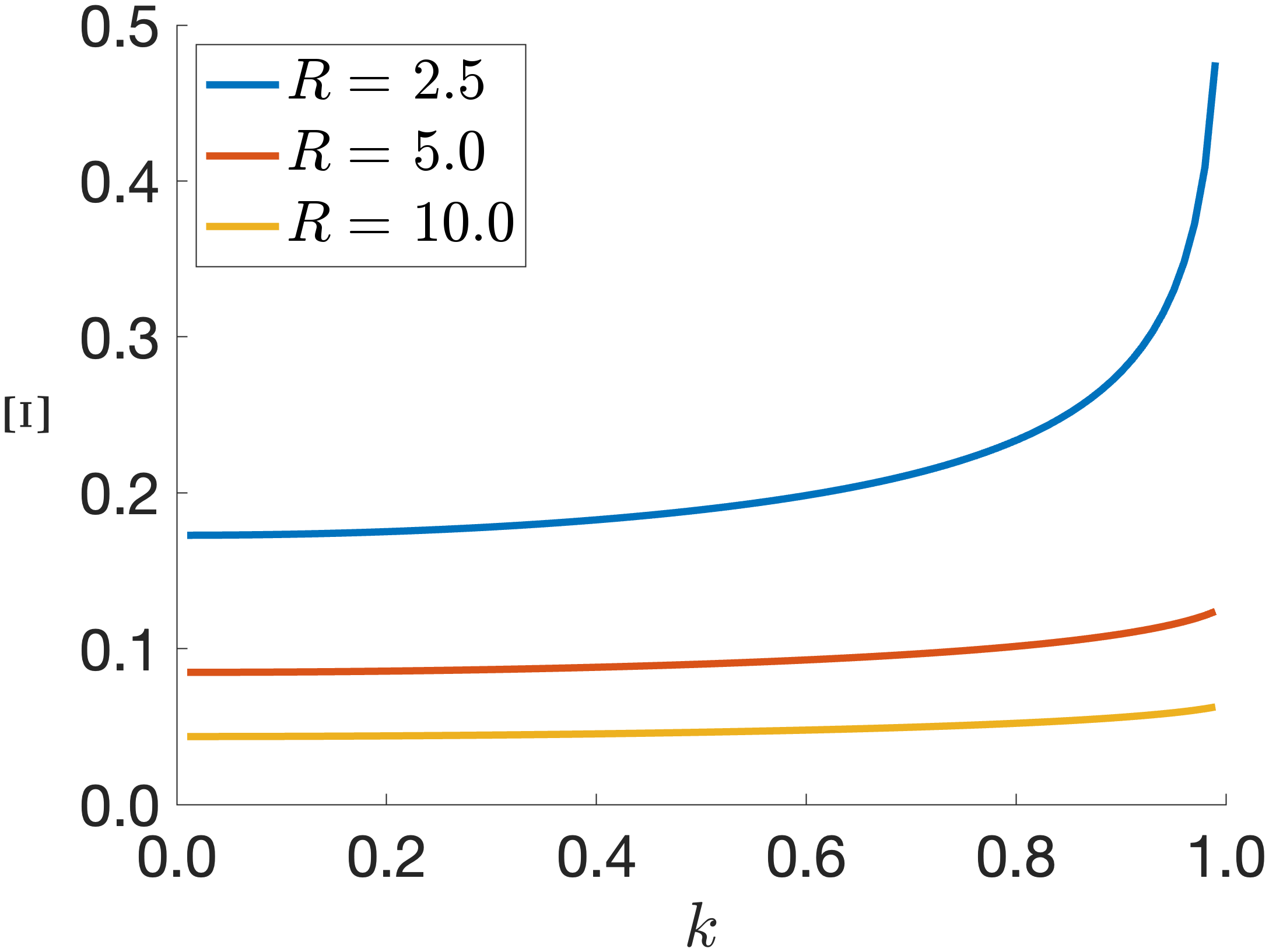}\label{fig:cassini_split_ex_b}}
    \caption{The relative splitting probabilities $\Xi$ defined in \eqref{ex:splitting_cassini} for the family of cassini domains \eqref{eq:cassini_domain} parameterized by $k\in(0,1)$. In panel (a), we show curves for fixed $R=5$ and various $\eps$ values. As $\eps\to0$, we see $\Xi\to0$ indicating that directional sensing through the splitting properties becomes infeasible in the small receptor limit. In panel (b), we show curves for fixed $\eps= 10^{-4}$ and various source distances $R$. \label{fig:cassini_split_ex}}
\end{figure}

Using our numerical methods, we build the linear system \eqref{eq:splitting_linear} for a range of geometries parameterized by $k\in(0,1)$. We place two receptors at positions $\bx_{1,2} = \pm(\sqrt{a^2+b^2},0)$ which corresponds to the rightmost ($\bx_1$) and leftmost $(\bx_2)$ points of the domain. We use the receptor extent $\eps$ to define the gauge parameter $\nu=-1/\log\eps$ (setting the logarithmic capacitance to unity). In Fig.~\ref{fig:cassini_split_ex_a} we hold the source distance fixed at $R=5$ and vary the scale parameter $\eps$. Recalling that a larger value of $\Xi$ is associated with stronger directional sensing, we observe that for each curve the dumbbell shape ($k\approx 1$) produces the largest value of $\Xi$. However, we additionally observe that $\Xi\to 0$ as $\eps\to0$ indicating that as  receptor extent shrinks, the directional signal also vanishes. This observation is corroborated by the result that (see \cite[Eqn.~(2.6)]{Lindsay2023} and \cite[Theorem 1.6]{lelievre2024spectralapproachnarrowescape})
\[
\lim_{\eps\to0}\phi_k^{*}(\bx) = \frac{1}{N}, \qquad k = 1,\ldots,N.
\]
This result states that, in the limit of vanishing receptor extent, the probability of reaching any particular receptor is uniform, independent of starting location, trap location, trap extent, and domain geometry. Consequently, splitting probabilities cannot be used for source inference in the vanishing receptor extent limit in 2D. The same result does not extend to 3D as shown in the sphere case \cite{Miles2020}.

\subsection{Optimization of trap orientations for internal elliptical absorbers.}

In a recent study \cite{CTL2025}, a high order expansion of the MFPT equation
\bsub\label{eq:MFPT_abs_main}
\begin{gather}
\label{eq:MFPT_abs_main_a}    D\Delta T +1 = 0, \qquad \bx\in\Omega;\qquad \partial_{\bn} T = 0, \qquad \bx\in\partial\Omega;\\[2pt]
\label{eq:MFPT_abs_main_b}    T = 0, \qquad \bx\in\partial\Omega_{\eps}\qquad j = 1,\ldots,N,
\end{gather}
\esub
was determined for a single elliptical all-absorbing trap centered at $\bx_1$. The trap has semi-major axis $\eps a$, semi-minor axis $\eps b$ and orientation $\phi$ with respect to the horizontal.  In the limit as $\eps\to0$, the two term GMFPT \eqref{eqn:GMFPT} was determined to be
\begin{equation}\label{eqn:intro_tau}
\tau = \frac{\tau_0}{D} + \frac{\eps^2}{D} \left[\frac{\pi ab}{|\Omega|}\tau_0 + \frac{a^2+b^2}{4} - \pi |\Omega| \frac{(a+b)^2}{2} (R_{x_1}^2 + R_{x_2}^2) + |\Omega| \frac{a^2-b^2}{4} \textbf{p} \cdot \begin{bmatrix}\cos 2\phi \\ \sin 2\phi \end{bmatrix} \right]
\end{equation}
Here the leading order term is 
\[
\tau_0 =  \frac{|\Omega|}{2\pi }\Big[ \frac{1}{\nu} + 2\pi R(\bx_1;\bx_1) \Big], \qquad \nu(\eps) = \frac{-1}{\log( \eps (a+b)/2)}.
\]
The vector $\textbf{p}$ is determined in terms of the derivatives of the regular part to be
\begin{equation}\label{eqn:intro_optimal_p}
    \textbf{p} = \begin{bmatrix} 
R_{x_1x_1}-R_{x_2x_2} - 2\pi (R_{x_1}^2 -R_{x_2}^2) \\[4pt]
2R_{x_1x_2} - 4\pi R_{x_1} R_{x_2}
\end{bmatrix}_{\bx=\bx_1}\, .
\end{equation}

\begin{figure}[htbp]
\centering
\subfigure[$a=1$, $b=1$.]{\includegraphics[height =1.2in]{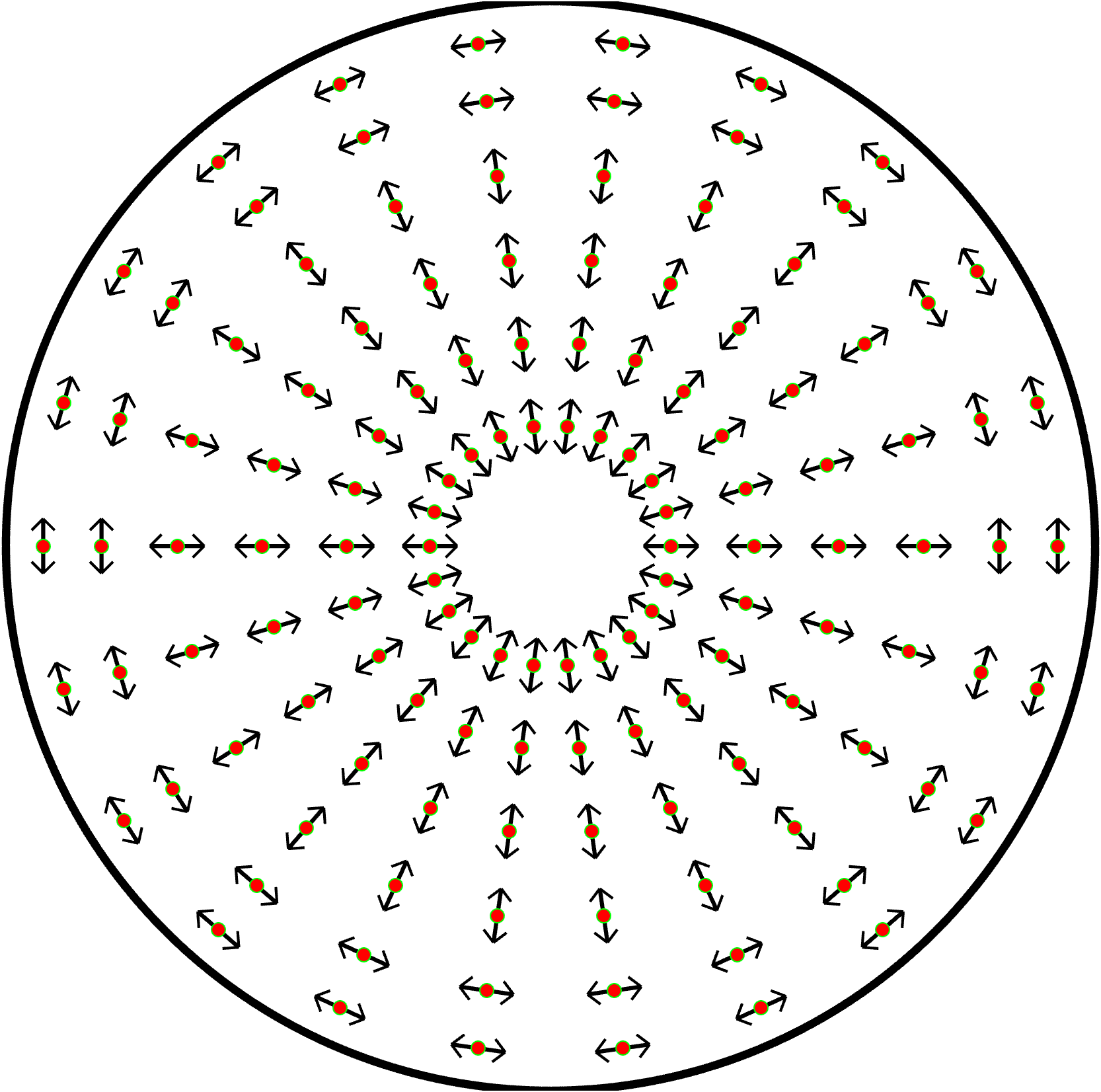}
\label{fig:EllipseDiskOrientation_a}}\quad
\subfigure[$a=1.1$, $b=1$.]{\includegraphics[height =1.2in]{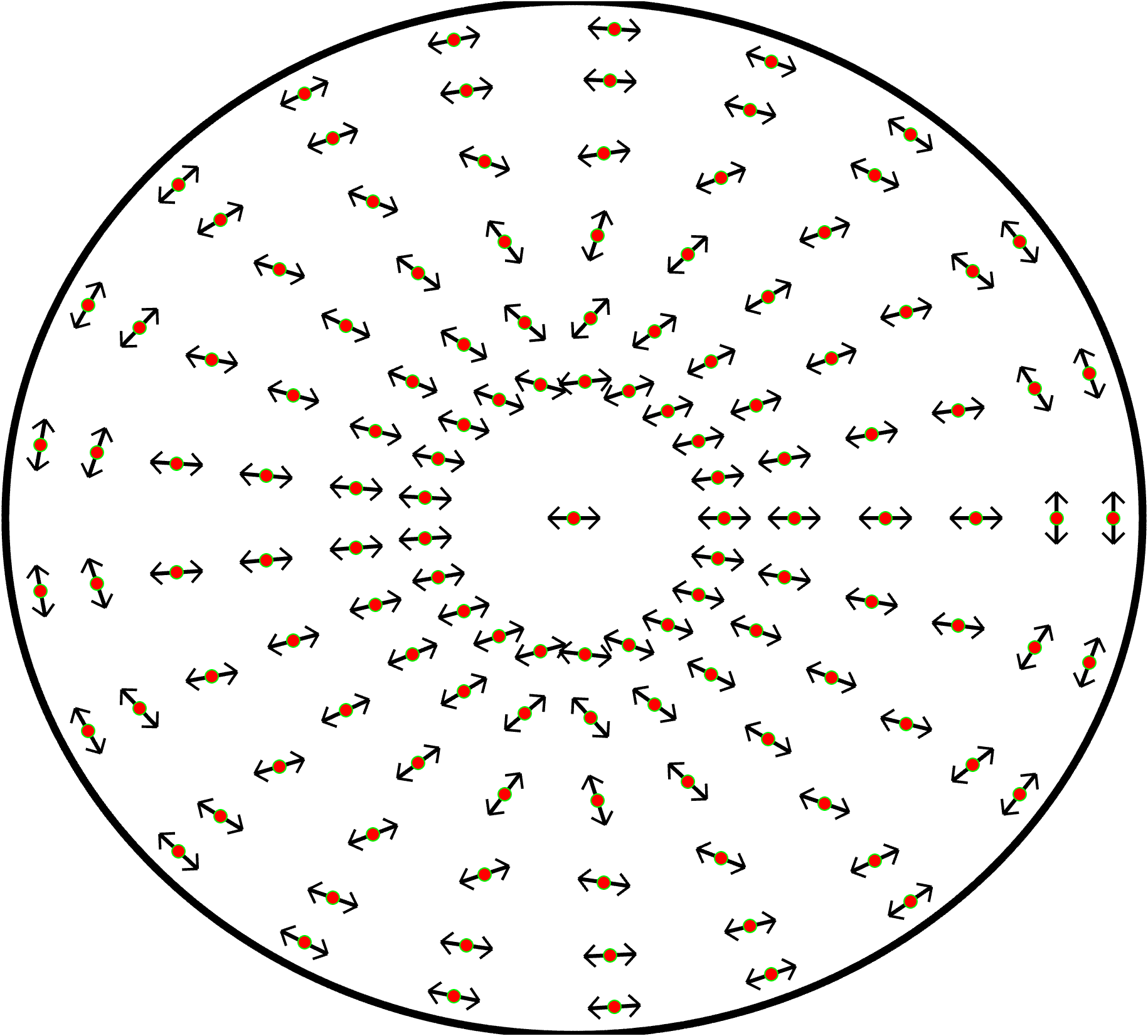}
\label{fig:EllipseDiskOrientation_b}}\quad
\subfigure[$a=1.5$, $b=1$.]{\includegraphics[height =1.2in]{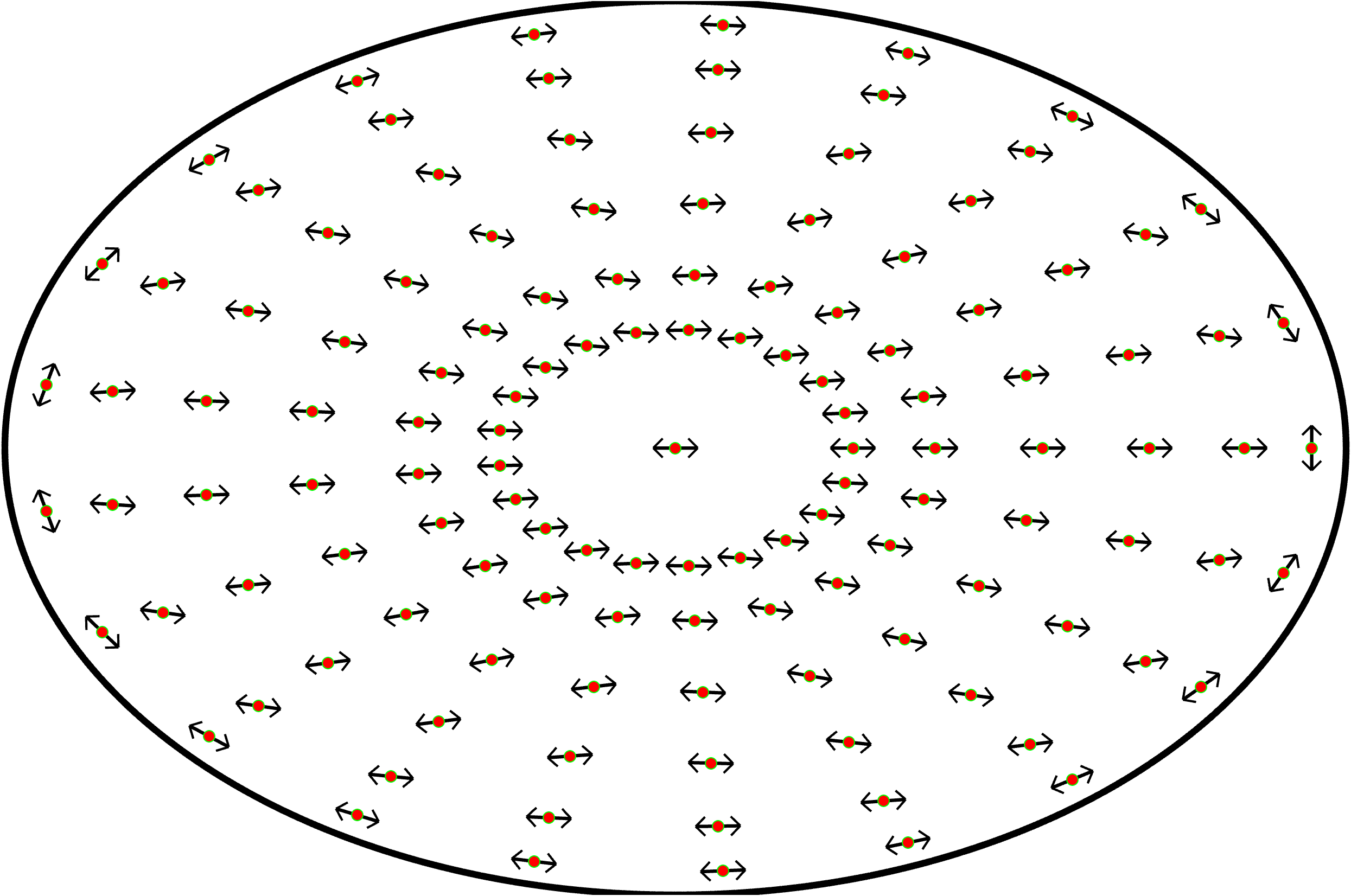}
\label{fig:EllipseDiskOrientation_c}}
\caption{The GMFPT minimizing direction for the semi-major axis of a small elliptical trap in a bounded domain. This figure is a reproduction of \cite[Fig.~6]{CTL2025}. In panel (a) the results for a disk show a bifurcation whereby the GMFPT is minimized when the trap is oriented in the radial direction for $0<|\emph{\bx}|<r_c$ and in the tangential direction for $r_c<|\emph{\bx}|<1$. In \cite{CTL2025} this critical value was determined to be $r_c= \sqrt{2-\sqrt{2}}$. In panels (b-c) we observe the smoothing out of this structure as the disk is perturbed into an ellipse for $a=1.1$ $b=1$ (panel (b)) $a=1.5$ $b=1$ (panel (c)). \label{fig:EllipseDiskOrientation}}
\end{figure}

The vector $\textbf{p}$ hence gives the direction along which the trap should be oriented to optimize the correction term $\tau_2$ of the GMFPT. We remark that to minimize the leading order term $\tau_0$, the trap should be located at points where $R_{x_1}=R_{x_2} = 0$ in which case $\textbf{p} = [R_{x_1x_1} -R_{x_2x_2}, 2R_{x_1x_2}]^T$. 

For a few domains, we solve for the regular part and directly calculate the first and second derivatives from the solution \eqref{eqn:Reg_Final}. First, we replicate in Fig.~\ref{fig:EllipseDiskOrientation} the results of the previous study \cite[Fig.~6]{CTL2025} on GMFPT minimizing trap orientations in disk and elliptical domains $\Omega$ for  which closed form solutions for the regular part are available. We confirm a bifurcation for the disk geometry whereby the GMFPT minimizing orientation for a trap centered at $\bx\in\Omega$ is in the radial direction for $0<|\bx|<r_c$ and in the tangential direction for $r_c<|\bx|<1$. The critical value was calculated \cite{CTL2025} to be $r_c = \sqrt{2-\sqrt{2}} \approx0.7654$. As the disk geometry is deformed into ellipses (cf.~Figs.~\ref{fig:EllipseDiskOrientation_b}-\ref{fig:EllipseDiskOrientation_c}), this sharp transition is gradually smoothed out.

\begin{figure}
\centering
\subfigure[Cassini domain for $k=0.99$.]{\includegraphics[height = 3cm]{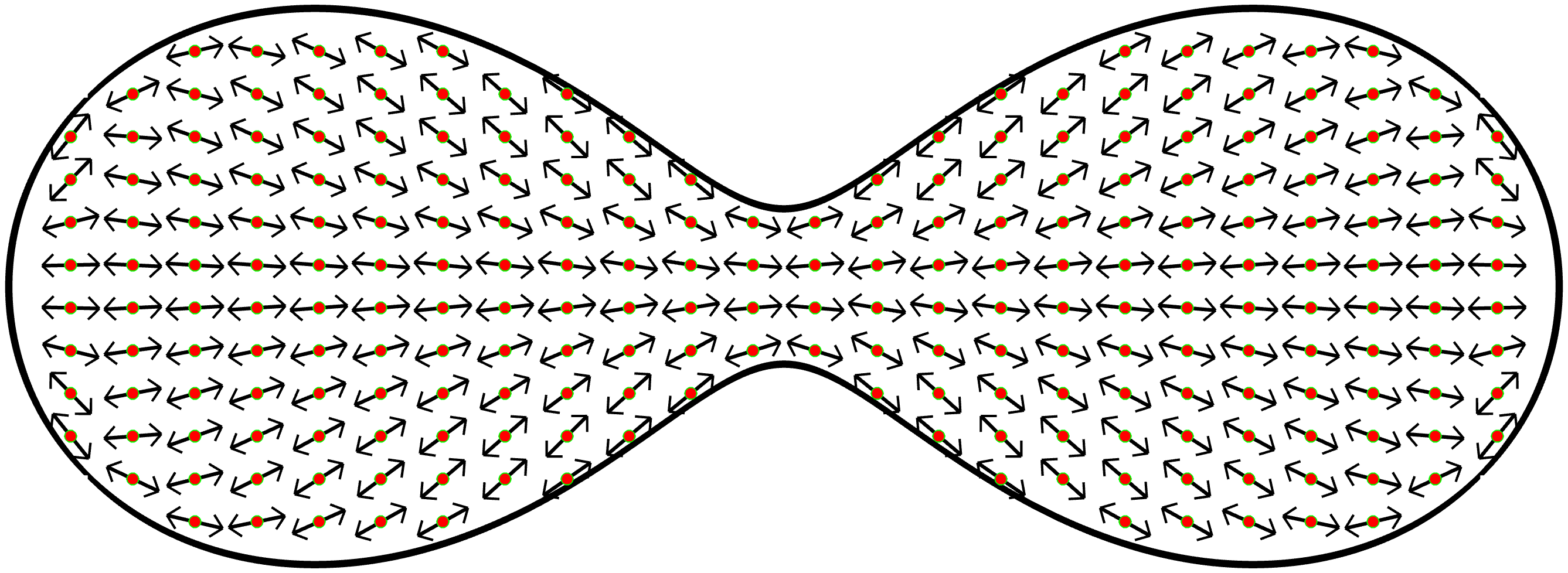}
\label{fig:orientation_examples_a}}\qquad
\subfigure[Star domain.]{\includegraphics[height = 3cm]{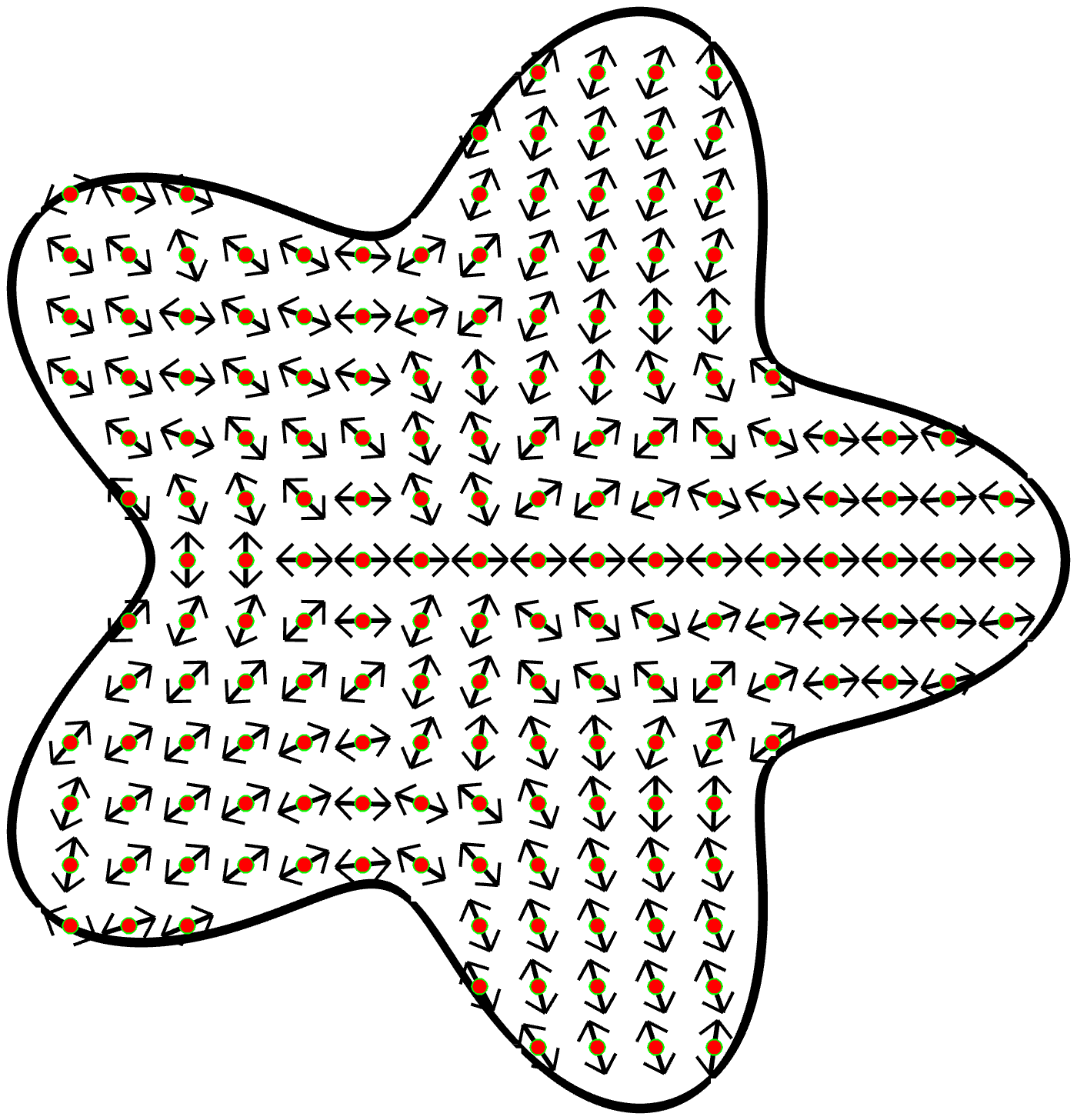}
\label{fig:orientation_examples_b}}\\[5pt]
\subfigure[Barbell domain.]{\includegraphics[height = 3cm]{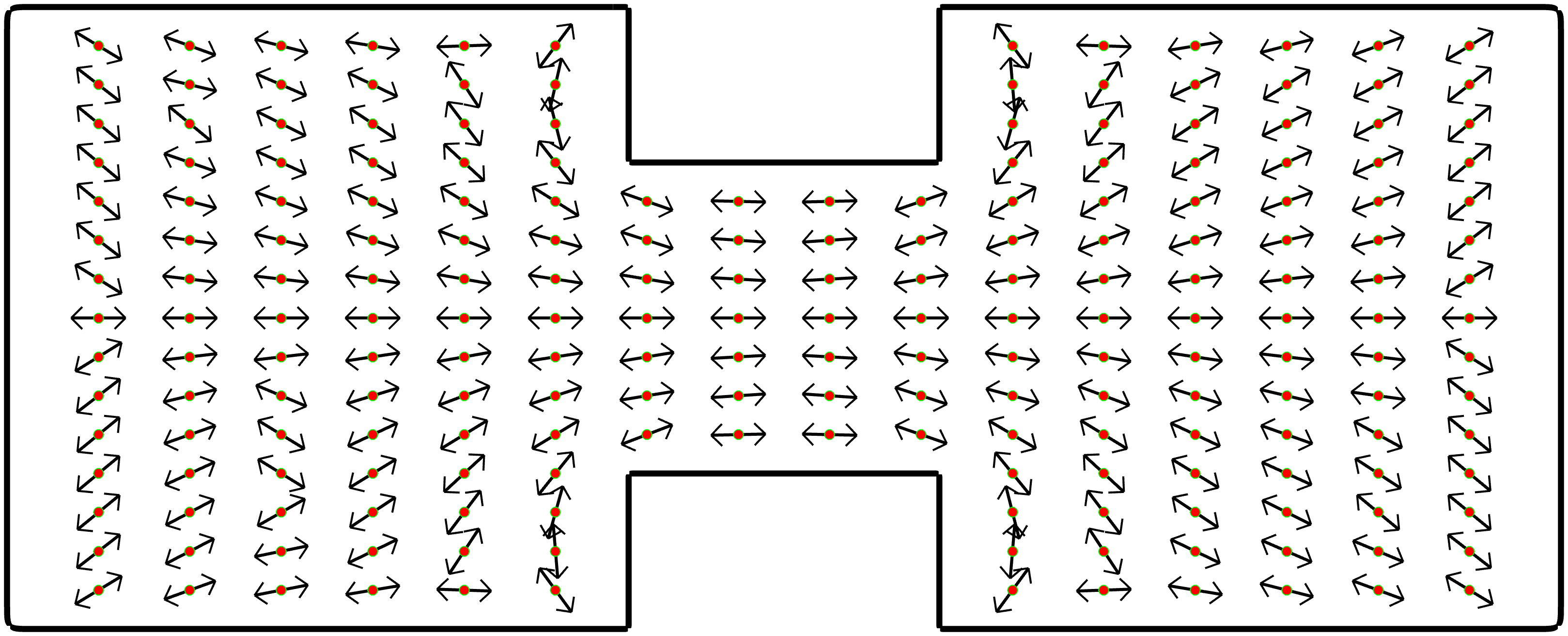}
\label{fig:orientation_examples_c}}\qquad
\subfigure[Random domain.]{\includegraphics[height = 3cm]{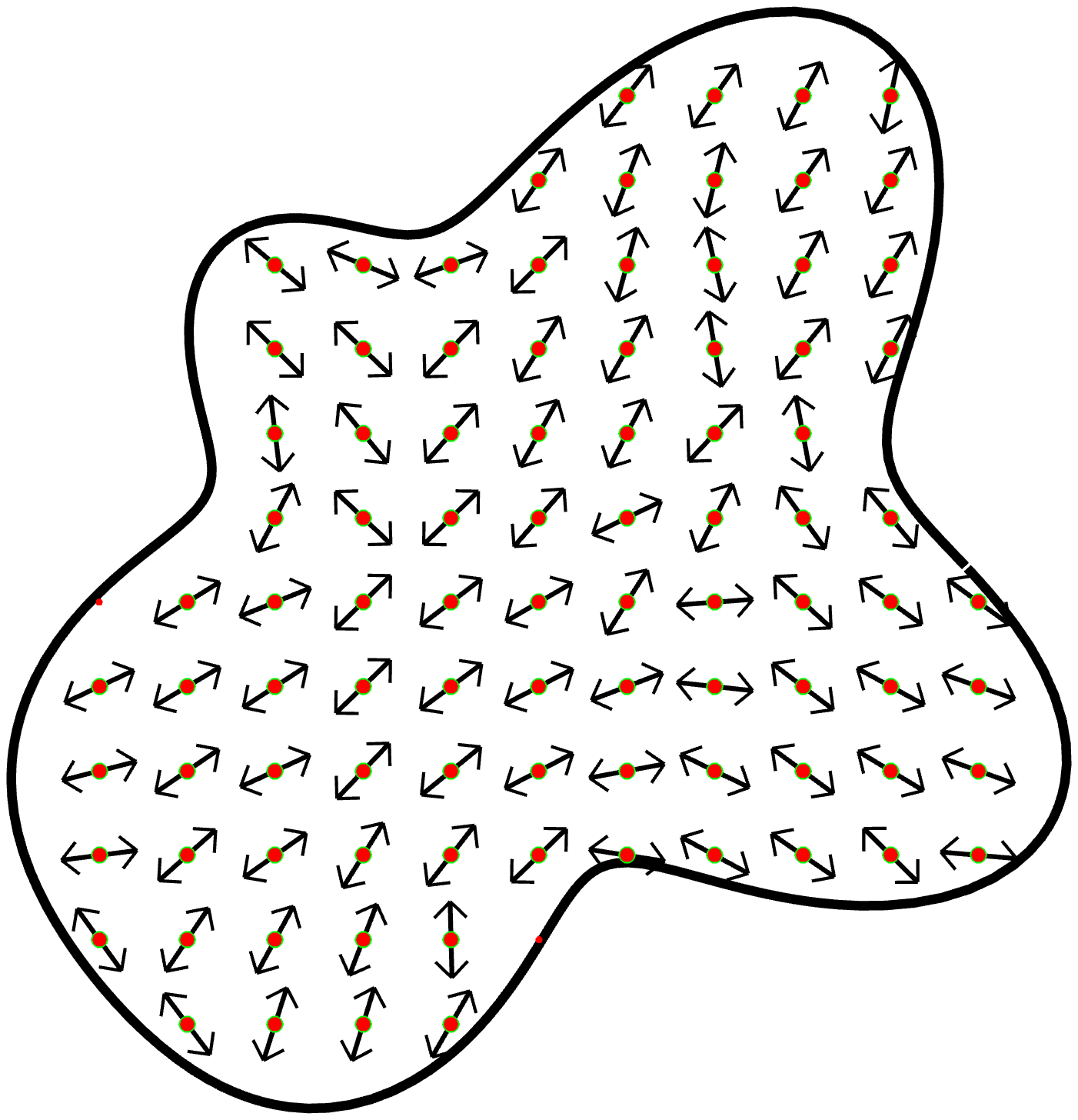}
\label{fig:orientation_examples_d}}
\caption{The GMFPT minimizing orientation for a single elliptical trap as predicted by \eqref{eqn:intro_optimal_p} for four domains. Panel (a): The Cassini oval for $k=0.99$.  Panel (b): A star domain with five arms. Panel (c): A barbell domain formed by two squares connected by a small rectangular neck. Panel (d): The random domain seen in previous examples.
\label{fig:orientation_examples}}
\end{figure}

In Fig.~\ref{fig:orientation_examples} we display optimal results for a few general geometries, including dumbbell (Fig.~\ref{fig:orientation_examples_a}) a star-shaped domain (Fig.~\ref{fig:orientation_examples_b}), a barbell domain (Fig.~\ref{fig:orientation_examples_c}) and the previously considered random domain (Fig.~\ref{fig:orientation_examples_d}). We observe a few general trends. First, the traps have a preference to be aligned in the longitudinal direction of the domain. Second, traps centered close to $\partial\Omega$ have a preference to be aligned parallel to the boundary.

\section{Discussion}\label{sec:Discussion}

In this work we have presented and validated a high order numerical method for the planar Green's function, a central quantity in narrow capture and escape problems and applications to chemical signaling. Previously available solutions were an exact formula in the case of a disk domain and rapidly convergent series for rectangles and ellipses. By interfacing with standard optimization routines, we have shown GMFPT minimizing arrangements of traps, either placed interior to the domain or on the boundary as escape windows. 

There are several avenues for extension of this work which warrant further study. One is extension to other operators such as the modified Helmholtz operator. For the interior problem in $\bx\in\Omega$, this equation satisfies
\bsub\label{eqn:MHH}
\begin{gather}
    \Delta G_h - s G_h= - \delta(\bx-\by), \qquad \bx\in\Omega;\\
    \partial_{\bn}G_h = 0, \quad \bx\in\partial\Omega, \qquad G_h (\bx;\by,s) = -\frac{1}{2\pi}\log|\bx-\by| + R_h(\bx;\by,s).
\end{gather}
\esub
The equation \eqref{eqn:MHH} arises when the Laplace Transform is applied to the heat equation. This approach is one route \cite{Lindsay2023,BL2025,CL2022,CHERRY2025,LindsayTzou2016,lacroix2025lightningmethodheatequation} to reproduce the full time-dependent statistics \cite{Grebenkov_2019} of arrivals at absorbing sites, i.e.~determining not just the mean statistics determined by \eqref{eq:MFPT_main}. These quantities are essential to determine extreme statistics for first passage time problems \cite{Linn_2022,Lawley2020,Morgan2023}.

In our presentation, we have reproduced existing optimization results for the minimizing locations of traps in a disk (cf.~Fig.~\ref{fig:disk_optimal}) and in an ellipse (cf.~Fig.~\ref{fig:EllipseN=4}). However, we expect that more sophisticated global routines may be able to identify lower energy configurations of \eqref{eq_Discrete}. Hence a priority for future studies is an integration of this  method with more sophisticated optimization strategies \cite{Gilbert_Cheviakov_2023,Ridgway2019,RIDGWAY2018}.

In future work, we will also address numerical algorithms for the corresponding three dimensional Green's function. In the case where $N$ small traps are located inside a bounded three dimensional region $\Omega$, similar asymptotic analysis have been performed in the limit as $\eps\to0$ \cite{CHEVIAKOV2011,BL2018}. In the case that each trap has a different size or shape, the GMFPT is given (see Principal Result 3.1 of \cite{CHEVIAKOV2011}) by
\begin{equation}
    \tau_0 = \frac{|\Omega|}{D}\left[\frac{1}{4\pi \eps N \bar{c} } + \frac{1}{N^2} \frac{\bc^T \mathcal{G} \bc}{\bar{c}^2} \right], \qquad \mathcal{G}_{i,j} = \left\{\begin{array}{rl}
        \Rint_b (\bx_i;\bx_i), & i = j;\\[2pt]
        \Gint_b (\bx_i;\bx_j), & i \neq j.
    \end{array} \right.
\end{equation}
The entries of the vector $\bc = (c_1,\ldots,c_N)$ are the capacitance of each trap and $\bar{c} = N^{-1}\sum_{j=1}^N c_j$. The capacitance of each trap is determined by their shape through the solution of the electrostatic problem
\bsub
\begin{gather}
    \Delta v_d = 0, \quad \bx  \in \mathbb{R}^3\setminus\Omega_k; \qquad v_d = 0, \quad \bx\in\partial\Omega_k;\\
    v_d  = \frac{c_d}{|\bx|} + \mathcal{O}(|\bx|^{-2}) \quad \mbox{as} \quad |\bx|\to\infty.
\end{gather}
\esub
The corresponding Green's function for the interior bulk problem is given by
\bsub\label{eqn_neumG3D}
\begin{gather}
 \label{eqn_neumG3D_a} \Delta \Gint_b = \frac{1}{|\Omega|}- \delta(\bx-\by), \quad \bx\in\Omega;\qquad  \partial_{\bn} \Gint_b \equiv  \nabla \Gint_b\cdot \textbf{n} = 0,\quad \bx\in\partial\Omega;\\
 \label{eqn_neumG3D_b} \int_{\Omega} \Gint_b(\bx;\by)d\bx = 0; \qquad
 \Gint_b(\bx;\by) \to \frac{1}{4\pi|\bx-\by|} + \Rint_b(\by)\,, \quad \mbox{as} \quad \bx\to\by.
\end{gather}
\esub
The singularity behavior and average constraint in \eqref{eqn_neumG3D_b} for the interior problem are relatively easy to accommodate using the strategy discussed in Sec.~\ref{sec:NumericalMethods}. The scenario of a surface source arises in many contexts \cite{10666768,BL2018,BL2025,Miles2020,BERNOFF2023,LBW2017,NURSULTANOV2021,CHEVIAKOV2011,Bressloff2023_3D,CheviakovWard2010}. In this case additional care is needed when dealing with the form of the surface behavior in $\Gint_s(\bx;\by)$ as $\bx\to\by$, in particular, there are lower order singularities that appear in the local behavior beyond the leading behavior $\Gint_s(\bx;\by) \sim (2\pi|\bx-\by|)^{-1}$ (see \cite{NURSULTANOV2021,Pop90,Silbergleit2003}). In future work, we will develop a numerical method that incorporates this local structure of the 3D surface Green's function.  

\section*{Acknowledgments} A.E.L. acknowledges support from the NSF under award DMS2052636. J.G.H. was supported in part by a Sloan Research Fellowship. This research was supported in part by grants from the NSF (DMS-2235451) and Simons Foundation (MPS-NITMB-00005320) to the NSF-Simons National Institute for Theory and Mathematics in Biology (NITMB).

\bibliographystyle{siam}
\bibliography{ref}

\end{document}